\documentclass[11pt]{article}
\usepackage{amssymb}
\usepackage{amsmath}
\usepackage{amsthm}
\newcommand{\cA}{{\cal A}}
\newcommand{\cB}{{\cal B}}

\newcommand{\cH}{{\cal H}}
\newcommand{\cE}{{\cal E}}

\newcommand{\cO}{{\cal O}}
\newcommand{\cL}{{\cal L}}

\newcommand{\cN}{{\cal N}}
\newcommand{\cF}{{\cal F}}

\newcommand{\cS}{{\cal S}}
\newcommand{\cT}{{\cal T}}
\newcommand{\cU}{{\cal U}}
\newcommand{\cV}{{\cal V}}

\newcommand{\cX}{{\cal X}}
\newcommand{\cY}{{\cal Y}}
\newcommand{\cZ}{{\cal Z}}

\renewcommand{\AA}{{\mathbb A}}

\newcommand{\ZZ}{{\mathbb Z}}

\newcommand{\PP}{{\mathbb P}}
\newcommand{\OO}{{\mathbb O}}

\newcommand{\gm}{\mathfrak{m}}    %gotic

\newcommand{\gp}{\mathfrak{p}}

\newcommand{\gt}{\mathfrak{t}}
\newcommand{\gr}{\mathfrak{r}}
\newcommand{\gs}{\mathfrak{s}}
\newcommand{\so}{\mathfrak{so}}
\newcommand{\gl}{\mathfrak{gl}}

\newcommand{\on}{\operatorname}

\newcommand{\Fl}{{\cal F}l}

\newcommand{\Rep}{{\on{Rep}}}

\newcommand{\Qlb}{\mathbb{\bar Q}_\ell}
\newcommand{\Gm}{\mathbb{G}_m}

\newcommand{\A}{\mathbb{A}}

\newcommand{\toup}[1]{\stackrel{#1}{\to}}
\newcommand{\hook}[1]{\stackrel{#1}{\hookrightarrow}}

\newcommand{\getsup}[1]{\stackrel{#1}{\gets}}
\newcommand{\Sp}{\on{\mathbb{S}p}}
\newcommand{\LieSp}{\on{\mathfrak{sp}}}

\newcommand{\IC}{\on{IC}}

\newcommand{\Hom}{\on{Hom}}

\newcommand{\End}{\on{End}}

\newcommand{\Sym}{\on{Sym}}
\newcommand{\SO}{\on{S\mathbb{O}}}
\newcommand{\Ker}{\on{Ker}}

\newcommand{\Aut}{\on{Aut}}

\newcommand{\RG}{\on{R\Gamma}}

\newcommand{\Pic}{\on{Pic}}

\newcommand{\Bun}{\on{Bun}}

\newcommand{\Bunt}{\on{\widetilde\Bun}}

\newcommand{\Spec}{\on{Spec}}

\newcommand{\supp}{\on{supp}}

\newcommand{\HOM}{{{\cal H}om}}

\newcommand{\Gr}{\on{Gr}}
\newcommand{\Grb}{\overline{\Gr}}

\newcommand{\GL}{\on{GL}}

\newcommand{\pr}{\on{pr}}
\newcommand{\id}{\on{id}}
\newcommand{\tr}{\on{tr}}
\newcommand{\QED}{$\square$} 
\newcommand{\Fq}{\mathbb{F}_q}  
\newcommand{\Fp}{\mathbb{F}_p}  % what for??????    
\newcommand{\iso}{{\widetilde\to}}

\newcommand{\comp}{\circ}
\newcommand{\Four}{\on{Four}}
\renewcommand{\H}{{\on{H}}}   %cohomologies
\newcommand{\DD}{\mathbb{D}}  %for duality
\newcommand{\D}{\on{D}}       %for derived categories     
\newcommand{\wt}{\widetilde}
\newcommand{\ov}[1]{\overline{#1}}
\newcommand{\select}[1]{{\it{#1}}}

\newcommand{\und}[1]{\underline{#1}}

\renewcommand{\P}{{\on{P}}}

\newcommand{\<}{\langle}
\renewcommand{\>}{\rangle}

\newcommand{\ev}{\mathit{ev}}

\newcommand{\Conv}{\on{Conv}}
\newcommand{\Loc}{\on{Loc}}
\newcommand{\Lie}{\on{Lie}}
\newcommand{\Sph}{\on{Sph}}
\newcommand{\Res}{\on{Res}}

\newcommand{\ttimes}{\tilde\times}

\newcommand{\tSp}{\wt\Sp}
\newcommand{\Funct}{\on{Funct}}
\newcommand{\Isom}{\on{Isom}}

\renewcommand{\Im}{\on{Im}}
\newcommand{\codim}{\on{codim}}
\newcommand{\Cov}{\on{Cov}}
\newcommand{\sign}{\on{sign}}
\newcommand{\SL}{\on{SL}}
\newcommand{\St}{\on{St}}
\newcommand{\tboxtimes}{\,\tilde\boxtimes\,}

\newcommand{\Vect}{\on{Vect}}
\newcommand{\Cl}{\on{Cl}}
\newcommand{\glob}{\on{glob}}
\newcommand{\Ind}{\on{Ind}}
\newcommand{\wh}{\widehat}
\newcommand{\ST}{\on{ST}}

\newtheorem{Lm}{Lemma}
\newtheorem{Th}{Theorem}
\newtheorem{Pp}{Proposition}
\newtheorem{Cor}{Corollary}

\theoremstyle{remark}
\newtheorem{Rem}{Remark}

\theoremstyle{definition}
\newtheorem{Def}{Definition}

\newenvironment{Prf}{\par\noindent {\it Proof }}{\QED}
\newcommand{\Step}[1]{\par\noindent{\bf Step {#1}}.}

\begin{document}

\author{Sergey Lysenko}
\title{Moduli of metaplectic bundles on curves and Theta-sheaves}
\date{}
\maketitle
\begin{abstract}
\noindent{\scshape Abstract}\hskip 0.8 em 
We give a geometric interpretation of the Weil representation of the metaplectic group, placing it in the framework of the geometric Langlands program. 

  For a smooth projective curve $X$ we introduce an algebraic stack $\Bunt_G$ of metaplectic bundles on $X$. It also has a local version $\wt\Gr_G$, which is a gerbe over the affine grassmanian of $G$. We define a categorical version of the (nonramified) Hecke algebra of the metaplectic group. This is a category $\Sph(\wt\Gr_G)$ of certain perverse sheaves on $\wt\Gr_G$, which act on $\Bunt_G$ by Hecke operators. A version of the Satake equivalence is proved describing $\Sph(\wt\Gr_G)$ as a tensor category. Further, we construct a perverse sheaf on $\Bunt_G$ corresponding to the Weil representation and show that it is a Hecke eigen-sheaf with respect to $\Sph(\wt\Gr_G)$.   
\end{abstract} 

\medskip

{\centerline{\scshape 1. Introduction}}

\bigskip\noindent
1.1  Historically $\theta$-series (such as, in one variable, $\sum q^{n^2}$)
have been one of the major methods of constructing automorphic forms. A representation-theoretic appoach to the theory of $\theta$-series, as discoved by A. Weil \cite{W} and extended by R. Howe \cite{H}, is based on the oscillator representation of the metaplectic group (cf. \cite{Pr} for a recent survey). In this paper we propose a geometric interpretation of this representation (in the nonramified case) placing it in the framework of the geometric Langlands program. 

 Let $k=\Fq$ be a finite field with $q$ odd. Set $K=k((t))$ and $\cO=k[[t]]$. Let $\Omega$ 
denote the completed module of relative differentials of $\cO$ over $k$. Let $M$ be a free $\cO$-module of rank $2n$ given with a nondegenerate
symplectic form $\wedge^2 M\to\Omega$. It is known that the continuous $\H^2(\Sp(M)(K), \{\pm 1\})\,\iso\, \ZZ/2\ZZ$ 
(\cite{Moo}, 10.4). As $\Sp(M)(K)$ is a perfect group, 
the corresponding metaplectic extension 
\begin{equation}
\label{meta_extension_small}
 1\to \{\pm 1\} \toup{i} \wh\Sp(M)(K)\to \Sp(M)(K)\to 1
 \end{equation}
is unique up to unique isomorphism. It can be constructed in two essentially different ways. 

   Recall the classical construction of A. Weil (\cite{W}). The Heisenberg group is $H(M)=M\oplus\Omega$ with operation 
$$
(m_1,\omega_1)(m_2,\omega_2)=(m_1+m_2, \, \omega_1+\omega_2+\frac{1}{2}\<m_1,m_2\>)
$$  
Fix a prime $\ell$ that does not divide $q$. Let $\psi: k\to \Qlb^*$ be a nontrivial additive character. Let $\chi:\Omega(K)\to \Qlb$ be given by $\chi(\omega)=\psi(\Res\omega)$. 
By the Stone and Von Neumann theorem (\cite{MVW}), there is a unique (up to isomorphism) smooth irreducible representation $(\rho, \cS_{\psi})$ of $H(M)(K)$ over $\Qlb$ with central character $\chi$. The group $\Sp(M)$ acts on $H(M)$ by group automorphisms
$
(m,\omega)\toup{g}(gm,\omega)
$
This gives rise to the group
\begin{multline*}
\;\;\;\tSp(M)(K)=\{(g,M[g])\mid g\in\Sp(M)(K), M[g]\in\Aut \cS_{\psi}\\
\!\!\!\!\rho(gm,\omega)\comp M[g]=M[g]\comp \rho(m,\omega)\;\; \mbox{for}\; (m,\omega)\in H(M)(K)\}\;\;\;
\end{multline*}
The group $\tSp(M)(K)$ is an extension of $\Sp(M)(K)$ by $\Qlb^*$. Its commutator subgroup is an extension of $\Sp(M)(K)$ by $\{\pm 1\}\hook{}\Qlb^*$,
uniquely isomorphic to (\ref{meta_extension_small}). 

  Another way is via Kac-Moody groups. Namely, view $\Sp(M)(K)$ as an ind-scheme over $k$. Let 
\begin{equation}
\label{Kac-Moody_extension_intro}
1\to \Gm\to \ov{\Sp}(M)(K)\to \Sp(M)(K)\to 1
\end{equation}
denote the canonical extension, here $\ov{\Sp}(M)(K)$ is an ind-scheme over $k$ (cf. \cite{F}). Passing to $k$-points we get an extension of abstract groups 
$1\to k^*\to \ov{\Sp}(M)(K)\to \Sp(M)(K)\to 1$. Then (\ref{meta_extension_small}) is the push-forward of this extension under $k^*\to k^*/(k^*)^2$.   
  
  The second construction underlies one of our main results, the tannakian description of the Langlands dual to the metaplectic group. Namely, the canonical splitting of (\ref{Kac-Moody_extension_intro}) over $\Sp(M)(\cO)$ yields a splitting of (\ref{meta_extension_small}) over $\Sp(M)(\cO)$. Consider the Hecke algebra
\begin{multline*}
\cH=\{f: \Sp(M)(\cO)\backslash \wh\Sp(M)(K)/\Sp(M)(\cO)\to\Qlb\mid  f(i(-1)g)=-f(g), \;\; g\in \wh\Sp(M)(K);\\
f\;\mbox{is of compact support}\}
\end{multline*}
The product is convolution, defined using the Haar measure on $\wh\Sp(M)(K)$ for which the inverse image of $\Sp(M)(\cO)$ has volume 1. 

  Set $G=\Sp(M)$. Let $\check{G}$ denote $\Sp_{2n}$ viewed as an algebraic group over $\Qlb$. Let $\Rep(\check{G})$ denote the category of finite-dimensional representations of $\check{G}$. Write $K(\Rep(\check{G}))$ for the Grothendieck ring of $\Rep(\check{G})$ over $\Qlb$. There is a canonical isomorphism of $\Qlb$-algebras 
$$
\cH\,\iso\, K(\Rep(\check{G}))
$$ 

 Actually, a categorical version of this isomorphism is proved. Consider the affine grassmanian $\Gr_G=G(K)/G(\cO)$, viewed as an ind-scheme over $k$. Let $W$ denote the nontrivial $\ell$-adic local system of rank one on $\Gm$ corresponding to the covering $\Gm\to\Gm$, $x\mapsto x^2$. Denote by $\Sph(\wt\Gr_G)$ the category of $G(\cO)$-equivariant perverse sheaves on $\ov{G}(K)/G(\cO)$, which are also $(\Gm, W)$-equivariant. Here $\wt\Gr_G$ denotes the stack quotient of $\ov{G}(K)$ by $\Gm$ with respect to the action $g\toup{x} x^2g$, $x\in \Gm, g\in \ov{G}(K)$. 
Actually, $\Sph(\wt\Gr_G)$ is a full subcategory of the category of perverse sheaves on $\wt\Gr_G$.
 
 Assuming for simplicity $k$ algebraically closed, we equip $\Sph(\wt\Gr_G)$ with the structure of a rigid tensor category. We establish a canonical equivalence of tensor categories 
$$
\Sph(\wt\Gr_G)\;\iso\; \Rep(\check{G})
$$   
 
\medskip\noindent
1.2  In the global setting let $X$ be a smooth projective curve over $k$. Let $G$ denote the sheaf of automorphisms of $\cO_X^n\oplus\Omega^n$ (now $\Omega$ is the canonical line bundle on $X$) preserving the symplectic form $\wedge^2(\cO_X^n\oplus\Omega^n)\to\Omega$. The stack $\Bun_G$ of $G$-bundles (=$G$-torsors) on $X$ classifies vector bundles $M$ of rank $2n$ on $X$, given with a nondegenerate symplectic form $\wedge^2 M\to\Omega$. We introduce an algebraic stack $\Bunt_G$ of metaplectic bundles on $X$. The stack $\wt\Gr_G$ is a local version of $\Bunt_G$. 
The category $\Sph(\wt\Gr_G)$ acts on $\D(\Bunt_G)$ by Hecke operators. 

  We construct a perverse sheaf $\Aut$ on $\Bunt_G$, a geometric analog of the Weil representation. We calculate the fibres of $\Aut$ and its constant terms for maximal parabolic subgroups of $G$. Finally, we argue that $\Aut$ is a Hecke eigensheaf on $\Bunt_G$ with eigenvalue 
$$
 \St=\RG(\PP^{2n-1},\Qlb)\otimes\Qlb[1](\frac{1}{2})^{\otimes 2n-1}
$$
viewed as a constant complex on $X$. Note that $\St$ is equipped with an action of $\SL_2$ of Arthur, the corresponding representation of $\SL_2$ is irreducible of dimension $2n$ and admits a unique, up to a multiple, symplectic form. One may imagine that $\Aut$ corresponds to a group homomorphism $\pi_1(X)\times\SL_2\to \check{G}$ trivial on $\pi_1(X)$. This agrees with Arthur's conjectures. 

\bigskip\medskip
\centerline{\scshape 2. Weil representation and motivations}

\bigskip\noindent
2.1 Let $X$ be a smooth projective absolutely irreducible curve over $k=\Fq$, 
$F=\Fq(X)$, $\AA$ be the adeles rings of $F$, $\cO\subset\AA$ be the entire
adeles. Assume that $q$ is odd. Fix a prime $\ell$ that does not divide $q$.
Let $\Omega$ denote the canonical line bundle on $X$.

 Let $M$ be a $2n$-dimensional vector space over $F$ with symplectic form $\wedge^2 M\to \Omega_F$, where $\Omega_F$ is the generic fibre of $\Omega$. The Heisenberg group
$H(M)=M\oplus \Omega_F$ with operation
$$
(m_1,\omega_1)(m_2,\omega_2)=(m_1+m_2, \, \omega_1+\omega_2+\frac{1}{2}\<m_1,m_2\>)
$$  
is algebraic over $F$. Fix a nontrivial additive character $\psi:\Fq\to \Qlb^*$. Then $H(M)(\AA)=M(\AA)\oplus\Omega(\AA)$ admits a canonical central character 
$\chi:\Omega(\AA)/\Omega(F)\to \Qlb^*$ given by 
$$
\chi(\omega)=\psi(\sum_{x\in X}
\tr_{k(x)/k}\Res \omega_x)
$$

 The Stone and Von Neumann theorem (\cite{MVW}) says that there is a unique (up to isomorphism) smooth irreducible representation $(\rho, \cS_{\psi})$ of $H(M)(\AA)$ over $\Qlb$ with central character $\chi$. The group $\Sp(M)$ acts on $H(M)$ by group automorphisms
$
(m,\omega)\toup{g}(gm,\omega)
$.
This defines the global metaplectic group\footnote{the notation $\tSp(M)(\AA)$ is ambiguous, these are not $\AA$-points of an algebraic group.} 
\begin{multline*}
\;\;\;\tSp(M)(\AA)=\{(g,M[g])\mid g\in\Sp(M)(\AA), M[g]\in\Aut \cS_{\psi}\\
\!\!\!\!\rho(gm,\omega)\comp M[g]=M[g]\comp \rho(m,\omega)\;\; \mbox{for}\; (m,\omega)\in H(M)(\AA)\}\;\;\;
\end{multline*}
included into an exact sequence 
\begin{equation}
\label{metaplectic_def}
1\to \Qlb^*\to\tSp(M)(\AA)\to \Sp(M)(\AA)\to 1
\end{equation}
The representation of $\tSp(M)(\AA)$ on $S_{\psi}$ is called the Weil (or oscillator) representation (\cite{W}). 

  For a subgroup $K\subset \Sp(M)(\AA)$ write $\tilde K$ for the preimage of $K$ in $\tSp(M)(\AA)$. Since $\chi$ is trivial on $\Omega_F$, one may talk about $H(M)$-invariant functionals on $S_{\psi}$, they are called theta-functionals. The space of theta-functionals is 1-dimensional and preserved by $\tSp(M)(F)$, so the action of $\tSp(M)(F)$ on this space defines a splitting of (\ref{metaplectic_def}) over $\Sp(M)(F)$. 
  
  View
$$
\Funct(\Sp(M)(F)\backslash \tSp(M)(\AA))=\{f: \Sp(M)(F)\backslash \tSp(M)(\AA))\to \Qlb\}
$$
as a representation of $\tSp(M)(\AA)$ by right translations. A theta-functional $\Theta: S_{\psi}\to\Qlb$ defines a morphism of $\tSp(M)(\AA)$-modules 
\begin{equation}
\label{hom_modules}
\cS_{\psi}\to \Funct(\Sp(M)(F)\backslash
\tSp(M)(\AA))
\end{equation}
sending $\phi$ to $\theta_{\phi}$ given by $\theta_{\phi}(g)=\Theta(g\phi)$ for $g\in \tSp(M)(\AA)$.
 
  Now assume that $M$ is actually a rank $2n$ vector bundle on $X$ with symplectic form $\wedge^2 M\to\Omega$. Then we get the subgroups $\Sp(M)(\cO)\subset \Sp(M)(\AA)$ and  
$M(\cO)\oplus\Omega(\cO)\subset H(M)(\AA)$. Moreover, the space of $M(\cO)\oplus \Omega(\cO)$-invariants in $\cS_{\psi}$ is 1-dimensional and preserved by $\tSp(M)(\cO)$. The action of $\tSp(M)(\cO)$ on this space yields a splitting of (\ref{metaplectic_def}) over $\Sp(M)(\cO)$. If $\phi_0\in \cS_{\psi}$ is a nonzero  $M(\cO)\oplus \Omega(\cO)$-invariant vector then its image under (\ref{hom_modules}) is the classical theta-function
$$
f_0: \Sp(M)(F)\backslash\tSp(M)(\AA))/\Sp(M)(\cO)\to\Qlb
$$
that we are going to geometrize. 

 Let $G$ denote the sheaf of automorphisms of $M$ preserving the form $\wedge^2 M\to\Omega$. This is a sheaf of groups (in flat topology) on $X$ locally in Zarisky topology isomorphic to $\Sp_{2n}$. 

\medskip\noindent
2.2  Assume $M=V\oplus (V^*\otimes\Omega)$ is a direct sum of lagrangian subbundles, the form being given by the canonical pairing $\<.,.\>$ between $V$ and $V^*$. Let 
$$
\chi_V: V(\AA)\oplus \Omega(\AA)\to\Qlb^*
$$ 
denote the character $\chi_V(v,\omega)=\chi(\omega)$. 

 We have the subgroup $V(\AA)\subset H(M)(\AA)$.
The space of $V(\AA)$-invariant functionals on $\cS_{\psi}$ is 1-dimensional. A choice of such functional identifies $\cS_{\psi}$ with the induced representation of $(V(\AA)\oplus\Omega(\AA), \chi_V)$ to  $H(M)(\AA)$. The latter  identifies with the Schwarz space $\cS(V^*\otimes\Omega(\AA))$ of locally constant compactly supported $\Qlb$-valued functions on $V^*\otimes\Omega(\AA)$, the corresponding functional on $\cS(V^*\otimes\Omega(\AA))$ becomes the evaluation at zero $\ev: \cS(V^*\otimes\Omega(\AA))\to\Qlb$. This is the Schr\" odinger model of $\cS_{\psi}$.

 Write $g\in\Sp(M)(\AA)$ as a matrix 
\begin{equation}
\label{def_g}
g=\left(
\begin{array}{cc}
a & b\\
c & d
\end{array}
\right),
\end{equation}
with $a\in\End(V)(\AA), b\in\Hom(V^*\otimes\Omega, V)(\AA), d\in \End(V^*)(\AA),
c\in\Hom(V,V^*\otimes\Omega)(\AA)$. Write $a^*$ for the transpose operator to $a$.

The defined up to a scalar automorphism $M[g]$ of $\cS(V^*\otimes\Omega(\AA))$ is described as follows. 
\begin{itemize}
\item For $a\in\GL(V)(\AA)$ we have $\left(
\begin{array}{cc}
a & 0\\
0 & {a^{*}}^{-1}
\end{array}
\right)\in\Sp(M)(\AA)$. Besides, $\left(
\begin{array}{cc}
1 & b\\
0 & 1
\end{array}
\right)\in\Sp(M)(\AA)$ if and only if $b\in (V\otimes V\otimes\Omega^{-1})(\AA)$ is symmetric. For $g$ given by (\ref{def_g}) with $c=0$ we have
\begin{equation}
\label{formula_interwine1}
(M[g]f)(v^*)=\chi(\frac{1}{2}\< a^*v^*, b^*v^*\>)f(a^*v^*), \;\;\; v^*\in V^*\otimes\Omega(\AA)
\end{equation}

\item if $b: V^*\otimes\Omega(\AA)\,\iso\,V(\AA)$ then $g=\left(
\begin{array}{cc}
0 & b\\
-{b^*}^{-1} & 0
\end{array}
\right)\in\Sp(M)(\AA)$ and  
\begin{equation}
\label{formula_interwine2}
(M[g]f)(v^*)=\int_{V(\AA)} \chi(\< v, v^*\>) f(b^{-1}v)dv, \;\;\; v^*\in V^*\otimes\Omega(\AA)
\end{equation}
for any Haar measure $dv$ on $V(\AA)$.
\end{itemize}

 Let $P\subset G$ denote the Siegel parabolic subgroup preserving $V$. 
The subgroup $\tilde P(\AA)$ preserves $\ev$ up to a multiple, so defining a splitting of (\ref{metaplectic_def}) over $P(\AA)$. This splitting coincides with the one given by (\ref{formula_interwine1}).

  Let $\phi_0\in\cS(V^*\otimes\Omega(\AA))$ denote the characteristic function of $V^*\otimes\Omega(\cO)$. Using (\ref{formula_interwine1}) and (\ref{formula_interwine2}) one shows that $\phi_0$ generates the space of $\Sp(M)(\cO)$-invariants in $\cS(V^*\otimes\Omega(\AA))$. In this model of $\cS_{\psi}$ the theta functional $\Theta: \cS(V^*\otimes\Omega(\AA))\to\Qlb$ is given by
$$
\Theta(\phi)=\sum_{v^*\in V^*\otimes\Omega(F)} \phi(v^*) \;\;\;\;\;\mbox{for}\;\; \phi\in \cS(V^*\otimes\Omega(\AA))
$$
Let $f_0$ denote the image of $\phi_0$ under the corresponding map (\ref{hom_modules}). 
Let us calculate the composition
$$
P(F)\backslash P(\AA)/P(\cO)\to \Sp(M)(F)\backslash \tSp(M)(\AA)/\Sp(M)(\cO)\toup{f_0}\Qlb
$$
denoted $f_P$. We used the fact that the splittings of (\ref{metaplectic_def}) over $P(\AA)$ and $\Sp(M)(\cO)$ are compatible over $P(\cO)$. 

 Denote by $\Bun_n$ the $k$-stack of rank $n$ vector bundles on $X$.
The set $\GL(V)(\AA)/\GL(V)(\cO)$ naturally identifies with the 
isomorphism classes of pairs $(L,\alpha)$, where $L\in\Bun_n(k)$ and 
$\alpha: L(F)\,\iso\, V(F)$. Here $L(F)$ is the generic fibre of $L$. 
 
 Let $a\in\GL(V)(\AA)$ and
$(L,\alpha)$ be the pair attached to $a\GL(V)(\cO)$. Then 
\begin{equation}
\label{iso_1}
\{v^*\in V^*\otimes\Omega(F)\mid a^*v^*\in V^*\otimes\Omega(\cO)\}\,\toup{\alpha^*}\, \Hom(L,\Omega)
\end{equation}
is an isomorphism.

 The group $P$ fits into an exact sequence $1\to (\Sym^2 V)\otimes\Omega^{-1} \to P\to \GL(V)\to 1$ 
of algebraic groups over $X$. For $g\in P(\AA)$ we get 
\begin{multline*}
f_P(g)=\Theta(g\phi_0)=\sum_{v^*\in V^*\otimes\Omega(F)}
(g\phi_0)(v^*)=
\sum_{v^*\in V^*\otimes\Omega(F)} \chi(\frac{1}{2}\<a^*v^*, b^* v^*\>)\phi_0(a^* v^*)=\\
\sum_{s\in \Hom(L,\Omega)} \chi(\frac{1}{2}\<s, ab^* s\>)
\end{multline*}
in view of (\ref{iso_1}).

 Let $\Bun_P$ be the $k$-stack of
$P$-bundles on $X$. Its $Y$-points for a scheme $Y$ is the category of $(Y\times X)\times_X P$-torsors over $Y\times X$. Then $\Bun_P$ classifies pairs $L\in\Bun_n$ together with an exact sequence on $X$
\begin{equation}
\label{sequence_sect2}
0\to \Sym^2 L\to ?\to\Omega\to 0
\end{equation}
(More generally, for a semidirect product of group schemes $1\to U\to P\to M\to 1$ providing a $P$-torsor $\cF_P$ is equivalent to providing a $M$-torsor $\cF_M$ and a $U_{\cF_M}$-torsor of isomorphisms $\Isom(\cF_P, \cF_M\times_M P)$ inducing a given one on the corresponding $M$-torsors).

 In view of the bijection $P(F)\backslash P(\AA)/P(\cO)\,\iso\, \Bun_P(k)$, the function $f_P$ on $\Bun_P(k)$ is described as follows. Let a $P$-torsor $\cF_P\in\Bun_P(k)$ be given by
$L\in\Bun_n(k)$ together with (\ref{sequence_sect2}). Consider the map
$q^{\cF_P}: \Hom(L,\Omega)\to k$
sending $s\in \Hom(L,\Omega)$ to the pairing of 
$$
s\otimes s\in\Hom(\Sym^2 L, \Omega^{\otimes 2})
$$ 
with the exact sequence (\ref{sequence_sect2}). Then
$$
f_P(\cF_P)=\sum_{s\in \Hom(L,\Omega)} \psi(q^{\cF_P}(s))
$$

  The function $f_P: \Bun_P(k)\to \Qlb$ is the trace of Frobenius of 
the following $\ell$-adic complex $S_{P,\psi}$ on $\Bun_P$. 
  
 Let $p: \cX\to\Bun_P$ be the stack over $\Bun_P$ with fibre $\Hom(L,\Omega)$. Let $q:\cX\to\A^1$ be the map sending $s\in \Hom(L,\Omega)$ to the pairing of (\ref{sequence_sect2}) with
$$
s\otimes s\in\Hom(\Sym^2 L, \Omega^{\otimes 2})
$$
The geometric analog of $f_P$ is the complex
$S_{P,\psi}=p_!q^*\cL_{\psi}\otimes\Qlb[1](\frac{1}{2})^{\otimes\dim\cX}$ on $\Bun_P$, here $\dim\cX$ denotes the dimension of the corresponding connected component of $\cX$.

\bigskip\bigskip
\centerline{\scshape 3. Main results}

\bigskip\noindent
3.1 {\scshape Notation\ } From now on $k$ denotes an algebraically closed field
of characteristic $p>2$, all the schemes (or stacks) we consider are defined
over $k$. 

 Let $X$ be a smooth projective connected curve. Write $\Omega$ for the canonical line bundle on $X$. Fix a prime $\ell\ne p$. For a scheme (or stack) $S$ write $\D(S)$ for the bounded derived category of $\ell$-adic \'etale sheaves on $S$, and $\P(S)\subset \D(S)$ for the category of perverse sheaves (the middle perversity function is always taken in absolut sense over $\Spec k$). 

 Fix a nontrivial character $\psi: \Fp\to\Qlb^*$ and denote by
$\cL_{\psi}$ the corresponding Artin-Shreier sheaf on $\A^1$. Fix a square root 
$\Qlb(\frac{1}{2})$ of the sheaf $\Qlb(1)$ on $\Spec \Fq$. Isomorphism classes of such correspond to square roots of $q$ in $\Qlb$. Fix an inclusion of fields $\Fq\hook{} k$.

 If $V\to S$ and $V^*\to S$ are dual rank $n$ vector bundles over a stack $S$, we normalize 
the Fourier trasform $\Four_{\psi}: \D(V)\to\D(V^*)$ by 
$\Four_{\psi}(K)=(p_{V^*})_!(\xi^*\cL_{\psi}\otimes p_V^*K)[n](\frac{n}{2})$,  
where $p_V, p_{V^*}$ are the projections, and $\xi: V\times_S V^*\to \A^1$ is the pairing.
  
 A $G$-torsor on a scheme $S$ is also referred to as a $G$-bundle on $S$. 
Write $\Vect^{\epsilon}$ for the tensor category of $\ZZ/2\ZZ$-graded vector spaces, our conventions about this category are those of (\cite{D}).
Write $\Vect\subset \Vect^{\epsilon}$ for its even component, i.e., the tensor category of vector spaces.  

\medskip\noindent
3.1.1  The sheaf (in flat topology) on the category of $k$-schemes
represented by $\mu_2:=\Ker(x\mapsto x^2: \Gm\to\Gm)$ is the constant sheaf $\{\pm 1\}$. 

 For a scheme $S$ and a line bundle $\cA$ on $S$ denote by $\tilde S$
the following $\mu_2$-gerbe over $S$. For an $S$-scheme $S'$, the category of
$S'$-points of $\tilde S$ is the category of pairs $(\cB,
\;\cB^2\,\iso\,\cA\mid_{S'})$, where $\cB$ is a line bundle on $S'$. Note that
$\tilde S\to S$ is \'etale.

 If $\tilde S\to S$ admits a section given by invertible $\cO_S$-module $\cB_0$ together
with $\cB_0^2\,\iso\,\cA$ then the gerbe is trivial, that is, $\tilde S\,\iso B(\mu_2/S)$
over $S$. In this case we get the $S_2$-covering $\Cov(\tilde S)\to \tilde S$,
whose fibre consists of isomorphisms $\cB\,\iso\,\cB_0$ whose square is the given one
$\cB^2\,\iso\,\cA$. This covering is locally trivial in \'etale topology, but not trivial even for $S=\Spec k$. Actually $S=\Cov(\tilde S)$.

\medskip\noindent
3.1.2 If in addition $\cA$ is a $\ZZ/2\ZZ$-graded line bundle on $S$ purely of degree zero, then by definition $\tilde S$ classifies a $\ZZ/2\ZZ$-graded line bundle $\cB$ purely of degree zero, given with a $\ZZ/2\ZZ$-graded isomorphism
$\cB^2\,\iso\, \cA$. If $\cB$ is a $\ZZ/2\ZZ$-graded line bundle on $S$ \select{of pure degree} (that is, placed in one degree only over each connected component) then a $\ZZ/2\ZZ$-graded isomorphism $\cB^2\,\iso\,\cA$ yields a (uniquely defined) section of $\tilde S$.  

\medskip\noindent
3.2 Let $\Bun_n$ be the stack of rank $n$ vector bundles on $X$. Let $G$ denote the sheaf of automorphisms of $\cO_X^n\oplus \Omega^n$ preserving the symplectic form $\wedge^2 (\cO_X^n\oplus \Omega^n)\to\Omega$. So, $G$ is a sheaf of groups in flat topology on the category of $X$-schemes. 
 
 The stack $\Bun_G$ of $G$-bundles on $X$ classifies
$M\in\Bun_{2n}$ together with a symplectic form $\wedge^2 M\to\Omega$. This is a smooth algebraic stack locally of finite type over $k$. Since $G$ is simply-connected, $\Bun_G$ is irreducible (\cite{BD}, 2.2.1). Let $d_G=\dim\Bun_G=(g-1)\dim \LieSp_{2n}$. To express the dependence on $n$ we write $G_n$, $\Bun_{G_n}$, $d_{G_n}$ and so on.

 Denote by $\cA$ the line bundle on $\Bun_G$ whose fibre at $M$ is
$\det\RG(X,M)$ (cf. \cite{D}). As $\chi(M)=0$, we view $\cA$ as a $\ZZ/2\ZZ$-graded line bundle placed in degree zero. It yields a $\mu_2$-gerbe
\begin{equation}
\label{gerbe_over_Bun}
\gr: \Bunt_G\to \Bun_G
\end{equation}
So, $S$-points of $\Bunt_G$ is the category: a line bundle $\cB$ on $S$, 
a vector bundle $M$ on $S\times X$ of rank $2n$ with symplectic form $\wedge^2
M\to \Omega_{S\times X/S}$, and an isomorphism of $\cO_S$-modules
$\cB^2\,\iso\, \det\RG(X,M)$. 

 The idea of using the determinant of cohomology was communicated to me by 
G. Laumon and goes back to P. Deligne \cite{D2}.

 Let $_i\Bun_G\hook{} \Bun_G$ be the locally closed substack
given by $\dim\H^0(X,M)=i$. Let $_i\Bunt_G$ denote the preimage of
$_i\Bun_G$ under $\gr$. 

\begin{Lm} 
\label{Lm_nonempty}
Each stratum $_i\Bun_G$ of $\Bun_G$ is nonempty.
\end{Lm}
\begin{Prf} 
 For $n=1$ take $M=\cA(D)\oplus (\cA^*\otimes\Omega(-D))$, where $D$ is an effective divisor of degree $i$ on $X$, and $\cA$ is a line bundle on $X$ of degree $g-1$ such that $\H^0(X,\cA)=\H^1(X, \cA)=0$. Such $\cA$ exist, because $\dim X^{(g-1)}=g-1$, and the dimension of the Picard scheme of $X$ is $g$. Then 
$\dim\H^0(X, M)=i$.

 For any $n$ construct $M\in{_i\Bun_G}$ as  $M=M_1\oplus\ldots\oplus M_n$ with
$M_j\in {_{i_j}\Bun_{G_1}}$ for some $i_1+\ldots+i_n=i$.
\end{Prf}

\medskip

 We have a line bundle $_i\cB$ on $_i\Bun_G$ whose fibre at $M\in\Bun_G$ is $\det\H^0(X,M)$. View it as a $\ZZ/2\ZZ$-graded placed in degree $\dim\H^0(X,M)$ modulo 2. Then for each $i$ we get a $\ZZ/2\ZZ$-graded isomorphism $_i\cB^2\,\iso\,\cA\mid_{_i\Bun_G}$. By 3.1.2,  the gerbe $_i\Bunt_G\to {_i\Bun_G}$ is trivial. So, we have the two-sheeted covering 
$$
_i\rho: \Cov(_i\Bunt_G)\to {_i\Bunt_G}
$$
The line bundles $_i\cB$ (viewed as ungraded) do not glue into a line bundle over $\Bun_G$ (the gerbe $\gr$ is nontrivial, because $\cA$ is a generator of the Picard group $\Pic(\Bun_G)\,\iso\,\ZZ$ by \cite{F}).

\begin{Def}
\label{Def_1}
 For each $i$ define a local system $_i\Aut$ on $_i\Bunt_G$ by 
$$
_i\Aut=\Hom_{S_2}(\sign, {_i\rho_!}\Qlb)
$$
Let $\Aut_g\in \P(\Bunt_G)$ (resp., $\Aut_s\in\P(\Bunt_G)$) denote the Goresky-MacPherson
extension of $_0\Aut\otimes\,\Qlb[d_G](\frac{d_G}{2})$ (resp., of
$_1\Aut\otimes\,\Qlb[d_G-1](\frac{d_G-1}{2})$) under
$_i\Bunt_G\hook{}\Bunt_G$. \footnote{Here `g' stands for generic and `s' for special.
We postpone to Proposition~\ref{Pp_up_to_minus_one} the proof of the fact that
$_1\Aut$ is a shifted perverse sheaf on $_1\Bunt_G$}  Set 
$$
\Aut=\Aut_g\oplus\Aut_s
$$ 
By construction, $\DD(\Aut)\,\iso\,\Aut$ canonically.
\end{Def}

 Here is our main result.

\begin{Th} 
\label{Th_1} For each $i$ the $*$-restriction $\Aut\mid_{_i\Bunt_G}$ identifies with
$$
\Aut\mid_{_i\Bunt_G}\,\iso\, {_i\Aut}\otimes\Qlb[1](\frac{1}{2})^{\otimes d_G-i},
$$
(once $\sqrt{-1}\in k$ is fixed, the corresponding isomorphism is well-defined up to a sign). The $*$-restriction of $\Aut_g$ (resp., of $\Aut_s$) to  $_i\Bunt_G$ vanishes for $i$ odd (resp., even). 
\end{Th}

\begin{Rem}
Classicaly, for two symplectic spaces $W,W'$
there is a natural map $\tSp(W)\times\tSp(W')\to \tSp(W\oplus W')$, and the restriction
of the metaplectic representation under this map is the tensor product of metaplectic
representations of the factors (\cite{Pr}, Remark~2.7).

 In geometric setting we have a map
$s_{n,m}:\Bun_{G_n}\times\Bun_{G_m}\to \Bun_{G_{n+m}}$ sending $M,M'$ to $M\oplus M'$. 
It extends to a map 
$$
\tilde s_{n,m}: \Bunt_{G_n}\times\Bunt_{G_m}\to \Bunt_{G_{n+m}}
$$
sending $(M,\cB, \cB^2\,\iso\, \det\RG(X,M))$ and $(M',\cB', \cB'^2\,\iso\,\det\RG(X,M'))$
to 
$$
(M\oplus M', \cB\otimes\cB', \cB^2\otimes\cB'^2\,\iso\,
\det\RG(X,M)\otimes\det\RG(X,M')\,\iso\,\det\RG(X,M\oplus M'))
$$ 
The restriction yields a map $s_{n,m}:
{_i\Bun_{G_n}}\times{_j\Bun_{G_m}}\to {_{i+j}\Bun_{G_{n+m}}}$ and we get  canonically
$s_{n,m}^*(_{i+j}\cB)\,\iso\, {_i\cB}\boxtimes {_j\cB}$. For any $i,j$ this yields an
isomorphism
$$
\tilde s_{n,m}^*(_{i+j}\Aut)\,\iso\, {_i\Aut}\boxtimes {_j\Aut}
$$
of local systems on $_i\Bunt_{G_n}\times {_j\Bunt_{G_m}}$. 
Thus, 
$$
\tilde s_{n,m}^*\Aut_g\otimes\Qlb[1](\frac{1}{2})^{\otimes d_{G_n}+d_{G_m}-d_{G_{n+m}}}
\,\iso\,
(\Aut_g\boxtimes\Aut_g)\oplus(\Aut_s\boxtimes\Aut_s)
$$
and
$$
\tilde s_{n,m}^*\Aut_s\otimes\Qlb[1](\frac{1}{2})^{\otimes
d_{G_n}+d_{G_m}-d_{G_{n+m}}}\,\iso\,
(\Aut_g\boxtimes\Aut_s)\oplus (\Aut_s\boxtimes\Aut_g)
$$
in the completed Grothendieck group $K(\Bunt_{G_n}\times\Bunt_{G_m})$ (the completion is with respect to the filtration given by the codimension of support). 
\end{Rem}

\medskip\noindent
3.3 For $1\le k\le n$ denote by $\Bun_{P_k}$ the stack classifying $M\in\Bun_G$
together with an isotropic subbundle $L_1\subset M$ of rank $k$. We write
$L_{-1}\subset M$ for the orthogonal complement of $L_1$, so a point of
$\Bun_{P_k}$ gives rise to a flag $(L_1\subset L_{-1}\subset M)$, and
$L_{-1}/L_1\in \Bun_{G_{n-k}}$ naturally. 

 Write $\nu_k:\Bun_{P_k}\to\Bun_G$ for the projection. Define the map
$$
\tilde\nu_k:
\Bunt_{G_{n-k}}\times_{\Bun_{G_{n-k}}}\Bun_{P_k}\to \Bunt_G
$$
as follows. An
$S$-point of the source is given by $(L_1\subset L_{-1}\subset M)\in\Bun_{P_k}(S)$
together with a ($\ZZ/2\ZZ$-graded of pure degree zero) invertible $\cO_S$-module $\cB$ and 
$\cB^2\,\iso\,\det\RG(X,L_{-1}/L_1)$. We have a canonical isomorphism of $\ZZ/2\ZZ$-graded lines
\begin{equation}
\label{equation_sect3.3}
\det\RG(X,L_1)\otimes\det\RG(X, L_{-1}/L_1)\otimes\det\RG(X,
L_1^*\otimes\Omega)\,\iso\, \det\RG(X,M)
\end{equation}
The map $\tilde\nu_k$ sends this point to $M\in\Bun_G$ together with an invertible
$\cO_S$-module $\cB'=\cB\otimes\det\RG(X,L_1)$ and $\cB'^2\,\iso\,\det\RG(X,M)$
given by (\ref{equation_sect3.3}). Since $\cB'$ is of pure degree as $\ZZ/2\ZZ$-graded, the map is well-defined by 3.1.2. 

 Let $\Bun_{Q_k}$ be the stack of collections: an exact sequence
$0\to L_1\to L_{-1}\to L_{-1}/L_1\to 0$ of vector bundles on $X$ with $L_1\in\Bun_k$
and $L_{-1}/L_1\in\Bun_{2n-2k}$, and a symplectic form $\wedge^2(L_{-1}/L_1)\to\Omega$
(thus, $L_{-1}/L_1\in\Bun_{G_{n-k}}$).  

 Let $\eta_k:\Bun_{P_k}\to\Bun_{Q_k}$ denote the natural projection.  Let
$^0\Bun_{Q_k}\subset\Bun_{Q_k}$ be the open substack given by $\H^0(X, \Sym^2 L_1)=0$.

\begin{Th} 
\label{Th_2}
For the diagram
$$
\Bunt_{G_{n-k}}\times_{\Bun_{G_{n-k}}}\Bun_{Q_k}\,\getsup{\id\times\eta_k}\;
\Bunt_{G_{n-k}}\times_{\Bun_{G_{n-k}}}\Bun_{P_k}\toup{\tilde\nu_k}\Bunt_G
$$
we have an isomorphism 
$$
(\id\times\eta_k)_!\tilde\nu_k^*\Aut\,\iso\, \Aut\boxtimes\,\Qlb[b](\frac{b}{2})
$$
over $\Bunt_{G_{n-k}}\times_{\Bun_{G_{n-k}}}{^0\Bun_{Q_k}}$. (Once $\sqrt{-1}\in k$ is fixed, the isomorphism is well-defined up to a sign on generic and special parts). Here $b(L_1)=d_G-d_{G_{n-k}}-
\chi(L_1)+2\chi(\Omega^{-1}\otimes\Sym^2 L_1)$ is a function of a connected component of $^0\Bun_{Q_k}$. 
 If $\chi(L_1)$ is even then, over the corresponding connected component, the above isomorphism preserves generic and special parts, othewise it interchanges them. 
\end{Th}

\medskip\noindent
3.4 In Sect.~8.1 we consider the affine Grassmanian $\Gr_G$ for $G$, it is equipped with a natural line bundle $\cL$ that generates the Picard group of $\Gr_G$. Let $\wt\Gr_G\to\Gr_G$ denote the $\mu_2$-gerbe of square roots of $\cL$. This is a local version of the gerbe (\ref{gerbe_over_Bun}). We introduce the category $\Sph(\wt\Gr_G)^{\flat}$ of \select{genuine spherical sheaves} on $\wt\Gr_G$ (cf. Definition~\ref{Def_genuine} and \ref{Def_Sph_flat}). 

 As for usual spherical sheaves on the affine Grassmanian, we equip $\Sph(\wt\Gr_G)^{\flat}$ with a structure of a rigid tensor category. Main result of Sect.~8 is the following version of the Satake equivalence.
 
\begin{Th} 
\label{Th_Satake}
The category $\Sph(\wt\Gr_G)^{\flat}$ is canonically equivalent, as a tensor category, to the category $\Rep(\Sp_{2n})$ of finite-dimensional $\Qlb$-representations of $\Sp_{2n}$.
\end{Th}

 In Sect.~9 we define for $K\in\Sph(\wt\Gr_G)^{\flat}$ Hecke operators $\H(K,\cdot): \D(\Bunt_G)\to \D(X\times\Bunt_G)$ compatible with the tensor structure on $\Sph(\wt\Gr_G)^{\flat}$. Finally, we prove Theorem~\ref{Th_Hecke} saying that 
$\Aut$ is a Hecke eigen-sheaf with eigenvalue 
$$
 \St=\RG(\PP^{2n-1},\Qlb)\otimes\Qlb[1](\frac{1}{2})^{\otimes 2n-1}
$$
viewed as a constant complex on $X$. 

\bigskip\noindent
\centerline{\scshape 4. Finite-dimensional model}

\medskip\noindent
4.1 Let $V$ be a $k$-vector space of dimension $d$. Write $\ST^2(V^*)$ for the space of symmetric tensors in $V^*\otimes V^*$, this is the space of symmetric bilinear forms on $V$. Think of $b\in \ST^2(V^*)$ as a map $b:V\to V^*$ such that $b^*=b$. 
Let $p: V\times\ST^2(V^*)\to \ST^2(V^*)$ denote the projection. Let
$\beta: V\times\ST^2(V^*)\to\A^1$ be the map that sends $(v,b)$ to $\<v, bv\>$.
Set 
$$
S_{\psi}=p_!\beta^*\cL_{\psi}\otimes\Qlb[1](\frac{1}{2})^{\otimes d+\frac{1}{2}d(d+1)}
$$
 
 Let $\pi:V\to \Sym^2V$ be the map $v\mapsto v\otimes v$. Then 
\begin{equation}
\label{def_S_psi}
S_{\psi}=\Four_{\psi}(\pi_!\Qlb[d](\frac{d}{2}))
\end{equation}
The map $\pi$ is finite, and $\pi_!\Qlb=
\cL_0+\cL_1$, where $\cL_0$ is the constant sheaf on the image $\Im\pi$ of $\pi$, and 
$\cL_1$ is a nontrivial local system of rank one on $\Im\pi-\{0\}$ extended by zero to
$\Im\pi$. So, $S_{\psi}$ is a direct
sum of two irreducible perverse sheaves.

\begin{Lm} $S_{\psi}$ is $\GL(V)$-equivariant.
\end{Lm}
\begin{Prf} Clearly, $\pi_!\Qlb$ is $\GL(V)$-equivariant. The Fourier transform preserves
$\GL(V)$-equivariance of a perverse sheaf. 
\end{Prf}

\medskip

 Stratify $\ST^2(V^*)$ by $Q_i(V)$, where $Q_i(V)$ is the locus of $b:V\to V^*$ such that $\dim\Ker b=i$. 
For $b\in \ST^2 (V^*)$ denote by $\beta_b:V\to\A^1$ the map $b\mapsto \<v,bv\>$.
We have a usual ambiguity in identifying $\ST^2(V^*)$ with $\Sym^2 (V^*)$: $b$ goes to $\beta_b$ or $\frac{1}{2}\beta_b$. Since $S_{\psi}$ is $\GL(V)$-equivariant, we can view it as a perverse sheaf on $\Sym^2(V^*)$ unambiguously.

\begin{Lm}
\label{Lm_gamma_function}
For $b\in Q_0(V)$ the complex 
$\RG_c(V, \beta_b^*\cL_{\psi})$ is a 1-dimensional vector space
placed in degree $d$.
\end{Lm}
\begin{Prf}
In some basis $\beta_b$ is given by $(x_1,\ldots, x_d)\mapsto x_1^2+\ldots+x_d^2$. Thus we
may assume $d=1$.
Consider the map $\pi: \A^1\to \A^1$ given by $\pi(x)=x^2$. As above
$\pi_!\Qlb\,\iso\, \cL_0\oplus \cL_1$ with $\cL_0=\Qlb$. 
We get $\RG_c(\A^1, \pi^*\cL_{\psi})\,\iso\, \RG_c(\Gm, \cL\otimes\cL_{\psi})$. The
latter is a vector space of dimension one placed in degree one 
(a gamma-function on $\Gm$). 
\end{Prf}

\medskip

 Let $\Cov(Q_0(V))\to Q_0(V)$ denote the two-sheeted covering of $Q_0(V)$ whose fibre over
$b:V\,\iso\, V^*$ is the set of trivializations $\det V\,\iso\, k$ whose square is the
one induced by $b$. 

 The group $\GL(V)$ acts transitively on $Q_0(V)$, so given $b\in Q_0(V)$ one gets an
identification $Q_0(V)\,\iso\, \GL(V)/\OO(V,b)$. Our covering becomes the map
$\GL(V)/\SO(V,b)\to \GL(V)/\OO(V)$.

\smallskip

 More generally, $\GL(V)$ acts transitively on $Q_i(V)$. For $b\in Q_i(V)$ with
$\Ker b=V_0$, we can consider $b$ as an element of $\Sym^2(V/V_0)^*$. We get a
parabolic
$P_0\subset \GL(V)$ of automorphisms of $V$ that preserve $V_0$. Let $\rm{St}_{V_0}$
be the preimage of $\OO(V/V_0,b)$ under $P_0\to \GL(V/V_0)$. Then
$\rm{St}_{V_0}$ is the stabilizer of $b\in Q_i(V)$ in $\GL(V)$. Since $\SO(V,b)$ is
connected, for $i<d$ there is exactly  one (up to isomorphism) nonconstant
$\GL(V)$-equivariant local system of rank one on $Q_i(V)$. It corresponds to the
$S_2$-covering
$\Cov(Q_i(V))\to Q_i(V)$ whose fibre over $b$ is the set of trivializations
$\det(V/V_0)\,\iso k$ compatible with $b$.

\begin{Pp} 
\label{Pp_linear_algebra1}
1) The $*$-restriction of $S_{\psi}$ to $Q_i(V)$ is a
$\GL(V)$-equivariant local system of rank one placed in degree
$i-\frac{1}{2}d(d+1)$. For $i<d$ this local system is nonconstant and
comes from the covering $\Cov(Q_i(V))\to Q_i(V)$.\\
2) $S_{\psi}=S_{\psi,g}\oplus S_{\psi,s}$ is a direct sum of two irreducible perverse
sheaves. Here $S_{\psi,g}$ is
the Goresky-MacPherson extension of $S_{\psi}\mid_{Q_0(V)}$, and $S_{\psi,s}$ is the
Goresky-MacPherson extension of $S_{\psi}\mid_{Q_1(V)}$ under $Q_1(V)\hook{} Q_{\ge
1}(V)$.

\smallskip\noindent
3) We have $\DD S_{\psi,g}\,\iso\, S_{\psi^{-1},g}$ and $\DD S_{\psi,s}\,\iso\,
S_{\psi^{-1},s}$ canonically.

\smallskip\noindent
4) If $V=V_1\oplus V_2$ is a direct sum of two vector spaces of dimensions $d_1$ and $d_2$
then the $*$-restriction of
$S_{\psi}\otimes\Qlb[1](\frac{1}{2})^{\otimes-\frac{1}{2}d(d+1)}$ to the subspace
$\Sym^2(V_1^*)\oplus\Sym^2(V_2^*)$
is canonically 
$$
(S_{\psi}\boxtimes
S_{\psi})\otimes\Qlb[1](\frac{1}{2})^{\otimes
-\frac{1}{2}d_1(d_1+1)-\frac{1}{2}d_2(d_2+1)}
$$
\end{Pp}
\begin{Prf}
2) A point of $Q_i(V)$ is given by a subspace $V_0\subset V$
of dimension $i$ together with nondegenerate form $b: V/V_0\to (V/V_0)^*$ such that
$b^*=b$. It follows that
$$
\dim Q_i(V)=\frac{1}{2}(d-i)(d+1-i)+(d-i)i=\frac{1}{2}(d-i)(d+1+i)
$$
From Lemma~\ref{Lm_gamma_function} applied to $V/V_0$ we deduce that
$S_{\psi}\mid_{Q_i(V)}$ is a local system of rank one placed in degree
$i-\frac{1}{2}d(d+1)$. From (\ref{def_S_psi}) we see that $\DD
S_{\psi}\iso S_{\psi^{-1}}$. For $0\le i\le d$ we have
$$
\dim Q_i(V)=\frac{1}{2}(d-i)(d+1+i)\le \frac{1}{2}d(d+1)-i,
$$
the equality holds only for $i=0$ and $i=1$. So, $S_{\psi}$ is the Goresky-MacPherson
extension from the open subscheme $Q_{\le 1}(V)$. 

 Let $S_{\psi,g}$ be the intermediate extension of $S_{\psi}\mid_{Q_0(V)}$ to $\Sym^2
V^*$. The $*$-restriction $S_{\psi,g}\mid_{Q_1(V)}$ vanishes. Indeed, it should be placed
in strictly negative perverse degrees, but $S_{\psi}\mid_{Q_1(V)}$ is a perverse sheaf. 
Part 2) follows.

\medskip\noindent
3) follows from (\ref{def_S_psi})

\medskip\noindent
4) The composition $V_1\oplus V_2\;\iso\; V\toup{\pi} \Sym^2 V\toup{a} \Sym^2
V_1\times\Sym^2 V_2$ equals $\pi\times\pi$. So, $a_!\pi_!\Qlb\;\iso\;
(\pi_!\Qlb\boxtimes\pi_!\Qlb)$. Fourier transform interchanges $a_!$ and the
$*$-restriction under the transpose map $a^*: \Sym^2 V^*_1\times\Sym^2 V^*_2\to\Sym^2 V^*$. 
  
\medskip\noindent
1) Since $S_{\psi}\mid_{Q_i(V)}$ is $\GL(V)$-equivariant, it remains to show it is
nonconstant for $i<d$. 

\medskip
\Step 1 Start with $d=1$ case, so $Q_0(V)\,\iso\,\Gm$. To show that $S_{\psi}$ is
nonconstant on $Q_0(V)$ in this case, it suffices to prove that
$\RG_c(\Gm, S_{\psi})=0$. 

 We will show that $\RG_c(\A^1\times\Gm, \beta^*\cL_{\psi})=0$, where
the map $\beta:\A^1\times\Gm\to\A^1$ sends $(v,b)$ to $bv^2$. Let $\tilde\beta: \A^1\times\Gm\to\A^1$ be the map that sends $(v,b)$ to $bv$. For the projection
$\pr_1: \A^1\times\Gm\to\A^1$ we have
$$
\pr_{1!}\tilde\beta^*\cL_{\psi}\;\iso\; j_*\Qlb[-1],
$$
where $j:\Gm\to \A^1$ is the open immersion (\cite{Lau}, Lemma~2.3). Let
$\pi:\A^1\to\A^1$ send $v$ to $v^2$. From the diagram
$$
\begin{array}{cccc}
\A^1\times\Gm & \toup{\pi\times\id} & \A^1\times\Gm & \toup{\tilde\beta} \A^1\\
\downarrow\lefteqn{\scriptstyle \pr_1} && \downarrow\lefteqn{\scriptstyle \pr_1} \\
\A^1 & \toup{\pi} & \A^1
\end{array}
$$
we learn that 
$$
\pr_{1!}\beta^*\cL_{\psi}\;\iso\; \pi^*\pr_{1!}\tilde\beta^*\cL_{\psi}
$$

 It suffices to show that $\RG_c(\A^1, \pi^*j_*\Qlb)=0$.
Recall that $\pi_!\Qlb\;\iso\;\Qlb\oplus \cL_1$, where $\cL_1$ is the local system on
$\Gm$ extended by zero to $\A^1$, which corresponds to the Galois covering
$\pi:\Gm\to\Gm$. We get 
$$
\RG_c(\A^1, \pi^*j_*\Qlb)\;\iso\;\RG_c(\A^1, \pi_!\Qlb\otimes j_*\Qlb)=0,
$$ 
because $\RG_c(\Gm,\cL_1)=0$ and $\RG_c(\A^1, j_*\Qlb)=0$.

\medskip
\Step 2 For any $d$ and $i<d$ choose a 
decomposition of $V$ as a direct sum $V=W\oplus V_1\oplus\ldots\oplus
V_{d-i}$, where $\dim V_j=1$ and $\dim W=i$. Then $Q_0(V_1)\times\ldots\times
Q_0(V_{d-i})\subset Q_i(V)$.
The restriction of $S_{\psi}$ to $Q_0(V_1)\times\ldots\times
Q_0(V_{d-i})$ is nonconstant by Step
1 combined with 4).
\end{Prf}

\medskip

\begin{Pp} 
\label{Pp_linear_algebra2}
A choice of a square root $i=\sqrt{-1}\in k$ yields for any $j$ an isomorphism
$$
S_{\psi}\otimes S_{\psi}\mid_{Q_j(V)}\,\iso\,
\Qlb[1](\frac{1}{2})^{\otimes -2j+d(d+1)}
$$
\end{Pp}
\begin{Prf}
Let $\beta_2: V\times V\times \Sym^2 V^*\to\A^1$ be the map sending $(v,u,b)$ to
$\<v,bv\>+\<u,bu\>$. Let $p_3: V\times V\times\Sym^2 V^*\to \Sym^2 V^*$ be the
projection. One checks that 
$$
S_{\psi}\otimes S_{\psi}\;\iso\;
p_{3!}\beta_2^*\cL_{\psi}\otimes\Qlb[1](\frac{1}{2})^{\otimes 2d+d(d+1)}
$$
The change of variables 
$$
\left\{
\begin{array}{c}
x=v+iu\\
y=v-iu
\end{array}
\right.
$$
makes $\beta_2$ to be the map sending $(x,y,b)$ to $\<x,by\>$. Sommate first over $x$ with
$y$ fixed, the assertion follows.
\end{Prf}

\begin{Pp} 
\label{Pp_linear_algebra3}
The $*$-restriction $\Four_{\psi}(\cL_i)\mid_{Q_j(V)}$ vanishes if and only if  
$j\ne i+d\mod 2$. In other words, if $i=d \mod 2$ then $\Four_{\psi}(\cL_i)$ has nontrivial fibres at 
$\cup_{j \;{\rm\scriptstyle even}} Q_j(V)$. If $i\ne d\mod 2$ then $\Four_{\psi}(\cL_i)$ has nontrivial fibres at $\cup_{j \;{\rm\scriptstyle
odd}} Q_j(V)$.

\smallskip

 In particular, $\Four_{\psi}(\cL_i)[d](\frac{d}{2})=S_{\psi,g}$ for $i=d \mod 2$ and
$\Four_{\psi}(\cL_i)[d](\frac{d}{2})=S_{\psi,s}$ for $i\ne d\mod 2$.
\end{Pp}

\smallskip
\begin{Prf}
For $d=1$ it is clear. Assume it is true for $d-1$. 

  The complex $\Four_{\psi}(\cL_j)$ is $\GL(V)$-equivariant, and 
$\GL(V)$ acts transitively on $Q_i(V)$. So, for each $i$ exactly one of two sheaves
$\Four_{\psi}(\cL_0)\mid_{Q_i(V)}$ or $\Four_{\psi}(\cL_1)\mid_{Q_i(V)}$ vanishes, and the other is a rank one (shifted) local system.

Write $V=V_1\oplus V_2$, where $\dim V_1=d-1$ and $\dim V_2=1$. Consider the natural map
$s: \Sym^2 V\to \Sym^2V_1\times\Sym^2 V_2$. We have 
$$
s_!(\cL_0)\,\iso\, (\cL_0\boxtimes\cL_0) \oplus (\cL_1\boxtimes\cL_1)
$$
and
$$
s_!(\cL_1)\,\iso\, (\cL_0\boxtimes\cL_1) \oplus (\cL_1\boxtimes\cL_0),
$$
where on the right hand side $\cL_i$ are those for $V_1$ and $V_2$. 

 Clearly, $Q_{i-1}(V_1)\times Q_1(V_2)\hook{} Q_i(V)$ and 
$Q_i(V_1)\times Q_0(V_2)\hook{} Q_i(V)$. Consider 
\begin{equation}
\label{restriction_1}
\Four_{\psi}(\cL_0)\mid_{Q_i(V_1)\times Q_0(V_2)}\;\iso\;
h^*(\Four_{\psi}(\cL_0)\boxtimes\Four_{\psi}(\cL_0))\oplus h^*(\Four_{\psi}(\cL_1)\boxtimes\Four_{\psi}(\cL_1)),
\end{equation}
where $h: Q_i(V_1)\times Q_0(V_2)\hook{} \Sym^2 V_1^*\times\Sym^2 V_2^*$. This
isomorphism holds up to a shift and a twist.

  If $i= d\mod 2$ then $h^*(\Four_{\psi}(\cL_1)\boxtimes\Four_{\psi}(\cL_1))$ is non zero by induction hypothesis, so the LHS of (\ref{restriction_1}) does not vanish, hence $\Four_{\psi}(\cL_0)\mid_{Q_i(V)}$ does not vanish either.

 If $i\ne d\mod 2$ then the RHS of (\ref{restriction_1}) vanishes by induction hypothesis,
so the LHS also vanishes. Thus,  $\Four_{\psi}(\cL_0)\mid_{Q_i(V)}$ vanishes.
\end{Prf}

\bigskip\noindent
4.2 Assume $d\ge 1$. Let $Y(V)$ be the moduli scheme of pairs: a one dimensional
subspace $V_0\subset V$ and $b\in\Sym^2(V/V_0)^*$. The projection
$Y(V)\to \Gr(1,V)$ is a vector bundle, where $\Gr(1,V)$ denotes the Grassmanian of
one-dimensional subspaces in $V$. 
Let $\alpha: Y(V)\to \Sym^2 V^*$ be the map
sending the above point to the composition
$$
V\to V/V_0\toup{b} (V/V_0)^*\hook{} V^*
$$
Clearly, $\alpha$ factors through $Q_{\ge 1}(V)\hook{} \Sym^2 V^*$. Note that $Y(V)$
is smooth.

\begin{Pp} The map $\alpha: Y(V)\to Q_{\ge 1}(V)$ is proper surjective and semi-small.
\end{Pp}
\begin{Prf} Stratify $Q_{\ge 1}(V)$ by $Q_i(V)$ for $i\ge 1$. The fibre of $\alpha$ 
over a point $b\in Q_i(V)$ is the projective space of 1-subspaces in $V'$, where $V'$
is the kernel of $b$. So, $\dim\alpha^{-1}(b)=i-1$ and $\dim
Q_i(V)=\frac{1}{2}(d-i)(d+1+i)$. We get
$$
2\dim\alpha^{-1}(b)\le \codim_{Q_{\ge 1}(V)} Q_i(V),
$$
the equality holds only for $i=1,2$. 
\end{Prf}

\medskip\noindent
4.3  {\scshape Relative version\ } Let now $S$ be a smooth scheme, $V\to S$ be a vector
bundle of rank $d$. Define $S_{\psi}\in\D(\Sym^2 V^*)$ by 
(\ref{def_S_psi}), so $S_{\psi}$ is a shifted perverse sheaf.  

 As above, $\Sym^2 V^*$ is stratified by locally closed subschemes $Q_i(V)$, they are
equipped with morphisms $Q_i(V)\to \Gr(i,V)$ over $S$. 

 We also have the $S_2$-coverings $\Cov(Q_i(V))\to Q_i(V)$. For an $S$-scheme $S'$, the
$S'$-points of $\Cov(Q_i(V))$ are collections: a rank $i$ subbundle $V_0\subset V\mid_{S'}$,
an isomorphism $b: V/V_0\to (V/V_0)^*$ of $\cO_{S'}$-modules with $b^*=b$, and a compatible
trivialization $\det (V/V_0)\,\iso\,\cO_{S'}$.

 Propositions~\ref{Pp_linear_algebra1}, \ref{Pp_linear_algebra2} and
\ref{Pp_linear_algebra3} hold in relative situation (one only changes a shift and a twist
in 3) of Proposition~\ref{Pp_linear_algebra1}). 

\medskip\noindent
4.4  {\scshape Finite-dimensional theta-sheaf}  \   This subsection is not used in the proofs and may be skipped.  

 Let $M$ be a symplectic $k$-vector space of dimension $2d$. Write $\cL(M)$ for the scheme of lagrangian subspaces of $M$. Set $Y=\cL(M)\times\cL(M)$. Consider the line bundle $\cA$ on $Y$ with fibre $(\det L_1)\otimes(\det L_2)$ over $(L_1,L_2)\in Y$. We view it as $\ZZ/2\ZZ$-graded purely of degree zero. Let $\tilde Y$ denote the stack of square roots of $\cA$. The $\mu_2$-gerbe $\tilde Y\to Y$ is nontrivial. The group $\Sp(M)$ acts naturally on $Y$, and $\cA$ is $\Sp(M)$-equivariant, so $\Sp(M)$ acts also on $\tilde Y$. 
  
   We are going to construct a $\Sp(M)$-equivariant perverse sheaf $S_M$ on $\tilde Y$ such that $-1\in\mu_2$ acts on $S_M$ as $-1$. 
  
   The $\Sp(M)$-orbits on $Y$ are indexed by $i=0,\ldots,d$. The orbit $Y_i$ is given by $\dim(L_1\cap L_2)=i$.  

\begin{Lm} The restriction of $\cA$ to each $Y_i$ admits a canonical $\Sp(M)$-equivariant square root. 
\end{Lm}
\begin{Prf} For $L_1,L_2\in\cL(M)$ let $(L_1\cap L_2)^{\perp}\subset M$ denotes the orthogonal complement to $L_1\cap L_2$. The symplectic form on $(L_1\cap L_2)^{\perp}/(L_1\cap L_2)$ induces an isomorphism
$L_2/(L_1\cap L_2)\,\iso\, (L_1/L_1\cap L_2)^*$. This yields a $\ZZ/2\ZZ$-graded isomorphism
$(\det L_1)\otimes (\det L_2)\,\iso\, \det(L_1\cap L_2)^{\otimes 2}$. By 3.1.2, we are done.
\end{Prf}

\medskip

  Let $W$ denote the nontrivial local system of rank one on $B(\mu_2)$ corresponding to the covering $\Spec k\to B(\mu_2)$. Let $\tilde Y_i$ denote the restriction of the gerbe
$\tilde Y\to Y$ to $Y_i$, so $\tilde Y_i\,\iso\, Y_i\times B(\mu_2)$ canonically. 

\begin{Def} Let $S_{M, g}$ (resp., $S_{M,s}$) denote the Goresky-MacPherson extension of 
$$
(\Qlb\boxtimes W)[\dim Y](\frac{\dim Y}{2})
$$ 
from $\tilde Y_0$ to $\tilde Y$ (resp., of $(\Qlb\boxtimes W)[\dim Y-1](\frac{\dim Y-1}{2})$ from $\tilde Y_1$ to $\tilde Y$). Set $S_M=S_{M,g}\oplus S_{M,s}$. 
\end{Def}

Denote by $\cY$ the stack quotient $Y/\Sp(M)$. Write $\tilde\cY\to\cY$ for the corresponding gerbe of square roots of $\cA$. We may view $S_M$ as a perverse sheaf on $\tilde \cY$.   

 Fix a lagrangian subspace $V\subset M$, let $P_V\subset \Sp(M)$ be the Seigel parabolic subgroup preserving $V$.  We have canonical isomorphisms of stacks 
$$
\cY\,\iso\, \cL(M)/P_V\,\iso\, P_V\backslash \Sp(M)/P_V
$$  
One may view $\cA$ as a line bundle on $\cL(M)/P_V$ with fibre $(\det V)\otimes(\det L)$. 

 Fix a splitting $V^*\to M$ of $0\to V\to M\to V^*\to 0$. Denote by $P^-_V\subset \Sp(M)$ the Seigel parabolic subgroup preserving $V^*\subset M$. Let $Z\subset \cL(M)$ be the open $P^-_V$-orbit, that is 
$$
Z=\{L\in\cL(M)\mid L\cap V^*=0\}
$$ 
The map $\Sym^2 V^*\to Z$ sending $b:V\to V^*$ to $L=\{v+bv\in M\mid v\in V\}$ is an isomorphism commuting with the action of $\GL(V)$. Denote by $\cZ$ the stack quotient $Z/\GL(V)$. View $S_{\psi}$ introduced in Sect. 4.1 as a perverse sheaf on $\cZ$. 

Denote by $\nu$ the composition (of an open immersion followed by a smooth map)
$$
\cZ\hook{} \cL(M)/\GL(V)\to \cL(M)/P_V=\cY
$$ 
The map $\nu:\cZ\to \cY$ is smooth, surjective and representable. It factors naturally as
$\cZ\toup{\tilde\nu} \tilde\cY\to\cY$. 

\begin{Pp} 
\label{Pp_linear_algebra_last}
There are isomorphisms of perverse sheaves on $\cZ$
$$
\tilde\nu^* S_{M,g}\otimes\Qlb[1](\frac{1}{2})^{\otimes \dim\cZ-\dim\cY}\;\iso\; S_{\psi,g}
$$
and 
$$
\tilde\nu^* S_{M,s}\otimes\Qlb[1](\frac{1}{2})^{\otimes \dim\cZ-\dim\cY}\;\iso\; S_{\psi,s}
$$
(Once $i=\sqrt{-1}\in k$ is fixed, such isomorphisms are well defined up to multiplication by $\pm1$).
\end{Pp}
\begin{Prf} The stack $\cZ$ is stratified by $\cZ_i=Q_i(V)/\GL(V)$, the quotient being taken in stack sense. Let $\cY_i$ denote the stack quotient $Y_i/\Sp(M)$. Note that $\cZ_i$ identifies with $\cZ\times_{\cY}\cY_i$ for $i=0,\ldots,d$.

 Let $\tilde\cY_i$ denote the restriction of the gerbe $\tilde\cY$ to $\cY_i$, so $\tilde\cY_i\,\iso\, \cY_i\times B(\mu_2)$ canonically. Remind the covering $\Cov(Q_i(V))\to Q_i(V)$ from Sect.~4.1.  It is $\GL(V)$-equivariant, so the stack quotient $\Cov(\cZ_i)=\Cov(Q_i(V))/\GL(V)$ is a two-sheeted covering of $\cZ_i$. For each $i$ we have a cartesian square
$$
\begin{array}{ccc}
\Cov(\cZ_i) & \to & \cY_i\\
\downarrow && \downarrow\\
\cZ_i & \toup{\tilde\nu} & \tilde\cY_i
\end{array}
$$
Our assertion follows now from Proposition~\ref{Pp_linear_algebra1}.
\end{Prf}

\medskip

\begin{Rem} Write $_MY$ (resp., $_M\tilde\cY$)  to express the dependence on $M$. If $M,M'$ are two symplectic spaces over $k$ of dimensions $d,d'$, consider the map 
$\tau_{M,M'}: {_MY}\times{_{M'}Y}\to {_{M\oplus M'}Y}$ sending $(L_1,L_2), (L'_1,L'_2)$ to
$(L_1\oplus L'_1, L_2\oplus L'_2)$. It yields a map 
$$
\tilde\tau_{M,M'}:  {_M\tilde\cY}\times{_{M'}\tilde\cY}\to {_{M\oplus M'}\tilde\cY}
$$
From 4) of Proposition~\ref{Pp_linear_algebra1} it follows that $\tilde\tau_{M,M'}^* S_{M\oplus M'}\;\iso\; S_M\boxtimes S_{M'}[2dd'](dd')$ canonically. 
\end{Rem}

\bigskip
\centerline{\scshape 5. Fourier coefficients of $\Aut$ for Siegel parabolic}

\bigskip\noindent
5.1 Write $\Bun_P=\Bun_{P_n}$. So,  $\Bun_P$ classifies $L\in\Bun_n$ together with an exact sequence $0\to \Sym^2 L\to
?\to \Omega\to 0$ on $X$. It induces an exact sequence 
\begin{equation}
\label{sequence_M}
0\to L\to M\to L^*\otimes\Omega\to 0,
\end{equation}
The map $\nu_n: \Bun_P\to\Bun_G$ is also denoted $\nu$. 

\begin{Lm} 
\label{Lm_tilde_nu}
The map $\nu: \Bun_P\to\Bun_G$ factors as the composition
$\Bun_P\toup{\tilde\nu} \Bunt_G\toup{\gr} \Bun_G$.
\end{Lm} 
\begin{Prf}
The sequence (\ref{sequence_M}) yields a $\ZZ/2\ZZ$-graded isomorphism
\begin{equation}
\label{iso_gerbe}
\det\RG(X,M)\,\iso\, \det\RG(X,L)\otimes\det\RG(X, L^*\otimes\Omega)\,\iso\,
\det\RG(X, L^*\otimes\Omega)^2
\end{equation}
Define $\tilde\nu$ by letting $\cB=\det\RG(X, L^*\otimes\Omega)$ together with $\cB^2\;\iso\;\cA$ given by (\ref{iso_gerbe}). By 3.1.2, $\tilde\nu$ is well-defined. 
\end{Prf}

\medskip\smallskip

 Let $^0\Bun_P\subset\Bun_P$ be the open substack given by $\H^0(X, \Sym^2 L)=0$. One checks that both $\nu: {^0\Bun_P}\to\Bun_G$ and $\tilde\nu: {^0\Bun_P}\to\Bunt_G$ are smooth.

\begin{Lm} 
\label{Lm_5_surjectivity}
The map $\nu: {^0\Bun_P}\to\Bun_G$ is surjective, so $\tilde\nu: {^0\Bun_P}\to\Bunt_G$ is
also surjective.
\end{Lm}
\begin{Prf} Let $M$ be a $k$-point of $\Bun_G$. It admits a line subbundle $L_1$ with
$\deg L_1<0$. Let $L_{-1}\subset M$ be the orthogonal complement to $L_1$, so $L_{-1}/L_1\in \Bun_{G_{n-1}}$ naturally. Continuing this procedure for $L_{-1}/L_1$ and so on, we get a flag of isotropic subbundles $L_1\subset\ldots\subset L_n\subset M$. Then $(L_n\subset M)$ is a $k$-point of $^0\Bun_P$.
\end{Prf}

\bigskip\noindent
5.2 {\scshape the sheaf $S_{P,\psi}$ on $\Bun_P$}

\smallskip\noindent
 Write $\Bun_n^d$ (resp., $\Bun_P^d$) for the connected component of the
corresponding stack given by $\deg L=d$. 

 Write $_c\Bun_n\subset \Bun_n$ for the open substack given by
$\H^0(X,L)=0$. Let $\cV\to \Bun_n$ be the stack whose fibre over $L\in\Bun_n$ is
$\Hom(L,\Omega)$. Let $_c\cV\to {_c\Bun_n}$ be the preimage of $_c\Bun_n$, over
$_c\Bun_n^d$ this is a vector bundle of rank $n(g-1)-d$.  

 Let $\cX=\cV\times_{\Bun_n}\Bun_P$ and $p:\cX\to\Bun_P$ be the projection. 
Let $q:\cX\to\A^1$ be the map sending $s\in\H^0(X, L^*\otimes\Omega)$ to the pairing of
$0\to\Sym^2 L\to ?\to\Omega\to 0$ with
$$
s\otimes s\in \H^0(X, (\Sym^2 L^*)\otimes\Omega^2)
$$
\begin{Def} Set $S_{P,\psi}=p_!q^*\cL_{\psi}\otimes\Qlb[1](\frac{1}{2})^{\otimes
\dim\cX}$, where $\dim\cX$ is the dimension of the corresponding connected component
of $\cX$.
\end{Def}

 Let $\cV_2\to \Bun_n$ be the stack whose fibre over $L\in\Bun_n$ is
$\Hom(\Sym^2 L, \Omega^2)$. Let $\pi_2: \cV\to\cV_2$ be the map sending 
$s\in\Hom(L,\Omega)$ to $s\otimes s$. Note that $\pi_2$ is finite, a $S_2$-covering
over the image $\Im\pi_2$ with removed zero section. By definition,
\begin{equation}
\label{second_def_S_P}
S_{P,\psi}\;\iso\; \Four_{\psi}(\pi_{2!}\Qlb)\otimes\Qlb[1](\frac{1}{2})^{\otimes \dim\cV},
\end{equation}
where $\Four_{\psi}: \D(\cV_2)\to \D(\Bun_P)$ is the Fourier transform functor. Note that $S_2$ acts on $S_{P,\psi}$. 

 Let $_c\Bun_P$ denote the preimage of $_c\Bun_n$ in $\Bun_P$. We see that over each connected component of $_c\Bun_P$, $S_{P,\psi}$ is a direct sum of two irreducible perverse sheaves and $\DD(S_{P,\psi})\,\iso\, S_{P,\psi^{-1}}$.  

\smallskip

 Let $\Sym^2{_c\cV}\to {_c\Bun_n}$ denote the symmetric square of the vector bundle
$_c\cV\to {_c\Bun_n}$. Let $\pi: {_c\cV}\to \Sym^2{_c\cV}$ be the map sending $s\in
\Hom(L,\Omega)$ to $s\otimes s$. Then $\pi_2$ decomposes as
$$
_c\cV\toup{\pi} \Sym^2{_c\cV}\toup{f^*}
{_c\cV_2}
$$
Given $L\in\Bun_n$, the transpose to the linear map
$\Sym^2 \H^0(X, L^*\otimes\Omega)\to \Hom(\Sym^2 L, \Omega^2)$ is
$$
\H^1(X, (\Sym^2 L)\otimes\Omega^{-1})\to \Sym^2 \H^1(X,L)
$$ 
It defines a morphism of stacks $f: {_c\Bun_P}\to \Sym^2{_c\cV^*}$ over $_c\Bun_n$. 

 We have the sheaf $S_{\psi}$ on $\Sym^2{_c\cV^*}$ defined in Sect.~4.3. From
(\ref{second_def_S_P}) we conclude that
\begin{equation}
\label{formula_S_P_psi}
S_{P,\psi}\;\iso\; f^*S_{\psi}\otimes\Qlb[1](\frac{1}{2})^{\otimes
\dim\cX-r-\frac{1}{2}r(r+1)}
\end{equation}
canonically over $_c\Bun_P$, where $r$ and $\dim\cX$ are
functions of the corresponding connected component with   
$r(_c\Bun^d_P)=n(g-1)-d$. 

 Denote by $S_{P,\psi,g}$ (resp., by $S_{P,\psi,s}$) the direct summand of
$S_{P,\psi}$ obtained by replacing $S_{\psi}$ by $S_{\psi,g}$ (resp., by $S_{\psi,s}$) in
(\ref{formula_S_P_psi}). Both $S_{P,\psi,g}$
and $S_{P,\psi,s}$ are irreducible perverse sheaves over each connected component of $_c\Bun_P$.

 Note that $^0\Bun_P\subset {_c\Bun_P}$.

\begin{Rem} 
\label{Rem_S_2-action}
Consider $\chi(L)$ as a function of a connected component of $_c\Bun_P$. By Proposition~\ref{Pp_linear_algebra3}, over a given connected component of $_c\Bun_P$, the $S_2$-invariants of $S_{P,\psi}$ are $S_{P,\psi,g}$
for $\chi(L)$ even and $S_{P,\psi,s}$ for $\chi(L)$ odd. 
\end{Rem} 

\medskip\noindent
5.3 \  Recall the stratification of $\Sym^2{_c\cV^*}$ by locally closed substacks
$Q_i(_c\cV)$ and the coverings $\Cov(Q_i(_c\cV))\to
Q_i(_c\cV)$ defined in Sect.~4.3.

 Set $_i\Bun_P=\nu^{-1}(_i\Bun_G)$ and $_{i,c}\Bun_P={_c\Bun_P}\cap {_i\Bun_P}$. 
For a point of $_c\Bun_P$ the exact sequence (\ref{sequence_M}) yields an exact sequence
\begin{equation}
\label{sequence_coh}
0\to\H^0(X,M)\to \H^0(X, L^*\otimes\Omega)\toup{b} \H^1(X,L)\to\H^1(X, M)\to 0
\end{equation}
Thus, we get a commutative diagram
$$
\begin{array}{ccc}
_{i,c}\Bun_P & \hook{} & {_c\Bun_P}\\
\downarrow\lefteqn{\scriptstyle f} && \downarrow\lefteqn{\scriptstyle f}\\
Q_i(_c\cV) & \hook{} & \Sym^2{_c\cV^*}
\end{array}
$$

 Let $_i\rho_P:\Cov(_{i,c}\Bun_P)\to {_{i,c}\Bun_P}$ be the covering obtained from
$\Cov(Q_i(_c\cV))\to Q_i(_c\cV)$ by the base change $f:{_{i,c}\Bun_P}\to Q_i(_c\cV)$. 

\begin{Pp} 
\label{Pp_cartesian_square}
For $i\ge 0$ there is a cartesian square
$$
\begin{array}{ccc}
\Cov(_{i,c}\Bun_P) & \to & \Cov(_i\Bunt_G)\\
\downarrow\lefteqn{\scriptstyle {_i\rho_P}} && \downarrow\lefteqn{\scriptstyle {_i\rho}}\\
_{i,c}\Bun_P & \toup{\tilde\nu} & {_i\Bunt_G}
\end{array}
$$
\end{Pp}
\begin{Prf}
Let $S$ be a scheme. Assume given an $S$-point of $_{i,c}\Bun_P$. It yields locally free
$\cO_S$-modules $V_0=\H^0(X,M)$ and $V=\H^0(X, L^*\otimes\Omega)$ included into an exact
sequence of $\cO_S$-modules (a relative version of (\ref{sequence_coh}))
$$
0\to V_0\to V\toup{b} V^*\to V_0^*\to 0
$$
with $b^*=b$. The $\cO_{S\times X}$-module $L$ together with the morphism of $\cO_S$-modules
$V\toup{b} V^*$ defines the corresponding $S$-point of $Q_i(_c\cV)$. 

  We have an isomorphism of $\cO_S$-modules $\cB=\det\RG(X, L^*\otimes\Omega)\,\iso\,
\det V$, because $\H^0(X,L)=0$. We also have an isomorphism
of $\cO_S$-modules $t: \cB^2\,\iso\,\det\RG(X,M)\,\iso(\det V_0)^2$ given by
(\ref{iso_gerbe}).

 A lifting of the corresponding $S$-point of $_i\Bunt_G$ to $\Cov(_i\Bunt_G)$ is an
isomorphism of $\cO_S$-modules $\cB\,\iso\, \det V_0$ whose square is $t$. 
The corresponding category is the category of $S$-points of $\Cov(_{i,c}\Bun_P)$.
\end{Prf} 

\smallskip

\begin{Pp} 
\label{Pp_up_to_minus_one}
There are isomorphisms of perverse sheaves on $^0\Bun_P$
$$
\tilde\nu^*\Aut_g\otimes\Qlb[1](\frac{1}{2})^{\otimes \dim\Bun_P-d_G}\;\iso\;
S_{P,\psi,g}
$$
and
$$
\tilde\nu^*\Aut_s\otimes\Qlb[1](\frac{1}{2})^{\otimes \dim\Bun_P-d_G}\;\iso\;
S_{P,\psi,s}
$$
Here $\dim\Bun_P$ denotes the dimension of the corresponding connected component of $\Bun_P$. (Once $\sqrt{-1}\in k$ is fixed, the above isomorphisms are well-defined up to a sign).
\end{Pp}
\begin{Prf} \  Recall that $S_{P,\psi,g}$ and $S_{P,\psi,s}$ are irreducible
perverse sheaves over each connected component of $_c\Bun_P$. By relative
version of Proposition~\ref{Pp_linear_algebra1},
$S_{P,\psi,g}$ over $_{0,c}\Bun_P$ (resp., $S_{P,\psi,s}$ over $_{1,c}\Bun_P$) is a
nonconstant local system of rank one corresponding to the covering
$\Cov(_{0,c}\Bun_P)\to{_{0,c}\Bun_P}$ (resp., $\Cov(_{1,c}\Bun_P)\to{_{1,c}\Bun_P}$). 
Moreover, for any $i$
$$
(S_{P,\psi}\otimes S_{P,\psi})\mid_{_{i,c}\Bun_P}
\,\iso\, \Qlb[2](1)^{\otimes \dim\Bun_P-i}
$$ 
by Proposition~\ref{Pp_linear_algebra2} (this requires a choice of $\sqrt{-1}\in k$). 

 By Proposition~\ref{Pp_cartesian_square}, for each $i$ we get isomorphisms
$$
\tilde\nu^*(_i\Aut)\mid_{_{i,c}\Bun_P}\,\iso\, \Hom_{S_2}(\sign,
(_i\rho_P)_!\Qlb))
$$
In particular,
$$
\tilde\nu^*(_i\Aut\otimes{_i\Aut})\mid_{_{i,c}\Bun_P}\;\iso\; \Qlb
$$

 Set $^0_i\Bun_P={^0\Bun_P}\cap{_i\Bun_P}$. By construction, $S_{P,\psi,s}$ is perverse over $_{1,c}\Bun_P$, hence also over
$^0_1\Bun_P$. Since $^0_1\Bun_P\to {_1\Bun_G}$ is smooth and
surjective, Propositions~\ref{Pp_linear_algebra1} and \ref{Pp_cartesian_square} imply that $_1\Aut\otimes\Qlb[1](\frac{1}{2})^{\otimes d_G-1}$ is 
perverse on $_1\Bun_G$. So, Defenition~Ê\ref{Def_1} makes sense.

 For each connected component $^0\Bun^d_P$ of $^0\Bun_P$ the
stack $^0\Bun^d_P\cap\, {_i\Bun_P}$ is non empty for $i=0,1$. Since $\tilde\nu:
{^0\Bun_P}\to \Bunt_G$ is smooth, our assertion follows.
\end{Prf}

\bigskip\smallskip\noindent
\begin{Prf}\select{of Theorem~\ref{Th_1}\ }
 For each $d$ the map $\tilde\nu: {^0\Bun^d_P}\to \Bunt_G$ is smooth with connected
fibres, and $\tilde\nu: {^0\Bun_P}\to \Bunt_G$ is surjective. So, by  Proposition~\ref{Pp_up_to_minus_one} it suffices to construct isomorphisms
$$
S_{P,\psi}\mid_{^0_i\Bun_P}\,\iso\, \tilde\nu^*(_i\Aut)\otimes\Qlb[1](\frac{1}{2})^{\otimes \dim\Bun_P -i}
$$
over $^0_i\Bun_P$. We have them by Proposition~\ref{Pp_cartesian_square} combined with relative version of Proposition~\ref{Pp_linear_algebra1}.  Proposition~\ref{Pp_linear_algebra3}
implies the second part of the theorem.
\end{Prf}

\bigskip\noindent
\begin{Rem} From Theorem~\ref{Th_1} it follows that
$\tilde\nu^*\Aut\otimes\Qlb[1](\frac{1}{2})^{\otimes \dim\Bun_P-d_G}$ equals
$S_{P,\psi}$ in the Grothendieck group $K(_c\Bun_P)$ over $_c\Bun_P$, which is
bigger than $^0\Bun_P$. We expect that actually the isomorphisms of
Proposition~\ref{Pp_up_to_minus_one} hold over $_c\Bun_P$.
\end{Rem}

\bigskip\medskip
\centerline{\scshape 6. Constant terms of $\Aut$ for maximal parabolics}

\bigskip\noindent
6.1  Recall the smooth map $\eta_k: \Bun_{P_k}\to\Bun_{Q_k}$ (cf. Sect.~3.3). Under
each of the two projections 
$
\Bun_{P_k}\times_{\Bun_{Q_k}}\Bun_{P_k}\to\Bun_{P_k}
$
the stack $\Bun_{P_k}\times_{\Bun_{Q_k}}\Bun_{P_k}$ identifies with the one classifying $(L_1\subset L_{-1}\subset M)\in\Bun_{P_k}$ together with an exact sequence $0\to\Sym^2 L_1\to ?\to\Omega\to 0$, the projection being the forgetful map.

 Let $\nu_{k,n}:\Bun_{P_{k,n}}\to\Bun_P$ be the stack classifying $(0\to\Sym^2 L\to ?\to\Omega\to 0)\in\Bun_P$ together with a subbundle
$L_1\subset L$ with $L_1\in\Bun_k$.

\begin{Lm} 
The map $\eta_k: \Bun_{P_k}\to\Bun_{Q_k}$ is surjective.
\end{Lm}
\begin{Prf}
Consider a $k$-point of $\Bun_{Q_k}$ given by a flag $(L_1\subset L_{-1})$ of vector bundles on $X$ with $L_{-1}/L_1\in\Bun_{G_{n-k}}$. Let show that the fibre of $\eta_k$ over it is nonempty.

 Pick a lagrangian subbundle $\cB\subset L_{-1}/L_1$ such that $\H^1(X, \cB^*\otimes L_1)=0$, it always exists. Let $L\subset L_{-1}$ be the preimage of $\cB$ under $L_{-1}\to L_{-1}/L_1$. The exact sequence $0\to L_1\to L\to\cB\to 0$ splits, we fix a splitting  $L\,\iso\, L_1\oplus\cB$. Then our $k$-point of $\Bun_{Q_k}$ becomes the data of two exact sequences
$$
0\to\Sym^2\cB\to ?\to\Omega\to 0
$$
and 
$$
0\to L_1\to ?\to \cB^*\otimes\Omega\to 0,
$$ 

 Pick any exact sequence $0\to\Sym^2 L_1\to ?\to\Omega\to 0$ and summate it with the above two. The result is an exact sequence $0\to\Sym^2 L\to ?\to\Omega\to 0$, the corresponding $P_{k,n}$-torsor induces a $P_k$-torsor lying in the fibre under consideration. 
\end{Prf}

\medskip

 Set $\Bun_{Q_{k,n}}=\Bun_{P(G_{n-k})}\times_{\Bun_{G_{n-k}}}\Bun_{Q_k}$, where $P(G_{n-k})\subset G_{n-k}$ is the Siegel parabolic. So, $\Bun_{Q_{k,n}}$ classifies 
a point $0\to L_1\to L_{-1}\to L_{-1}/L_1\to 0$ of $\Bun_{Q_k}$ together with a lagrangian subbundle $L/L_1\subset L_{-1}/L_1$. Consider the diagram
$$
\begin{array}{cc}
\Bun_P\getsup{\nu_{k,n}}\Bun_{P_{k,n}}\toup{\eta_{k,n}} \!\! \!\!& \Bun_{Q_{k,n}}\\
& \downarrow\lefteqn{\scriptstyle r_k}\\
& \Bun_{P(G_{n-k})},
\end{array}
$$
where $r_k$ and $\eta_{k,n}$ denote the projections.

 Write $S_{P(G_n),\psi}$ to express the dependence of $S_{P,\psi}$ on $n$. Note that $\Bun_{P(G_0)}=\Spec k$, $S_{P(G_0),\psi,g}=\Qlb$ and $S_{P(G_0),\psi,s}=0$.

\begin{Pp}
\label{Pp_7_commuting_with_S_2}
 We have a canonicall isomorphism commuting with $S_2$-action
$$
(\eta_{k,n})_!\nu_{k,n}^*S_{P,\psi}\;\iso\; r_k^*S_{P(G_{n-k}),\psi}[a](\frac{a}{2}),
$$ 
where $a\in\ZZ$ is the function of a connected component of $\Bun_{Q_{k,n}}$ given by 
$$
a=\dim\Bun_n-\dim\Bun_{n-k} - \chi(L_1) + \chi(\Omega^{-1}\otimes\Sym^2 L_1)-\chi(\Omega^{-1}\otimes L_1\otimes(L/L_1))
$$
\end{Pp}
\begin{Prf}
Consider the map 
$$
\cX\times_{\Bun_P}\Bun_{P_{k,n}}=\cV\times_{\Bun_n}\Bun_{P_{k,n}}\toup{\id\times\eta_{k,n}} \cV\times_{\Bun_n}\Bun_{Q_{k,n}}
$$ 
Write $\A^1\getsup{q_n}\cX_{G_n}\toup{p_n} \Bun_{P(G_n)}$ to express the dependence on $n$ of the diagram $\A^1\getsup{q}\cX\toup{p}\Bun_P$ introduced in Sect.~5.2.
 
 Denote temporary by  $i: \cX_{G_{n-k}}\times_{\Bun_{P(G_{n-k})}}\Bun_{Q_{k,n}}\hook{} \cV\times_{\Bun_n}\Bun_{Q_{k,n}}$ the closed embedding given by the condition that $s\in\Hom(L,\Omega)$ lies in $\Hom(L/L_1,\Omega)$. 

 Set $a_0=-\chi(\Omega^{-1}\otimes\Sym^2L_1)$ viewed as a function of a connected component of $\Bun_{Q_{k,n}}$. Let us establish a canonical isomorphism
\begin{equation}
\label{iso_temporary_Pp7}
(\id\times\eta_{k,n})_!(q^*\cL_{\psi}\boxtimes\Qlb)\;\iso\; i_!(q_{n-k}^*\cL_{\psi}\boxtimes\Qlb)[-2a_0](-a_0)
\end{equation}
Consider a $k$-point of $\cV\times_{\Bun_n}\Bun_{P_{k,n}}$ given by 
$(L_1\subset L\subset L_{-1}\subset M)$ and $s\in\Hom(L,\Omega)$. The fibre, say $Y$, of $\id\times\eta_{k,n}$ over its image in $\cV\times_{\Bun_n}\Bun_{Q_{k,n}}$ identifies with the stack of exact sequences 
\begin{equation}
\label{sequence_Pp7_proof}
0\to \Sym^2 L_1\to ?\to\Omega\to 0
\end{equation}
 on $X$. The restriction of $q^*\cL_{\psi}\boxtimes\Qlb$ to Y is (up to a tensoring by a 1-dimensional vector space) the restriction of $\cL_{\psi}$ under the map $Y\to \A^1$
pairing $\Sym^2 L_1\hook{} \Sym^2 L\toup{s\otimes s} \Omega^2$ with (\ref{sequence_Pp7_proof}). 
 
 So, the fibre of the LHS of (\ref{iso_temporary_Pp7}) vanishes unless $s\in\Hom(L/L_1,\Omega)$. The isomorphism (\ref{iso_temporary_Pp7}) follows, here $a_0=\dim Y$. 

 For the projection $\pr:\cV\times_{\Bun_n}\Bun_{Q_{k,n}}\to\Bun_{Q_{k,n}}$ we get
$$
\pr_! (\id\times\eta_{k,n})_!(q^*\cL_{\psi}\boxtimes\Qlb)\otimes\Qlb[1] (\frac{1}{2})^{\otimes\dim\cX} \,\iso\, (\eta_{k,n})_!\nu_{k,n}^*S_{P,\psi}
$$
Our assertion follows, because $a=\dim\cX_{G_n}-\dim\cX_{G_{n-k}}-2a_0$.
\end{Prf}

\bigskip

\begin{Prf}\select{of Theorem~\ref{Th_2}} \ \ 
We have the diagram
$$
\begin{array}{ccc}
\Bun_P & \toup{\tilde\nu} & \Bunt_G\\
\uparrow\lefteqn{\scriptstyle \nu_{k,n}} && \uparrow\lefteqn{\scriptstyle
\tilde\nu_k}\\
\Bun_{P_{k,n}} & \to & \Bunt_{G_{n-k}}\times_{\Bun_{G_{n-k}}}\Bun_{P_k}\\
\downarrow\lefteqn{\scriptstyle \eta_{k,n}} && \downarrow\lefteqn{\scriptstyle
\id\times\eta_k}\\
\Bun_{Q_{k,n}} & \toup{\tilde\nu\times\id} &
\Bunt_{G_{n-k}}\times_{\Bun_{G_{n-k}}}\Bun_{Q_k}\\
\downarrow\lefteqn{\scriptstyle r_k} && \downarrow\\
\Bun_{P(G_{n-k})} & \toup{\tilde\nu} & \Bunt_{G_{n-k}},
\end{array}
$$
where the middle square is cartesian. So,
$$
(\tilde\nu\times\id)^*(\id\times\eta_k)_!\tilde\nu_k^*\Aut\;\iso\;
(\eta_{k,n})_!\nu_{k,n}^*\tilde\nu^*\Aut
$$

 Let $^0\Bun_{Q_{k,n}}\subset\Bun_{Q_{k,n}}$ be the open substack given by three
conditions: $\H^0(X,\Sym^2 L_1)=0$, $\H^0(X,L_1\otimes L/L_1)=0$, and
$\H^0(X,\Sym^2(L/L_1))=0$. As in Lemma~\ref{Lm_5_surjectivity}, one checks that
\begin{equation}
\label{map_inside_Th2}
^0\Bun_{Q_{k,n}}\toup{\tilde\nu\times\id} \;
\Bunt_{G_{n-k}}\times_{\Bun_{G_{n-k}}}{^0\Bun_{Q_k}}
\end{equation}
is smooth and surjective. Since $\eta_{k,n}^{-1}(^0\Bun_{Q_{k,n}})\subset
\nu_{k,n}^{-1}(^0\Bun_P)$, from Propositions~\ref{Pp_up_to_minus_one} and
\ref{Pp_7_commuting_with_S_2} we get
\begin{equation}
\label{iso_inside_Th2}
(\tilde\nu\times\id)^*(\id\times\eta_k)_!\tilde\nu_k^*\Aut\;\iso\;
r_k^* S_{P(G_{n-k}),\psi}\otimes\Qlb[1](\frac{1}{2})^{\otimes d_G-\dim\Bun_P+a}
\end{equation}
over $^0\Bun_{Q_{k,n}}$. The restriction of $r_k$ to $^0\Bun_{Q_{k,n}}$ factors as
$$
^0\Bun_{Q_{k,n}}\toup{r_k} {^0\Bun_{P(G_{n-k})}}\hook{}\Bun_{P(G_{n-k})}
$$
So, by Proposition~\ref{Pp_up_to_minus_one} applied to $G_{n-k}$, the RHS of
(\ref{iso_inside_Th2}) identifies with
$$
(\tilde\nu\times\id)^*(\Aut\boxtimes\Qlb)\otimes\Qlb[1](\frac{1}{2})^{\otimes
d_G-\dim\Bun_P+a+\dim\Bun_{P(G_{n-k})}-d_{G_{n-k}}}
$$
We have
$b(L_1)=d_G-\dim\Bun_P+a+\dim\Bun_{P(G_{n-k})}-d_{G_{n-k}}$. 
Since $\Bun_{Q_k}\to\Bun_{G_{n-k}}$ is smooth, $\Aut\boxtimes\,\Qlb$ is a shifted
perverse sheaf on $\Bunt_{G_{n-k}}\times_{\Bun_{G_{n-k}}}\Bun_{Q_k}$.

 Since the restriction of the map (\ref{map_inside_Th2}) to each connected component of  $^0\Bun_{Q_{k,n}}$ has connected fibres, we get the desired isomorphism.

 The second assertion follows from Remark~\ref{Rem_S_2-action} combined with Proposition~\ref{Pp_7_commuting_with_S_2}.
\end{Prf}  

\bigskip\medskip
\centerline{\scshape 7. Towards geometric $\theta$-lifting}

\bigskip\noindent
This section is not used in the proofs and may be skipped. Let $\tau_{n,m}:\Bun_{G_n}\times\Bun_{\SO_m}\to \Bun_{G_{nm}}$ be the following map. Given $\SO_m$-torsor $\cF_W$, let $W$ denote the vector bundle induced from it via the standard representation of $\SO_m$. Given in addition $M\in\Bun_{G_n}$ we get naturally a symplectic form $\wedge^2 (M\otimes W)\to\Omega$. The map $\tau_{n,m}$ sends $(M,W)$ to $M\otimes W$.

Let $\cA_{\SO_m}$ denote the (naturally $\ZZ/2\ZZ$-graded) line bundle on $\Bun_{\SO_m}$, whose fibre at $\cF_W$ is $\det\RG(X,W)$. Write $\cA_{G_n}$ to express the dependence on $n$ of the determinant of cohomology on $\Bun_{G_n}$. 

\begin{Lm} For $m\ge 3$ we have a $\ZZ/2\ZZ$-graded canonical isomorphism 
over $\Bun_{G_n}\times\Bun_{\SO_m}$
$$
\tau^*_{n,m}\cA_{G_{nm}}\,\iso\, (\cA^m_{G_n}\boxtimes \cA^{2n}_{\SO_m})\otimes\det\RG(X,\cO)^{\otimes -2nm}
$$ 
\end{Lm}
\begin{Prf}
\Step 1 Let us show that for any $M\in\Bun_{G_n}$, $V\in\Bun_{\SL_2}$ we have canonically
$$
\det\RG(X, M\otimes V)\,\iso\, \det\RG(X,M)^2\otimes\det\RG(X, V)^{2n}\otimes\det\RG(X,\cO)^{-4n}
$$
Indeed, for $V=\cO^2$ we have $\det\RG(X,M\otimes V)\,\iso\, \det\RG(X,M)^2$.
Further, for $M=\cO^n\oplus\Omega^n$ 
$$
\det\RG(X, M\otimes V)\,\iso\, \det\RG(X,V)^n\otimes \det\RG(X,V\otimes\Omega)^n\,\iso\, \det\RG(X,V)^{2n}
$$ 
Since $\H^0(\Bun_{G_n},\cO)=\H^0(\Bun_{\SL_2},\cO)=k$, the assertion follows.

\Step 2
Let $\cF^0_W$ be the trivial $\SO_m$-torsor on $X$. 
Restricting $\tau_{n,m}^*\cA_{G_{nm}}$ under $\Bun_{G_n}\,\toup{\id\times\cF^0_W}\,
\Bun_{G_n}\times\Bun_{\SO_m}$, we get $\cA_{G_n}^m$ canonically.

 For $a\in\ZZ/2\ZZ$ denote by $\Bun^a_{\SO_m}$ the corresponding connected component of $\Bun_{\SO_m}$. Let $\cF^0_{G_n}$ be the $G_n$-bundle $\cO^n\oplus\Omega^n$ on $X$. The restriction of $\tau_{n,m}^*\cA_{G_{nm}}$ under $\cF^0_{G_n}\times\id:\Bun_{\SO_m}\to\Bun_{G_n}\times\Bun_{\SO_m}$ is
$\cA_{\SO_m}^{2n}$ canonically. This yields the desired isomorphism over $\Bun_{G_n}\times\Bun_{\SO_m}^0$. 

 If $\cE$ is a line bundle on $X$ of odd degree then $W=\cE\oplus\cE^*\oplus \cO^{m-2}\in \Bun_{\SO_m}^1$. For this $W$ we get
$$
\det\RG(X, M\otimes W)\,\iso\,\det\RG(M\otimes (\cE\oplus\cE^*))\otimes\det\RG(X,M)^{m-2}
$$
By Step 1, 
$$
\det\RG(M\otimes (\cE\oplus\cE^*))\,\iso\, \det\RG(X,M)^2\otimes\det\RG(X, \cE\oplus\cE^*)^{2n}\otimes\det\RG(X,\cO)^{-4n}
$$ 
The desired isomorphism over $\Bun_{G_n}\times\Bun_{\SO_m}^1$ follows. 
\end{Prf}

\bigskip
 
 By the lemma combined with 3.1.2, for $m$ even there is a canonical map 
$$
\tilde\tau_{n,m}: \Bun_{G_n}\times\Bun_{\SO_m}\to \Bunt_{G_{nm}}
$$ 
extending $\tau_{n,m}$. For $m$ odd there is a canonical map 
$$
\tilde\tau_{n,m}: \Bunt_{G_n}\times\Bun_{\SO_m}\to \Bunt_{G_{nm}}
$$ 
extending $\tau_{n,m}$. 

The complex $\tilde\tau^*_{n,m}\Aut$ viewed as a kernel of intergal operators gives rise to a pair of functors between the categories $\D(\Bunt_{G_n})$ and $\D(\Bun_{\SO_m})$ (for $m$ even one may replace $\Bunt_{G_n}$ by $\Bun_{G_n}$). These functors are the geometric counterpart of the classical theta-lifting (in the nonramified case) for the dual reductive pair $\Sp_{2n}, \SO_m$ (cf., for example, \cite{Pr}, Sect.~8), we will study them separately.  

\bigskip\bigskip
\centerline{\scshape 8. Genuine spherical sheaves on $\wt\Gr_G$}

\bigskip\noindent
8.1 Let $\cO=k[[t]]$ and $K=k((t))$. Let $\Omega_{\cO}$ denote the completed module of relative differentials of $\cO$ over $k$. Pick a free $\cO$-module $M_0$ of rank $2n$ with symplectic form $\wedge^2 M_0\to\Omega_{\cO}$.

 In Sect.~8.1-8.2 $G$ will denote the sheaf of automorphisms of $M_0$ preserving the symplectic form. One associates to $G$ the affine grassmanian $\Gr_G$ (cf. \cite{BD}, p. 172 or \cite{F}), which is an ind-scheme over $k$, the fpqc quotient $\Gr_G=G(K)/G(\cO)$. Here $G(\cO)$ (resp., $G(K)$) is the functor associating to a $k$-algebra $R$ the group of automorphisms of $M_{0,R}:=M_0\otimes_{\cO} R[[t]]$ (resp., of $M_0\otimes_{\cO} R((t))$) preserving all the structures.

  Recall that the Picard group of $\Gr_G$ is $\ZZ$ (cf. \cite{F}), let us introduce the notation for the generator. We have the affine grassmanian $\Gr_{\SL(M_0)}$. Its $R$-points are projective $R[[t]]$-modules of finite type $M\subset M_0\otimes_{\cO} R((t))$ with
\begin{itemize}
\item $t^m M_{0,R}\subset M\subset t^{-m} M_{0,R}$ for some $m>>0$;
\item $\det_{R[[t]]} M=\det_{R[[t]]}M_{0,R}$ as a subspace of $(\det_{R[[t]]}M_{0,R})\otimes_{R[[t]]} R((t))$
\end{itemize}

 We postpone to Lemma~\ref{Lm_covering_Gr} the proof of the fact that $M/t^mM_{0,R}$ is a projective $R$-module for $m>>0$. This allows to introduce the line bundle $\cL$ on $\Gr_{\SL(M_0)}$ whose fibre at $M$ is 
$$
\det(M_0:M):={\det}_R(M_0/t^mM_0)\otimes{\det}_R(M/t^mM_0)^{-1},
$$  
independent of $m$ such that $t^mM_0\subset M$. View it as $\ZZ/2\ZZ$-graded purely of degree zero. 

 The standard representation of $G$ yields a map $\Gr_G\to\Gr_{\SL(M_0)}$, and we also write $\cL$ for the restriction of this line bundle to $\Gr_G$. Then $\cL$ generates the Picard group of $\Gr_G$. Recall that $\cL$ is $G(\cO)$-equivariant on $\Gr_G$. 
  Let $\wt\Gr_G\to \Gr_G$ denote the $\mu_2$-gerbe of square roots of $\cL$. Then
$G(\cO)$ acts on $\wt\Gr_G$ extending the action on $\Gr_G$ (cf. A.3).

\begin{Def} 
\label{Def_genuine}
Let $\Sph(\wt\Gr_G)$ be the category of $G(\cO)$-equivariant perverse sheaves on $\wt\Gr_G$ on which $-1\in\mu_2$ acts as $-1$. We call it the category of \select{genuine spherical sheaves} on $\wt\Gr_G$.
\end{Def}

 A $\theta$-characteristic is a free $\cO$-module $\cN$ of rank 1 together with $\cN\otimes_{\cO}\cN\,\iso\,\Omega_{\cO}$. A choice of a $\theta$-characteristic yields an isomorphism of group schemes $G(\cO)\,\iso \Sp(M_0\otimes_{\cO} \cN^{-1})$ over $k$. A further choice of a symplectic base in $M_0\otimes_{\cO} \cN^{-1}$ over $\cO$ identifies $G(\cO)$ with $\Sp_{2n}(\cO)$. So, we may view  the standard maximal torus and Borel  $T\subset B\subset \Sp_{2n}\subset \Sp_{2n}(\cO)$ as subgroups of $G(\cO)$. Write $\Lambda^+$ for the set of dominant coweights of $\Sp_{2n}$. 

 We have a stratification of $\Gr_G$ by $G(\cO)$-orbits indexed by $\Lambda^+$, write $\Gr_G^{\lambda}$ for the $G(\cO)$-orbit passing by $\lambda(t)\in T(K)$ (\cite{BD}, p. 180). Let $\wt\Gr_G^{\lambda}$ be the preimage of $\Gr_G^{\lambda}$ in $\wt\Gr_G$. 

\begin{Pp}
\label{Pp_euqivariant_trivialization}
 For any $\lambda\in\Lambda^+$ there is a $G(\cO)$-equivariant trivialization $\wt\Gr_G^{\lambda}\,\iso\,\Gr_G^{\lambda}\times B(\mu_2)$, the $G(\cO)$-action on the RHS being the product of the action on $\Gr_G^{\lambda}$ and the trivial action on $B(\mu_2)$.
\end{Pp}
\begin{Prf}
\Step 1  For $\lambda\in\Lambda^+$ denote by $\St_{\lambda}$ the stabilizor of $\lambda(t)\in\Gr_G$ in $G(\cO)$. Let $M_{\lambda}=\lambda(t)M_0$ and $M'=M_0+M_{\lambda}$, $M''=M_0\cap M_{\lambda}$. 
 
 The symplectic form $\wedge^2(M_0\otimes_{\cO} K)\to \Omega(K)=\Omega_{\cO}\otimes_{\cO}K$
induces a map $(M'/M_0)\otimes (M'/M_{\lambda})\,\iso\,(M_{\lambda}/M'')\otimes(M_0/M'')\to \Omega(K)/\Omega_{\cO}$. Composing further with the residue map, we get a pairing between $k$-vector spaces $M'/M_0$ and $M'/M_{\lambda}$. We'll check in Step 2 that the pairing is perfect.
So, the fibre of $\cL$ at $M_{\lambda}$ is 
$$
\cL_{M_{\lambda}}\;\iso\; \det(M_0 : M_{\lambda})\;\iso\;
\frac{\det(M'/M_{\lambda})}{\det(M'/M_0)}\;\iso\; \det(M'/M_{\lambda})^{\otimes 2}
$$
The group $\St_{\lambda}$ acts on $\det(M'/M_{\lambda})$ by some character 
$\chi:\St_{\lambda}\to\Gm$. So, $\St_{\lambda}$ acts on $\cL_{M_{\lambda}}$ by $\chi^2$. Let $\cB$ be the $G(\cO)$-equivariant line bundle on $\Gr^{\lambda}_G$ corresponding to $\chi$. Then we have a $G(\cO)$-equivariant isomorphism $\cB^2\iso \cL\mid_{\Gr_G^{\lambda}}$, and our assertion follows from Lemma~\ref{Lm_A_equivariant}.

\Step 2 Realize $\Sp_{2n}$ as the subgroup of $\SL_{2n}$ preserving the form on $k^{2n}$ given by the matrix 
$$
\left(
\begin{array}{cc}
0 & E_n\\
-E_n & 0
\end{array}
\right),
$$
where $E_n$ is the identity matrix in $\SL_n$. Let $T\subset\Sp_{2n}$ be the maximal torus of diagonal matrices. A coweight $\lambda=(a_1,\ldots, a_n; -a_1,\ldots,-a_n)$ of $T$ is dominant iff $a_1\ge\ldots\ge a_n\ge 0$. Pick a trivialization $\cN\,\iso\cO$ and a symplectic base $e_i$  in $M_0$. Then 
$$
M_{\lambda}=t^{a_1}\cO e_1\oplus\ldots\oplus t^{a_n}\cO e_n\oplus
t^{-a_1}\cO e_{n+1}\oplus\ldots\oplus t^{-a_n}\cO e_{2n}
$$
and $M'=\cO e_1\oplus\ldots\oplus \cO e_n\oplus
t^{-a_1}\cO e_{n+1}\oplus\ldots\oplus t^{-a_n}\cO e_{2n}$. Since
$$
\begin{array}{l}
M'/M_0\,\iso\, t^{-a_1}\cO e_{n+1}\oplus\ldots\oplus t^{-a_n}\cO e_{2n}/\cO e_{n+1}\oplus\ldots\oplus \cO e_{2n}\\
\\
M'/M_{\lambda}\,\iso\, \cO e_1\oplus\ldots\oplus \cO e_n/
t^{a_1}\cO e_1\oplus\ldots\oplus t^{a_n}\cO e_n,
\end{array}
$$
the pairing is perfect.
\end{Prf}

\bigskip

  Let $W$ denote the nontrivial local system of rank one on $B(\mu_2)$ corresponding to the covering $\Spec k\to B(\mu_2)$. For $\lambda\in\Lambda^+$ there is a unique irreducible $G(\cO)$-equivariant perverse sheaf on $\wt\Gr_G^{\lambda}$,  on which $-1\in \mu_2$ acts as $-1$, namely $(\Qlb\boxtimes W)\otimes \Qlb[1](\frac{1}{2})^{\otimes \dim\Gr_G^{\lambda}}$. Denote by $\cA_{\lambda}$ its Goresky-MacPherson extension to $\wt\Gr_G$. By Proposition~\ref{Pp_euqivariant_trivialization}, the irreducible objects of the category $\Sph(\wt\Gr_G)$ are exactly $\cA_{\lambda}, \lambda\in\Lambda^+$.

 Note that $\Sph(\wt\Gr_G)$ is closed under extensions in $\P(\wt\Gr_G)$ (if $-1\in\mu_2$ acts as $-1$ on perverse sheaves $K_1,K_2$ then it acts as $-1$ on any extension of $K_1$ by $K_2$). Since $\DD(\cA_{\lambda})\iso\cA_{\lambda}$ canonically, $\Sph(\wt\Gr_G)$ is preserved by Verdier duality.

  Consider the action of the torus $T\subset G(\cO)$ on $\Gr_G$. 
The following will be used in Sect.~8.4.

\begin{Lm} 
\label{Lm_covering_Gr}
i) There is a covering of $\Gr_G$ by $T$-invariant open ind-schemes $U_i$ and $T$-equivariant trivializations $\cL\mid_{U_i}\,\iso\,\cO_{U_i}$.\\
ii) For an $R$-point $M\subset M_0\otimes_{\cO} R((t))$ of $\Gr_{\SL(M_0)}$ and $m>>0$ the $R$-module $M/t^m M_{0,R}$ is projective. 
\end{Lm}
\begin{Prf} i) Pick a trivialization $\cN\,\iso\,\cO$, so that our base of $M_0\otimes\cN^{-1}$ gives rise to a base $\{e_1,\ldots, e_{2n}\}$ of $M_0$. Consider the corresponding maximal torus $T'$ of $\SL(M_0)$. Set $M^-=Ae_1\oplus\ldots\oplus Ae_{2n}$ with $A=t^{-1}k[t^{-1}]$. For a coweight $\lambda:\Gm\to T'$ of $\SL(M_0)$ denote by  $U_{\lambda}\subset \Gr_{\SL(M_0)}$ the open locus classifying lattices $M\subset M_0\otimes_{\cO} K$ such that $M\oplus \lambda(t)M^-=M_0\otimes_{\cO} K$. Here $\lambda=(b_1,\ldots,b_{2n})$ with $\sum b_i=0$ and $\lambda(t)M^{-}=At^{b_1}e_1\oplus\ldots\oplus At^{b_{2n}}e_{2n}$. 

 One checks that  the union of $U_{\lambda}$ is $\Gr_{\SL(M_0)}$. Clearly, $U_{\lambda}$ is  
$T'$-invariant. As shown by Faltings (\cite{F}, Sect.~2), for each $\lambda$ there is a trivialization $\cL\mid_{U_{\lambda}}\,\iso\,\cO_{U_{\lambda}}$ equivariant under the stabilizor of $\lambda(t)M^-$ in $\SL(M_0)(K)$. This stabilizor contains $T'$, so the trivializations are $T'$-equivariant.

Restricting everything under the map $\Gr_G\to \Gr_{\SL(M_0)}$ corresponding to the standard representation, one concludes the proof.

\smallskip\noindent
ii) (argument due to the unknown referee) \  
Localizing in Zarisky topology of $R$, pick a coweight $\lambda$ of $\SL(M_0)$ such that $M\oplus \lambda(t)M^-_R=M_0\otimes_{\cO} R((t))$. Here $M^-_R=A_Re_1\oplus\ldots\oplus A_Re_{2n}$ and $A_R=t^{-1}R[t^{-1}]$. For $m>>0$ we get $t^{-m}M_{0,R}=M\oplus U$, where $U=\lambda(t)M^-_R\cap t^{-m} M_{0,R}$, and 
$$
(M/t^m M_{0,R})\oplus U\,\iso\, t^{-m}M_{0,R}/t^m M_{0,R}
$$
\end{Prf}

\bigskip\noindent
8.2 {\scshape The convolution product.} Following \cite{MV}, consider the diagram 
$$
\Gr_G\times\Gr_G\getsup{p_G\times\id} G(K)\times\Gr_G \toup{q_G} G(K)\times_{G(\cO)}\Gr_G \toup{m} \Gr_G,
$$
Here $p_G: G(K)\to\Gr_G$ is the projection, $G(K)\times_{G(\cO)}\Gr_G$ is the quotient of $G(K)\times\Gr_G$ by $G(\cO)$, where the action is given by
$x(g,hG(\cO))=(gx^{-1}, xhG(\cO))$ for $x\in G(\cO)$, and $m$ is the product map.

 The map $p_G\times m: G(K)\times_{G(\cO)}\Gr_G\to \Gr_G\times\Gr_G$ sending $(g, hG(\cO))$ to $(gG(\cO), ghG(\cO))$ is an isomorphism.

 We have a canonical isomorphism $q_G^*m^*\cL\,\iso\, p_G^*\cL\boxtimes\cL$. Moreover, the above $G(\cO)$-action on $G(K)\times\Gr_G$ lifts to a $G(\cO)$-equivariant structure on $p_G^*\cL\boxtimes\cL$ giving rise to the line bundle $p_G^*\cL\tboxtimes \cL$ on $G(K)\times_{G(\cO)}\Gr_G$. Thus, $m^*\cL\,\iso\, p_G^*\cL\tboxtimes \cL$ canonically.
 
 Set $\wt{G(K)}=G(K)\times_{\Gr_G}\wt\Gr_G$. Both actions of $G(\cO)$ on $G(K)$ by left and right translations extend naturally to actions on $\wt{G(K)}$. We'll refer to them again as actions by left and right translations, by abuse of terminology. Under the action on $\wt{G(K)}$ by right translations, the projection $\tilde p_G:\wt{G(K)}\to\wt\Gr_G$ is a $G(\cO)$-torsor (cf. A.2).

 Taking the tensor product of square roots of $p_G^*\cL$ and of $\cL$, we get a map $\tilde m$ that fits into the diagram
$$
\begin{array}{ccc}
\wt{G(K)}\times\wt\Gr_G & \toup{\tilde m} & \wt\Gr_G\\
\downarrow && \downarrow\\
G(K)\times\Gr_G & \toup{m\comp q_G} & \Gr_G
\end{array}
$$

 One checks that 
\begin{equation}
\label{map_G(O)-torsor}
\tilde p_G\times\tilde m: \wt{G(K)}\times\wt\Gr_G\to 
\wt\Gr_G\times\wt\Gr_G
\end{equation}
is a $G(\cO)$-torsor, where $G(\cO)$ acts on 
$\wt{G(K)}\times\wt\Gr_G$ as the product of the action by right translations on $\wt{G(K)}$ with the action on $\wt\Gr_G$. 

 Consider the diagram
$$
\wt\Gr_G\times\wt\Gr_G \,\getsup{\tilde p_G\times\id}\,
\wt{G(K)}\times\wt\Gr_G \,\toup{\tilde p_G\times\tilde m} \,\wt\Gr_G\times\wt\Gr_G\toup{pr_2}\wt\Gr_G
$$

\begin{Def} For $K_1,K_2\in \Sph(\wt\Gr_G)$ define the convolution product $K_1\ast K_2\in \D(\wt\Gr_G)$ by 
$$
K_1\ast K_2=\pr_{2!} K,
$$ 
where $K$ is a perverse sheaf on $\wt\Gr_G\times\wt\Gr_G$ such that $(\tilde p_G\times\tilde m)^*K\,\iso\, \tilde p_G^*K_1\boxtimes K_2$. Since (\ref{map_G(O)-torsor}) is a $G(\cO)$-torsor and $ \tilde p_G^*K_1\boxtimes K_2$ is equivariant under the corresponding $G(\cO)$-action on $\wt{G(K)}\times\wt\Gr_G$, $K$ is defined up to a unique isomorphism (cf. A.2).
 
 For $(a, b)\in\mu_2\times\mu_2$ the image under $\tilde p_G\times\tilde m$ of the corresponding 2-automorphism of $\wt{G(K)}\times\wt\Gr_G$ is the 2-automorphism
$(a, ab)$ of $\wt\Gr_G\times\wt\Gr_G$. So, by Lemma~\ref{Lm_Lemma_descent_on_gerb}, $K$ descends to a perverse sheaf $K'$ on $\Gr_G\times\wt\Gr_G$ (such $K'$ is defined up to a unique isomorphism). Since
$\RG_c(B(\mu_2),\Qlb)=\Qlb$, we see that $K_1\ast K_2\,\iso\, \pr_{2!} K'$, where
$\pr_2:  \Gr_G\times\wt\Gr_G\to\wt\Gr_G$ is the projection. Moreover, $-1\in\mu_2$ acts on $K_1\ast K_2$ as $-1$.
\end{Def}

\begin{Pp} For $K_1,K_2\in\Sph(\wt\Gr_G)$ we have $K_1\ast K_2\in\Sph(\wt\Gr_G)$.
\end{Pp}
\begin{Prf} Following \cite{MV}, stratify $\Gr_G\times\wt\Gr_G$ by locally closed substacks $\wt\Gr_G^{\lambda,\mu}$, $\lambda,\mu\in\Lambda^+$, where $\wt\Gr_G^{\lambda,\mu}$ is the preimage of $(p_G\times m)(p_G^{-1}(\Gr_G^{\lambda})\times_{G(\cO)}\Gr_G^{\mu})$ under $\Gr_G\times\wt\Gr_G\to\Gr_G\times\Gr_G$. 

 Stratify also $\wt\Gr_G$ by $\wt\Gr_G^{\lambda}$, $\lambda\in\Lambda^+$. By Lemma~4.3 of \select{loc.cit.}, $\pr_2:\Gr_G\times\wt\Gr_G\to\wt\Gr_G$ is stratified semi-small map. Our assertion follows from Lemma~4.2 of \select{loc.cit.}
\end{Prf}

\medskip

 In a similar way one defines a convolution product $K_1\ast K_2\ast K_3$ of three sheaves $K_1,K_2,K_3\in\Sph(\wt\Gr_G)$. Moreover, $(K_1\ast K_2)\ast K_3\,\iso\,
K_1\ast K_2\ast K_3\,\iso\, K_1\ast (K_2\ast K_3)$ canonically, and $\cA_0$ is a unit object. So, $\Sph(\wt\Gr_G)$ is an associative tensor category (a category with tensor functor and an associativity constraint).

 Observe that for each $\lambda\in\Lambda^+$ the $G(\cO)$-orbit $\Gr_G^{\lambda}$ is even-dimensional. 

\begin{Pp} 
\label{Pp_semisimplicity}
1) For $\lambda\in\Lambda^+$ the odd cohomology sheaves of $\cA_{\lambda}$ (with respect to the usual t-structure) vanish. \\ 
2) The category $\Sph(\wt\Gr_G)$ is semisimple.
\end{Pp}
\begin{Prf}
1a) Given $\lambda_1,\ldots,\lambda_r\in\Lambda^+$, consider the convolution diagram
$$
m: \Conv^{\lambda_1,\ldots,\lambda_r}\to \Grb^{\lambda_1+\ldots+\lambda_r}_G,
$$ 
where we have set $\Conv^{\lambda_1,\ldots,\lambda_r}=\Gr^{\lambda_1}_G\ttimes\ldots\ttimes\Gr^{\lambda_r}_G$. Let $\wt\Conv^{\lambda_1,\ldots,\lambda_r}$ be the restriction of the gerbe $\wt\Gr_G$ under the above map $m$. The canonical section $s: \Gr_G^{\lambda_1+\ldots+\lambda_r}\to\wt\Gr_G^{\lambda_1+\ldots+\lambda_r}$ yields a section 
$m^{-1}(s)$ of the gerbe $\wt\Conv^{\lambda_1,\ldots,\lambda_r}$ over $m^{-1}(\Gr_G^{\lambda_1+\ldots+\lambda_r})$. One checks that this section extends to a section
$\Conv^{\lambda_1,\ldots,\lambda_r}\to \wt\Conv^{\lambda_1,\ldots,\lambda_r}$. 

\medskip\noindent
1b) We adopt Gaitsgory's proof of a theorem of Lusztig to our situation (\cite{G}, A.7). Namely, let $\Fl$ denote the affine flag variety. This is the ind-scheme classifying a $G$-bundle $\cF_G$ on $\Spec\cO$ with trivialization $\cF_G\,\iso\,\cF^0_G\mid_{\Spec K}$ and a reduction of $\cF_G\mid_{\Spec \cO/(t)}$ to the Borel subgroup $B$. 

 Let $\wt\Fl$ denote the restriction of the gerbe $\wt\Gr_G$ under the (smooth) projection $\Fl\to\Gr_G$. Let $I\subset G(\cO)$ be the Iwahory subgroup. For an element $w$ of the affine Weil group of $G$, let $\Fl^w$ denote the corresponding $I$-orbit on $\Fl$. Set $\wt\Fl^w=\Fl^w\times_{\Fl}\wt\Fl$. 
 
  Let $\mu\in\Lambda^+$ be such that the projection $\Fl^w\to\Gr_G$ factors through $\Gr_G^{\mu}$. The canonical section $\Gr_G^{\mu}\to\wt\Gr_G^{\mu}$ yields a section $s: \Fl^w\to \wt\Fl^w$ of the gerbe $\wt\Fl^w\to\Fl^w$. Let $\cA_w$ denote the irreducible perverse sheaf on the closure of $\wt\Fl^w$ on which $-1\in\mu_2$ acts as $-1$ and whose restriction under $s$ is $\IC_{\Fl^w}$. It suffices to show the parity vanishing for stalks of $\cA_w$. 
  
   Let $w=s_1\cdot\ldots\cdot s_r$ be a reduced decomposition of $w$ into a product of simple reflections. Denote by $p: \Conv^{s_1,\ldots,s_r}_{\Fl}\to \ov{\Fl}^w$ the Bott-Samelson resolution (\select{loc.cit.} or \cite{F}, Sect.~3, where it is referred to as Demazure resolution).  
Let $\wt\Conv^{s_1,\ldots,s_r}_{\Fl}$ be the restriction of our gerbe to  $\Conv^{s_1,\ldots,s_r}_{\Fl}$. By 1a), the section 
$$
p^{-1}(\Fl^w)\to p^{-1}(\wt\Fl^w)
$$ induced by $s$ extends to a global section 
$\Conv^{s_1,\ldots,s_r}_{\Fl}\to \wt\Conv^{s_1,\ldots,s_r}_{\Fl}$.  
The desired assertion follows, because the fibres of $p$ have cohomology with compact support in even degrees only (\cite{G}, A.7).

\medskip\noindent
2) Follows from 1) as in (\cite{BD}, 5.3.3). This uses the fact that each $\Gr_G^{\lambda}$ has cohomology only in even degrees (5.3.2 of \select{loc.cit.}).
\end{Prf}

\begin{Rem} 
\label{Rem_AutO}
The group of automorphisms of the $k$-algebra $\cO$ is naturally the group of $k$-points of a (reduced) affine group scheme $\Aut^0\cO$ over $k$ (\cite{BD}, 2.6.5). Assume that $M_0=\cO^n\oplus\Omega^n_{\cO}$ with standard symplectic form. Then $\Aut^0\cO$ acts on $M_0$ and, hence, on $\Gr_G$. Moreover, $\cL$ is naturally equivariant under this action. It follows that $\Aut^0\cO$ acts on $\wt\Gr_G$. Proposition~\ref{Pp_euqivariant_trivialization} then can be strengthened saying that the gerbe $\wt\Gr_G^{\lambda}\to\Gr_G^{\lambda}$ admits a $G(\cO)\rtimes\Aut^0\cO$-equivariant trivialization.  

 We also see that each $\cA_{\lambda}$ is $G(\cO)\rtimes\Aut^0\cO$-equivariant (this property is true for the constant sheaf over $\Gr_G^{\lambda}$ and is preserved under intermediate extension). By Proposition~\ref{Pp_semisimplicity}, each $K\in\Sph(\wt\Gr_G)$ is $\Aut^0\cO$-equivariant. Moreover, such equivariant structure is unique (because the stabilizer of a point is connected) and compatible with any morphism in $\Sph(\wt\Gr_G)$.  
\end{Rem} 

\bigskip
\noindent
8.3 {\scshape The fusion product}  Following \cite{MV}, we will show that the convolution product defined above can be interpreted as a `fusion' product, thus leading to a commutativity constraint on $\Sph(\wt\Gr_G)$. The original idea of this interpretation for spherical sheaves on $\Gr_G$ is due to V.~Drinfeld.

 Let $G$ denote the sheaf of groups on $X$ introduced in Sect.~3.2. For $x\in X(k)$ write $\cO_x$ for the completed local ring at $x$ and $K_x$ for its fraction field. Write $\Gr_{G,x}=G(K_x)/G(\cO_x)$ for the corresponding version of the affine grassmanian. 

 Write $\cF^0_G$ for the `trivial' $G$-torsor on $X$ given by $M_0=\cO_X^n\oplus\Omega^n$ with standard symplectic form $\wedge^2 M_0\to \Omega$.
 
 For a $k$-algebra $R$ write $X_R=X\times\Spec R$ and $X^*_R=(X-x)\times\Spec R$.
By \cite{BL1,BL2},  $\Gr_{G,x}$ is the functor on the category of $k$-algebras sending $R$ to the set  of isomorphism classes of $\{\cF_G, \nu\}$, where $\cF_G$ is a $G$-torsor on $X_R$ and $\nu: \cF_G\,\iso\, \cF^0_G\mid_{X^*_R}$ is a trivialization outside $x$.   
 
 Let $M$ denote the vector bundle on $X$ induced from $\cF_G$ via the standard representation of $G$. Set $M_x=M\otimes\cO_x$ and $M_{0,x}=M_0\otimes \cO_x$. Then $M_x\subset M_{0,x}\otimes_{\cO_x} K_x$ is a sublattice, and we continue to denote by $\cL$ the line bundle on $\Gr_{G,x}$ with fibre $\det(M_{0,x}: M_x)$. Then $\wt\Gr_{G,x}$ and $\Sph(\wt\Gr_{G,x})$ are defined as in Sect.~8.1.

 Write $\Gr_{G,X^d}$ for the functor associating to a $k$-algebra $R$ the set 
$$
\{(x_1,\ldots,x_d)\in X^d(R),  \;\mbox{a}\; G\!-\!\mbox{torsor}\; \cF_G\; \mbox{on}\; X_R,\;
\cF_G\,\iso\,\cF^0_G\mid_{X_R-\cup x_i}\}
$$
Here $x_i\in X(R)$ are thought of as subschemes in $X_R$ by taking their graphs.

 Let $G_{X^d}$ denote the functor sending a $k$-algebra $R$ to the set $\{(x_1,\ldots,x_d)\in X^d(R), \mu)$, where $\mu$ is an automorphism of $\cF^0_G$ restricted to the formal neighborhood 
$\widehat X_{R,D}$ of $D=x_1\cup\ldots\cup x_d$ in $X_R$. So, $G_{X^d}$ is a group scheme over $X^d$, whose fibre over $(x_1,\ldots,x_d)$ is $\prod_i G(\cO_{y_i})$ with $\{y_1,\ldots,y_s\}=\{x_1,\ldots,x_d\}$ and $y_i$ pairwise distinct. 

 Let $\cL$ be the line bundle on $\Gr_{G,X^n}$ whose fibre is $\det\RG(X, M_0)\otimes\det\RG(X,M)^{-1}$, where $M$ is the vector bundle on $X$ induced from $\cF_G$ via the standard representation of $G$.
  
\begin{Lm} 
\label{Lm_factorisation_property}
For a $k$-point $(x_1,\ldots,x_d,\, \cF_G)$ of $\Gr_{G,X^d}$ let $\{y_1,\ldots,y_s\}=\{x_1,\ldots,x_d\}$ with $y_i$ pairwise distinct. The fibre of $\cL$ at this $k$-point is canonically isomorphic (as $\ZZ/2\ZZ$-graded) to 
$$
\otimes_{i=1}^s \det(M_{0,y_i}: M_{y_i}) \eqno{\square}
$$
\end{Lm}

\smallskip

 One checks that the natural action of $G_{X^d}$ on $\Gr_{G,X^d}$ lifts to a $G_{X^d}$-equivariant structure on $\cL$. We have $\wt\Gr_{G,X^d}$ and $\Sph(\wt\Gr_{G,X^d})$ defined as above.

\bigskip\noindent
8.3.1  Consider the diagram of stacks over $X^2$, where the left and right square is cartesian
$$
\begin{array}{ccccccc}
\wt\Gr_{G,X}\times\wt\Gr_{G,X} & \getsup{\tilde p_{G,X}} & \wt{C}_{G,X} &\toup{\tilde q_{G,X}}& \wt\Conv_{G,X} & \toup{\tilde m_X} & \wt\Gr_{G,X^2}\\
\downarrow && \downarrow &&  \downarrow && \downarrow \\
\Gr_{G,X}\times\Gr_{G,X} & \getsup{p_{G,X}} & C_{G,X} & \toup{q_{G,X}}  & \Conv_{G,X} & \toup{m_X} & \Gr_{G,X^2}
\end{array}
$$
Here the low row is the usual convolution diagram \cite{MV}, (5.2). Namely, $C_{G,X}$ is the ind-scheme classifying collections:
\begin{equation}
\label{equation_C_GX}
\left\{
\begin{array}{l}
x_1,x_2\in X, \; G-\!\mbox{torsors} \; \cF^1_G, \cF^2_G \;\mbox{on}\;  X \; \mbox{with trivializations} \; \nu_i:\cF^i_G\,\iso\,\cF^0_G\mid_{X-x_i}, \\ 
\mu_1: \cF^1_G\,\iso\,\cF^0_G\mid_{\widehat{X}_{x_2}},
\end{array} 
\right.
\end{equation}
where $\widehat{X}_{x_2}$ is  the formal neighborhood of $x_2$ in $X$. The map $p_{G,X}$ forgets $\mu_1$.

 The ind-scheme $\Conv_{G,X}$ classifies collections:
\begin{equation}
\label{equation_conv_GX}
\left\{
\begin{array}{l}
x_1,x_2\in X, \; G-\!\mbox{torsors} \; \cF^1_G, \cF_G \;\mbox{on}\;  X, \\
\mbox{isomorphisms} \; \nu_1:\cF^1_G\,\iso\,\cF^0_G\mid_{X-x_1}, \;\mbox{and}\; 
\eta: \cF^1_G\,\iso\,\cF_G\mid_{X-x_2}
\end{array} 
\right.
\end{equation}
The map $m_X$ sends this collection to $(x_1,x_2,\cF_G)$ together with the trivialization $\eta\comp\nu_1^{-1}: \cF^0_G\,\iso\,\cF_G\mid_{X-x_1-x_2}$. 

 The map $q_{G,X}$ sends (\ref{equation_C_GX}) to the collection (\ref{equation_conv_GX}), where $\cF_G$ is obtained by gluing $\cF^1_G$ on $X-x_2$ and $\cF^2_G$ on $\widehat{X}_{x_2}$ using their identification over $(X-x_2)\cap\widehat{X}_{x_2}$ via $\nu_2^{-1}\comp \mu_1$. 
 
 The canonical isomorphism
$$
q_{G,X}^*m_X^*\cL\,\iso\, p_{G,X}^*(\cL\boxtimes\cL)
$$ 
allows to define $\tilde q_{G,X}$ as follows.
Write $M_i$ (resp., $M$) for the vector bundle induced from $\cF_G^i$ (resp., $\cF_G$) via the standard representation of $G$. 

 A point of $\tilde C_{G,X}$ is given by (\ref{equation_C_GX}) together with 1-dimensional vector spaces $\cB_1,\cB_2$ and $\cB_i^2\,\iso\, \cL_{\cF^i_G}$. 
By Lemma~\ref{Lm_factorisation_property}, $\cL_{\cF^i_G}\,\iso\, \det( M_{0,x_i}:\det M_{i,x_i})$.

 A point of $\wt\Conv_{G,X}$ is given by (\ref{equation_conv_GX}) together with 1-dimensional vector space $\cB$ and
$\cB^2\,\iso\,\cL_{\cF_G}$. We have 
$$
\cL_{\cF_G}\,\iso\, \frac{\det\RG(X, M_0)}{\det\RG(X, M_1)}\otimes
\frac{\det\RG(X, M_1)}{\det\RG(X, M)}\,\iso
\; \det(M_{0,x_1}: M_{1,x_1}) \otimes
\det(M_{1,x_2}: M_{x_2})\;\iso\,\cL_{\cF^1_G}\otimes\cL_{\cF^2_G},
$$
the last isomorphism being given by $\mu_1: \det(M_{1,x_2})\,\iso\,\det(M_{0,x_2})$ and $M_{x_2}\,\iso\, M_{2,x_2}$. Define $\tilde q_{G,X}$ by setting $\cB=\cB_1\otimes\cB_2$. 
 
 As in Sect.~8.2 one checks that for $K_1,K_2\in\Sph(\wt\Gr_{G,X})$ there is a (defined up to a unique isomoprhism) perverse sheaf $K_{12}$ on $\wt\Conv_{G,X}$ with
$\tilde q_{G,X}^*K_{12}\,\iso\, \tilde p_{G,X}^*(K_1\boxtimes K_2)$. Moreover, $-1\in\mu_2$ acts on $K_{12}$ as $-1$. We then let
$$
K_1 \ast_X K_2=\tilde m_{X !} K_{12}
$$
Let $U\subset X^2$ be the complement to the diagonal. Let $j: \wt\Gr_{G,X^2}(U)\hook{} \wt\Gr_{G,X^2}$ be the preimage of $U$. Recall that $m_X$ is stratified small, an isomorphism over the preimage of $U$ (\cite{MV}). So, the same holds for the representable map $\tilde m_X$. Thus, $K_1\ast_X K_2$ is a perverse sheaf, the Goresky-MacPherson from $\wt\Gr_{G,X^2}(U)$. Besides, $-1\in\mu_2$ acts on it as $-1$. Moreover, $K_1\ast_X K_2\in\Sph(\wt\Gr_{G,X^2})$, because $G_{X^2}$-equivariance is clear over $\wt\Gr_{G,X^2}(U)$ and is preserved under the intermediate
extension.

 Recall the group ind-scheme $\Aut^0\cO$ (cf. Remark~\ref{Rem_AutO}). Let $\hat X\to X$ be the $\Aut^0\cO$-torsor whose fibre is the set of all trivializations $\cO_x\,\iso\cO$. It is known that $\Gr_{G,X}\,\iso\, \hat X\times_{\Aut^0\cO}\Gr_G$ (\cite{BD}, 5.3.11). The line bundle $\cL$ on 
$\Gr_{G,X}$ identifies with the descent of the $\Aut^0\cO$-equivariant line bundle $\cO\boxtimes \cL$ under $\hat X\times \Gr_G\to \Gr_{G,X}$. 
Since any $K\in\Sph(\wt\Gr_G)$ is $\Aut^0\cO$-equivariant, we have a natural (fully faithful) functor 
\begin{equation}
\label{functor_tau_0}
\tau^0: \Sph(\wt\Gr_G)\to\Sph(\wt\Gr_{G,X})[-1]
\end{equation}
Let $\glob: \Sph(\wt\Gr_G)\to\Sph(\wt\Gr_{G,X})$ denote the functor $\glob=\tau^0[1]$. 

 Now define the commutativity constraint following \cite{MV}. Let $i:\wt\Gr_{G,X}\to \wt\Gr_{G,X^2}$ be the preimage of the diagonal in $X^2$. For $F_1,F_2\in\Sph(\wt\Gr_G)$ letting $K_i=\tau^0 F_i$ define 
$$
K_{12}\mid_U:=K_{12}\mid_{\wt\Gr_{G,X^2}(U)}
$$ 
as above (but now it is placed in perverse degree 2). We get
\begin{eqnarray}
K_1\ast_X K_2\,\iso\, j_{!*}(K_{12}\mid_U)\\
\label{iso_on_diagonal}
\tau^0(F_1\ast F_2)\,\iso\, i^*(K_1\ast_X K_2)
\end{eqnarray}
So, the involution $\sigma$ of $\wt\Gr_{G,X^2}$ interchanging $x_i$ yields
$$
\tau^0(F_1\ast F_2)\,\iso\, i^* j_{!*}(K_{12}\mid_U)\,\iso\, i^* j_{!*}(K_{21}\mid_U)\,\iso\, 
\tau^0(F_2\ast F_1),
$$
because $\sigma^*(K_{12}\mid_U)\,\iso\, K_{21}\mid_U$. (We used the functor $\tau^0$ instead of $\glob$ to avoid the signs ambiguity in the commutativity constraints).
 
 To show that the associativity and commutativity constraints are compatible, argue as in (\cite{BD}, 5.3.13-5.3.17). Namely, one defines for a non-empty finite set $I$ a category $\otimes_I\Sph(\wt\Gr_G)$ and for any surjection $h:I\to I'$ a functor $\ast_h: \otimes_I\Sph(\wt\Gr_G)\to \otimes_{I'}\Sph(\wt\Gr_G)$. They are compatible in the sense of (\select{loc.cit.}, (266)). Thus, $\Sph(\wt\Gr_G)$ is a tensor category. 

\begin{Rem} Fix $x\in X(k)$. Consider the Hecke stack $_x\cH_G$ classifying two $G$-bundles $\cF_G,\cF'_G$ on $X$ together with an isomorphism $\cF_G\,\iso\, \cF'_G\mid_{X-x}$.  Let $p$ (resp., $p'$) be the projection $_x\cH_G\to \Bun_G$ sending the above collection to $\cF_G$ (resp., $\cF'_G$). Write $\Bun_G^x$ for the stack classifying a $G$-torsor $\cF_G$ on $X$ together with a trivialization $\cF_G\,\iso\,\cF^0_G\mid_{D_x}$ over the formal disk $D_x$ around $x$. 

 Let $\gamma$ (resp., $\gamma'$) be the isomorphism $\Bun_G^x\times_{G(\cO_x)}\Gr_{G,x}\,\iso\,{_x\cH_G}$  such that the projection to the first term corresponds to $p$ (resp., to $p'$). 
Write $M$ (resp., $M'$) for the vector bundle corresponding to $\cF_G$ (resp., to $\cF'_G$) via the standard representation of $G$. Write $\cL$ for the ($\ZZ/2\ZZ$-graded) line bundle on $_x\cH_G$ with fibre $\det\RG(X, M)\otimes\det\RG(X, M')^{-1}$. Let $_x\tilde{\cH}_G$ be the gerbe of square roots of $\cL$. Both $\gamma$ and $\gamma'$ extend to $G(\cO_x)$-torsors 
$$
\tilde\gamma, \tilde\gamma':
\Bun_G^x\times\wt\Gr_{G,x}\to {_x{\tilde\cH}_G}
$$

 For $\cS\in\Sph(\wt\Gr_{G,x})$ denote by $\Qlb\tboxtimes\cS$ (resp., by $\Qlb\tboxtimes'\cS$) the twisted tensor product viewed as a perverse sheaf on $_x\tilde\cH_G$ via $\tilde\gamma$ (resp., $\tilde\gamma'$). Given $\cS\in\Sph(\wt\Gr_{G,x})$ there is a (defined up to a unique isomorphism) $\cT\in\Sph(\wt\Gr_{G,x})$ equipped with an isomorphism
$
\Qlb\tboxtimes\cS\,\iso\, \Qlb\tboxtimes'\cT
$. 
This defines a covariant involution functor $\star$ on the category $\Sph(\wt\Gr_{G,x})$ 
By Remark~\ref{Rem_AutO}, we may view $\star$ as an involution functor on $\Sph(\wt\Gr_G)$ independently of a choice of a trivialization $\cO_x\,\iso\,\cO$.

 In the same way as for usual spherical sheaves on $\Gr_G$, one checks that for $K_1,K_2,K_3\in\Sph(\wt\Gr_G)$ we have canonically $R\Hom(K_1\ast K_2, K_3)\,\iso\, R\Hom(K_1, K_3\ast\DD(\star K_2))$. So, $K_3\ast\DD(\star K_2)$ represents the internal $\HOM(K_2,K_3)$  in the sense of the tensor structure on $\Sph(\wt\Gr_G)$. Besides, $\star(K_1\ast K_2)\,\iso\, (\star K_2)\ast (\star K_1)$ canonically. We also have $\DD(\star\cA_{\lambda})\,\iso\,\star\cA_{\lambda}\,\iso\,\cA_{\lambda}$ for each $\lambda\in\Lambda_+$. 
\end{Rem} 

% soslat'sja na LN in Math Deligne? Hom predstavimy? Reflexivity?

\bigskip\noindent
8.4 {\scshape Functors $F^{\theta}$.\ }  Let $P\subset G$ denote the Siegel parabolic preserving  $\cO^n_X\subset \cO^n_X\oplus\Omega^n$. Write $Q$ for the Levi quotient, so $Q\iso\GL_n$ canonically. Let $\check{\Lambda}_{G,P}$ denote the lattice of characters of $P/[P,P]=Q/[Q,Q]$ and $\Lambda_{G,P}$ the dual lattice. Let $\check{\omega}_n\in\check{\Lambda}_{G,P}$ denote the fundamental weight of $G$ corresponding to the unique simple coroot which is not a coroot of $Q$. So, $\check{\omega}_n$ is the
highest weight of an irreducible subrepresentation in $\wedge^n M$, where $M$ is the standard representation of $G$. Then $\check{\omega}_n$ is a free generator of $\check{\Lambda}_{G,P}$. 

 The connected components of $\Gr_{Q,x}$ are indexed by $\Lambda_{G,P}$, the component $\Gr_{Q,x}^{\theta}$ classifies $(L\in\Bun_n, \nu: L\,\iso\,\cO^n\mid_{X-x})$ such that $\deg L=-\<\theta,\check{\omega}_n\> $. 
The reduced part $\Gr_{Q,x,red}^{\theta}\hook{} \Gr_{Q,x}^{\theta}$ is the ind-scheme classifying $(L\in\Bun_n, \nu: L\,\iso\,\cO^n\mid_{X-x})$ that induce  
an isomorphism
\begin{equation}
\label{iso_for_Gr_Q_theta}
\det L\,\iso\, \cO(-\<\theta,\check{\omega}_n\>x)
\end{equation}

 Following \cite{BG}, for $\theta\in\Lambda_{G,P}$ let $S^{\theta}_P$ denote the ind-scheme classifying: $(\cF_P,\nu)$, where $\cF_P$ is a $P$-torsor on $X$ and $\nu:\cF_P\,\iso\,\cF^0_P\mid_{X-x}$ is a trivialization such that $(\cF_P\times_P Q,\nu)$ lies in $\Gr_{Q,x}^{\theta}$. 
In other words,  $S^{\theta}_P$ classifies a $P$-torsor given by an exact sequence $0\to \Sym^2L\to ?\to\Omega\to 0$ on $X$ with $L\in\Bun_n$, a splitting of this sequence over $X-x$, and a trivialization $\nu: L\,\iso\,\cO^n\mid_{X-x}$ with $\deg L=-\<\theta,\check{\omega}_n\>$. 
The reduced part $(S^{\theta}_P)_{red}$ is given by the addidtional condition that $\nu$ induces an isomorphism (\ref{iso_for_Gr_Q_theta}). 
 
 We have a map $\gs^{\theta}_P: S^{\theta}_P\to \Gr_{G,x}$ sending $(\cF_P,\nu)$ to $(\cF_P\times_P \, G,\nu)$, its restriction $(S^{\theta}_P)_{red}\hook{} \Gr_{G,x}$ is a locally closed immersion. 

  The map $\gs^{\theta}_{\bar P}: S^{\theta}_{\bar P}\to\Gr_{G,x}$ is defined in a similar way using the lagrangian subbundle $\Omega^n\subset \cO^n_X\oplus\Omega^n$ that defines the opposite parabolic subgroup $\bar P\subset G$. 
 
  Write $\gt^{\theta}_P: S_P^{\theta}\to \Gr_{Q,x}^{\theta}$ for  the projection sending $(\cF_P,\nu)$ to $(\cF_P\times_P Q,\nu)$ and $\gr^{\theta}_P: \Gr_{Q,x}^{\theta}\hook{} S_P^{\theta}$ for  the natural section, similarly for $\bar P$. 
    
 Fix an isomorphism $\Gm\,\iso\,Z(Q)$, where $Z(Q)$ is the center of $Q$, in such a way that $\Gm\,\iso\,Z(Q)$ acts adjointly on the unipotent radical $U(P)\subset P$ with strictly positive weights. The subscheme of $Z(Q)$-fixed points in $\Gr_G$ is $Q(K)G(\cO)/G(\cO)$, its connected components are $\Gr_{Q,red}^{\theta}$, $\theta\in\Lambda_{G,P}$. 
One checks that 
$$
\{x\in\Gr_{G,x}\mid \lim_{t\to 0} tx\in \Gr_{Q,x,red}^{\theta}\}=(S_P^{\theta})_{red} \;\;\;\; \mbox{and}
$$
$$
\{x\in\Gr_{G,x}\mid \lim_{t\to\infty} tx\in \Gr_{Q,x,red}^{\theta}\}=(S_{\bar P}^{\theta})_{red}
$$

 Consider the diagram
$$
\begin{array}{ccc}
\wt S^{\theta}_P &  \toup{\tilde\gs^{\theta}_P} & \wt\Gr_{G,x}\\
\uparrow\lefteqn{\scriptstyle \tilde\gr^{\theta}_P} && \uparrow\lefteqn{\scriptstyle \tilde\gs^{\theta}_{\bar P}}\\
\wt\Gr^{\theta}_{Q,x} & \toup{\tilde\gr^{\theta}_{\bar P}} & \wt S^{\theta}_{\bar P}
\end{array}
$$
obtained by restricting the gerbe $\wt\Gr_{G,x}\to\Gr_{G,x}$ with respect to the corresponding maps.

\begin{Lm} 
\label{Lm_section_over_unip_orbits}
There exists a canonical $P(\cO_x)$-equivariant section $i^{\theta}_P: S^{\theta}_P\to \wt S^{\theta}_P$ of the gerbe
$\wt S^{\theta}_P\to S^{\theta}_P$.
\end{Lm}
\begin{Prf}
Remind the line bundle $\cL$ on $\Gr_{G,x}$ introduced in 8.3.
Consider the map $\Gr_{G,x}\to\Bun_G$ sending $(\cF_G,\nu: \cF_G\,\iso\,\cF^0_G\mid_{X-x})$ to $\cF_G$. The restriction of $\cA$ under this map identifies canonically with $\cL^{-1}\otimes \det\RG(X,M_0)$, where $M_0=\cO_X^n\oplus\Omega^n$. Since  $\det\RG(X,M_0)\,\iso\, \det\RG(X,\cO)^{\otimes 2n}$, we get a cartesian square
$$
\begin{array}{ccc}
\wt\Gr_{G,x} & \to & \wt\Bun_G\\
\downarrow && \downarrow\\
\Gr_{G,x} & \to & \Bun_G
\end{array}
$$
Remind the map $\tilde \nu$ defined in Lemma~\ref{Lm_tilde_nu}. Now the diagram 
$$
\begin{array}{ccccc}
S^{\theta}_P & \to & \Bun_P & \toup{\tilde\nu} & \wt\Bun_G\\
\downarrow && \downarrow & \swarrow\lefteqn{\scriptstyle \gr}\\
\Gr_{G,x} &\to &  \Bun_G
\end{array}
$$
yields the section $i^{\theta}_P$.

 To see that it is $P(\cO_x)$-equivariant, rewrite it in local terms as follows. On $\Gr_{Q,x}^{\theta}$ we have the $\ZZ/2\ZZ$-graded $Q(\cO_x)$-equivariant line bundle, say $_{\theta}\cL$, whose fibre at $(L, L\,\iso\,\cO^n\mid_{X-x})$ is 
$$
\det(L_0\otimes\cO_x : L\otimes\cO_x)
$$ 
with $L_0=\cO^n_X$. 
Hence $(\gt^{\theta}_P)^*{_{\theta}\cL}$ is a $P(\cO_x)$-equivariant line bundle on $S^{\theta}_P$. The canonical $\ZZ/2\ZZ$-graded $P(\cO_x)$-equivariant isomorphism $(\gs^{\theta}_P)^*\cL\,\iso\, (\gt^{\theta}_P)^*(_{\theta}\cL)^{\otimes 2}$ defines the section $i^{\theta}_P$
via 3.1.2. 
\end{Prf}

\bigskip

  Define the functors $F^{\theta}, F'^{\theta}: \Sph(\wt\Gr_{G,x})\to \D(\Gr^{\theta}_{Q,x})$ by
 $$
F'^{\theta}(K)=(\gt^{\theta}_P)_!(i^{\theta}_P)^*(\tilde \gs^{\theta}_P)^*K \;\;\;\;\;
\mbox{and}\;\;\;\;\;F^{\theta}(K)=F'^{\theta}(K)\otimes\Qlb[1](\frac{1}{2})^{\otimes \<\theta, 2\check{\rho}-2\check{\rho}_Q\>}
$$
We have used the fact that $2(\check{\rho}-\check{\rho}_Q)\in \check{\Lambda}_{G,P}$. 

\begin{Rem} We could replace in the definition of $F^{\theta}$ and  $F'^{\theta}$ the ind-schemes $S^{\theta}_P$ and $\Gr^{\theta}_{Q,x}$ by their reduced parts, the corresponding functors would be canonically isomorphic to the old ones. In some geometric questions we work rather with the corresponding reduced ind-schemes (without indicating that explicitly, for example in Proposition~\ref{Pp_exactness} and \ref{Pp_multiplicities}, Corolary~\ref{Cor1} and so on). 
\end{Rem}

\begin{Pp} 
\label{Pp_exactness}
The functor $F^{\theta}(K)$ maps $\Sph(\wt\Gr_{G,x})$ to the category $\Sph(\Gr_{Q,x}^{\theta})$ of $Q(\cO_x)$-equivariant perverse sheaves on $\Gr_{Q,x}^{\theta}$. 
In particular, it is exact.
\end{Pp}
\begin{Prf}
 By Lemma~\ref{Lm_covering_Gr} combined with Proposition~\ref{Pp_hyperbolic_loc}, we get 
the hyperbolic localization functors $\Sph(\wt\Gr_{G,x})\to \D(\wt\Gr^{\theta}_{Q,x})$
given by
\begin{equation}
\label{eq_hyp_loc}
K\mapsto (\tilde\gr^{\theta}_{\bar P})^*(\tilde \gs^{\theta}_{\bar P})^!K
\;\iso\; (\tilde\gr^{\theta}_P)^!(\tilde \gs^{\theta}_P)^*K=K^{!*}
\end{equation}
By Lemma~\ref{Lm_section_over_unip_orbits}, we have moreover $K^{!*}\,\iso\, (\gt^{\theta}_P\times\id)_!(\tilde s^{\theta}_P)^*K$, where 
$$
\gt^{\theta}_P\times\id: \tilde S^{\theta}_P=S^{\theta}_P\times B(\mu_2)\to \Gr_{Q,x}^{\theta}\times B(\mu_2)=\wt\Gr_{Q,x}^{\theta}
$$ 
The complex $K^{!*}$ is $Q(\cO_x)$-equivariant, because both $\tilde \gs^{\theta}_P$ and $\tilde\gr^{\theta}_P$ are $Q(\cO_x)$-equivariant. The dimension estimates given in (\cite{BG}, Proposition~4.3.3) show that $F^{\theta}(K)$ is placed in non-positive perverse degrees. Now (\ref{eq_hyp_loc}) garantees that $F^{\theta}(K)$ is placed in non-negative perverse degrees.
\end{Prf}

\medskip

 Let $w_0$ (resp., $w_0^Q$) denote the longest element of the Weil group $W$ of $G$ (resp., $W_Q$ of $Q$). 
 
\begin{Cor} 
\label{Cor1}
i) Let $\lambda\in\Lambda^+$ and $\theta$ be the image of $\lambda$ in $\Lambda_{G,P}$. Then
$\cA_{Q,\lambda}$ (resp., $\cA_{Q,-w_0^Q(\lambda)}$) appears with multiplicity one in $F^{\theta}(\cA_{\lambda})$ (resp., in $F^{-\theta}(\cA_{\lambda})$).

\smallskip\noindent
ii) The functor $F: \Sph(\wt\Gr_{G,x})\to \Sph(\Gr_{Q,x})$ given by $F=\mathop{\oplus}\limits_{\theta\in \Lambda_{G,P}} F^{\theta}$ is exact and faithful.
\end{Cor}
\begin{Prf}
i) Note that  $S^{\theta}_P\cap \Gr_G^{\lambda}$ is open in $\Gr_G^{\lambda}$. Moreover, $\Gr_Q^{\theta}\cap \Gr_G^{\lambda}=\Gr_Q^{\lambda}$. Since $P/Q$ is affine, $\Gr_Q\hook{} S_P$ is a closed immersion. So, $\Gr_Q^{\theta}\cap \Gr_G^{\lambda}\hook{} S^{\theta}_P\cap \Gr_G^{\lambda}$ is a smooth closed subscheme. It follows that 
$(\tilde\gr^{\theta}_P)^!(\tilde \gs^{\theta}_P)^*\cA_{\lambda}$ is a shifted constant sheaf over $\Gr_Q^{\lambda}$. The first assertion follows.

 For the second, note that $\Gr_Q^{-\theta}\cap \Gr_G^{\lambda}=\Gr_Q^{-w_0^Q(\lambda)}$, and the map 
$$
\gt^{-\theta}_P: S^{-\theta}_P\cap \Gr_G^{\lambda}\to\Gr_Q^{-\theta}
$$ 
is an isomorphism over the $Q(\cO)$-orbit $\Gr_Q^{-w_0^Q(\lambda)}$. 

\smallskip\noindent
ii) Since $F$ is exact, to show faithfulness, it suffices to prove that $F$ does not annihilate a nonzero object. To this end, it suffices to show that $F(\cA_{\lambda})\ne 0$ for any dominant coweight $\lambda$, which follows from i).
\end{Prf}

\bigskip\noindent
8.5 {\scshape Example: explicit calculation\  }
Let $\alpha\in\Lambda^+$ denote the coroot of $\Sp_{2n}$ corresponding to the maximal root $\check{\alpha}_{max}$ of $\Sp_{2n}$. So, $\alpha$ is the highest weight of the standard representation of the Langlands dual group $\SO_{2n+1}$ of $\Sp_{2n}$. For this subsection take $G$ to be that of 8.1 for $M_0=\cO^n\oplus\Omega_{\cO}^n$. 
Following (\cite{BD}, Sect.~4.5.12) the closure $\ov{\Gr}_G^{\alpha}$ of $\Gr_G^{\alpha}$ in $\Gr_G$ is described as follows. 

 The $G(k)$-orbit $V$ in $\Gr_G$ passing through $\alpha(t)G(\cO)$ is identified with the projective space $V\,\iso\, \PP^{2n-1}$, and $\Gr_G^{\alpha}$ is the total space of the line bundle $\cO(2)$ over $V$. 
 
  Let $V=\PP^{2n-1}\hook{} \PP^{n(2n+1)-1}$ be the Veronese map. Write $x_1,\ldots,x_{2n}$ for the homogeneous coordinates in $\PP^{2n-1}$ and $t_{ij}$ with $1\le i\le j\le 2n$ for the homogeneous coordinates in $\PP^{n(2n+1)-1}$. Then the inclusion is given by $t_{ij}=x_ix_j$. Its image is the subscheme defined by homogeneous equations
 \begin{equation}
 \label{eq_Veronese}
 t_{ij}t_{kl}=t_{ik}t_{jl}
 \end{equation}
for all $i,j,k,l$ whenever  this makes sense. 

 One may identify the Lie algebra of $\Sp_{2n}$ with $\AA^{n(2n+1)}$ in such a way that the set
 $Z$ of elements $\Sp_{2n}$-conjugate to a multiple of the maximal root becomes the subscheme $Z\subset \AA^{n(2n+1)}=\Spec k[t_{ij}]$ given by equations (\ref{eq_Veronese}). Let $A\in Z$ denote the origin of this cone. Let $\bar Z\subset \PP^{n(2n+1)}$ be the projective closure of $Z$. Then $\ov{\Gr}_G^{\alpha}=\bar Z$ and $\Gr_G^{\alpha}=\bar Z-A$.  
 
  The projection $\pi: \bar Z-A\to V$ is an affine fibration on which $\cO(2)$ acts transitively and freely (and the corresponding torsor is trivial). So, $\pi^*$ yields a diagram of isomorphisms
$$
\begin{array}{ccccc}
\Cl(V) & \iso & \Cl(\bar Z-A) & \iso &\Cl(\bar Z)\\
\downarrow && \downarrow \\
\Pic(V) & \iso & \Pic(\bar Z-A) & \iso & \ZZ,
\end{array}
$$ 
where for a variety $S$ we denote by $\Cl(S)$ the Weil divisors class group.  

 Write $(t_{ij},w)$ for the homogeneous coordinates in $\PP^{n(2n+1)}$. Let the subscheme  $V\subset \bar Z$ be given by $w=0$, it is a section of $\pi$. We have $Z=\bar Z-V$. 
 
  The image in $\Cl(V)$ of the hyperplane section of $\PP^{n(2n+1)-1}$ is 2. It follows that the image of $V$ in  $\Cl(\bar Z)$ is 2 and $\Cl(Z)\,\iso\,\ZZ/2\ZZ$. 
 
 Let $L\subset Z$ denote the preimage under $\pi$ of the subscheme of $V$ given by $x_1=0$. 
Denote again by $L$ the corresponding Weil divisor on $\bar Z$. Then $L$ is not locally principal in $\cO_{Z,A}$. Indeed, let $\gp\subset \cO_{Z,A}$ denote the ideal corresponding to $L$ and $\gm_{Z,A}\subset\cO_{Z,A}$ the maximal ideal. Then
$t_{ij}$ ($1\le i\le j\le n$) form a base in the cotangent space $\gm_{Z,A}/\gm_{Z,A}^2$, and 
the elements $t_{1j}\in\gp$ ($1\le j\le n$) are linearly independent in $\gm_{Z,A}/\gm_{Z,A}^2$. So, $\Pic Z=0$, and $\cO_{\bar Z}(V)$ generates $\Pic(\bar Z)$. The image of $\cO_{\bar Z}(V)$ under the composition
$$
\Pic(\bar Z)\hook{} \Cl(\bar Z)\iso \Cl(\bar Z-A)\,\iso\,\Pic(\bar Z-A)\,\iso\,\ZZ
$$
is 2. In other words, $\cO_{\bar Z-A}(L)$ does not extend to a line bundle on $\bar Z$.

 The line bundle $\cL\mid_{\ov{\Gr}_G^{\alpha}}$ identifies with $\cO_{\PP^{n(2n+1)}}(1)\mid_{\bar Z}$. Let $\tilde Z\to\bar Z$ denote the $\mu_2$-gerbe of square roots of this bundle. We see that 
this gerbe is nontrivial, though trivial over $\bar Z-A$. 
 
  Set $Y=\AA^{2n}=\Spec k[x_i]$. Let $\tau: Y\to Z$ be the map given by $t_{ij}=x_ix_j$. Clearly, $Y-\tau^{-1}(A)\to Z-A$ is a $S_2$-Galois covering. 

 For a coweight $\lambda$ of $Q$ denote by $\cA_{Q,\lambda}$ the intersection cohomology sheaf of the $Q(\cO)$-orbit on $\Gr_Q$ passing through $\lambda(t)Q(\cO)$. 

\begin{Pp} 
\label{Pp_explicit_calculation}
1) The sheaf $\cA_{\alpha}$ is the extension by zero from $\bar Z-A$.\\
2) We have $F^0(\cA_{\alpha})=0$. For $\theta\in\Lambda_{G,P}$ such that $\<\theta,\check{\omega}_n\>=1$ we have $F^{\theta}(\cA_{\alpha})\,\iso\,\cA_{Q,\alpha}$ and
$F^{-\theta}(\cA_{\alpha})\,\iso\,\cA_{Q,-\alpha}$.
\end{Pp}
\begin{Prf}
1) Note that $\cO_{Z-A}(L)$ generates the group $\Pic(Z-A)\,\iso\,\Cl(Z-A)\,\iso\,\Cl(Z)\,\iso\,\ZZ/2\ZZ$. 
The gerbe $\tilde Z$ is obtained by gluing together
trivial gerbes $Z\times B(\mu_2)$ and $(\bar Z-A)\times B(\mu_2)$ over
$Z-A$. The gluing data is an automorphism of the gerbe $(Z-A)\times
B(\mu_2)$ which can be described as follows.

 An $S$-point of $(Z-A)\times B(\mu_2)$ is a line bundle $\cB$ on $S$ together
with $\cB^2\iso \cO_S$ and a map $S\to (Z-A)$. Our automorphism sends
this point to the same map $S\to (Z-A)$ and replaces $\cB$ by $\cB$
tensored with the restriction of $\cO_{Z-A}(L)$ to $S$. 

 We have the $\mu_2$-torsor over $Z-A$ consisting of those sections of
 $\cO_{Z-A}(L)$ whose square is 1. This is exactly the Galois covering
 $Y-\tau^{-1}(A)\to Z-A$. 

 Let $W$ denote the nontrivial rank one local system on $B(\mu_2)$ corresponding to the covering $\Spec k\to B(\mu_2)$. If we identify our gerbe over $Z$ with $Z\times B(\mu_2)$ then
 over that locus $\cA_{\alpha}$ becomes the exteriour product $N\boxtimes W$,
 where $N$ is the nontrivial local system on $Z-A$ extended by zero to
 $A$ and  corresponding to the covering $Y-\tau^{-1}(A)\to Z-A$.
 
 \medskip\noindent
 2) Considering $\Gr_Q^0$ as a subscheme of $\Gr_G$, one checks that $\Gr_Q^0\cap \ov{\Gr}_G^{\alpha}$ is the point scheme $1\in\Gr_{G}$. Consider the $*$-restriction 
 $N\mid_{Z\cap L}$. Since the $!$-fibre at $A$ of $N\mid_{Z\cap L}$ vanishes, we get $F^0(\cA_{\alpha})=0$.
 
   Let $\theta\in\Lambda_{G,P}$ be such that $\<\theta,\check{\omega}_n\>=1$. Recall the map $\pi: \bar Z-A\to V$. We have 
$$
\Gr_G^{\alpha}\cap S^{\theta}_P=\pi^{-1}(V_0),
$$ 
where $V_0\subset V=\PP(M_0(x)/M_0)$ is the complement to $\PP(L_0(x)/L_0)$. In other words, $\Gr_G^{\alpha}\cap S^{\theta}_P\subset \Gr_G^{\alpha}$ is the open subscheme given by the condition that the line $(M+M_0)/M_0$ is not contained in $L_0(x)/L_0$. 
Further, $\Gr^{\alpha}_G\cap \Gr_Q^{\theta}=\Gr_Q^{\alpha}$. The isomorphism  $F^{\theta}(\cA_{\alpha})\,\iso\,\cA_{Q,\alpha}$ follows.

 We have $\Gr_G^{\alpha}\cap S^{-\theta}_P=\Gr_Q^{-\alpha}$. This yields the last isomorphism.
\end{Prf}

\medskip

\begin{Rem} 
\label{Rem_estimate}
Let $\lambda\in\Lambda^+$ and $\theta\in\Lambda_{G,P}$. If $F^{\theta}(\cA_{\lambda})\ne 0$ then 
\begin{equation}
\label{ineq_1}
-\<\lambda,\check{\omega}_n\>\le \<\theta,\check{\omega}_n\>\le \<\lambda,\check{\omega}_n\>
\end{equation}
Indeed, if $S^{\theta}_P\cap\ov{\Gr}^{\lambda}_G\ne \emptyset$ then (\ref{ineq_1}) holds. 
  More generally, for a reductive group $G$ and its parabolic subgroup $P$ the condition $S^{\theta}_P\cap\ov{\Gr}^{\lambda}_G\ne \emptyset$ implies $ \<\lambda,w_0(\check{\lambda})\>
\le  \<\theta,\check{\lambda}\>\le \<\lambda,\check{\lambda}\>$ for any $\check{\lambda}\in\check{\Lambda}_{G,P}$ which is dominant for $G$.
\end{Rem}

\bigskip\noindent
8.6 {\scshape The functors $F^{\theta}_{X^d}$ } 

\smallskip\noindent
 Let $\Gr_{Q,X^d}$ denote the ind-scheme
classifying $(x_1,\ldots,x_d)\in X^d$ and  $L\in\Bun_n$ with trivialization $L\,\iso\, \cO^n\mid_{X-x_1\cup\ldots\cup x_d}$. Its connected components are indexed by $\Lambda_{G,P}$, the component $\Gr_{Q,X^d}^{\theta}$ is given by $\deg L=-\<\theta,\check{\omega}_n\>$. We have a natural map  $\Gr_{Q, X^d}\to \Gr_{G,X^d}$ sending the above point to $L\oplus (L^*\otimes\Omega)$ with the induced trivialization outside $x_i$. The composition
$$
(\Gr_{Q, X^d})_{red}\hook{}\Gr_{Q, X^d}\to \Gr_{G,X^d}
$$
is a closed immersion. 

For $\theta\in\Lambda_{G,P}$ denote by
$S_{P, X^d}^{\theta}$ the ind-scheme classifying collections: $(x_1,\ldots,x_d)\in X^d$, 
a $P$-torsor $\cF_P$ on $X$ with trivialization $\nu: \cF_P\,\iso\,\cF^0_P\mid_{X-x_1\cup\ldots\cup x_d}$ such that the induced $Q$-torsor $\cF_P\times_P Q$ lies in  $\Gr_{Q,X^d}^{\theta}$. Here $\cF^0_P$ is the $G$-torsor $\cF^0_G=\cO_X^n\oplus\Omega^n$ with $P$-structure corresponding to the lagrangian subbundle $\cO_X^n$. 

 Considering $\cF^0_{\bar P}$ as $\cF^0_G$ with $\bar P$-structure given by $\Omega^n$, one similarly defines the ind-scheme $S_{\bar P, X^d}^{\theta}$. 
As in 8.4, one defines a diagram
\begin{equation}
\label{diag_P_Q_X_d}
\begin{array}{ccc}
S_{P, X^d}^{\theta} & \toup{\gs^{\theta}_{P,X^d}} & \Gr_{G,X^d}\\
\uparrow\lefteqn{\scriptstyle\gr^{\theta}_{P,X^d}} && \uparrow\\
\Gr_{Q,X^d}^{\theta} & \to & S_{\bar P, X^d}^{\theta}
\end{array}
\end{equation}
Both $(S_{P, X^d}^{\theta})_{red}$ and $(S_{\bar P, X^d}^{\theta})_{red}$ are locally closed in $\Gr_{G,X^d}$, and their intersection is 
$(\Gr_{Q,X^d}^{\theta})_{red}$. 
  
  For a $k$-point $(x_1,\ldots,x_d)\in X^d$ with $\{x_1,\ldots,x_d\}=\{y_1,\ldots,y_s\}$ and $y_i$ pairwise distinct, the fibre of the diagram (\ref{diag_P_Q_X_d}) over $(x_1,\ldots,x_d)\in X^d$ is
$$
\begin{array}{ccc}
\mathop{\cup}\limits_{\theta_1+\ldots+\theta_s=\theta} (\prod_i S^{\theta_i}_P) & \to & \prod\limits_{i=1}^s \Gr_{G,y_i}\\
\uparrow && \uparrow\\
\mathop{\cup}\limits_{\theta_1+\ldots+\theta_s=\theta} (\prod_i \Gr^{\theta_i}_{Q,y_i}) & \to &
\mathop{\cup}\limits_{\theta_1+\ldots+\theta_s=\theta} (\prod_i S^{\theta_i}_{\bar P})
\end{array}
$$

 Similarly to $G_{X^d}$, one defines a group scheme $Q_{X^d}$ (resp., $P_{X^d}$) over $X^d$, it acts naturally on $\Gr_{Q,X^d}^{\theta}$ (resp., on $S^{\theta}_{P,X^d}$). Denote by $\Sph(\Gr_{Q,X^d}^{\theta})$ the category of $Q_{X^d}$-equivariant perverse sheaves on $\Gr_{Q,X^d}^{\theta}$. Let us define the functors
$$
F^{\theta}_{X^d}, F'^{\theta}_{X^d}: \Sph(\wt\Gr_{G,X^d})\to \D(\Gr^{\theta}_{Q,X^d})
$$ 
  
  Let $\tilde\gs^{\theta}_{P,X^d}: \tilde S^{\theta}_{P,X^d}\to \wt\Gr_{G,X^d}$ 
be the map obtained by the base change $\wt\Gr_{G,X^d}\to\Gr_{G,X^d}$ 
from (\ref{diag_P_Q_X_d}). As in Lemma~\ref{Lm_section_over_unip_orbits}, one defines a $P_{X^d}$-equivariant section
$i^{\theta}_{P,X^d}: S^{\theta}_{P,X^d}\to \tilde S^{\theta}_{P,X^d}$ of the gerbe $\tilde S^{\theta}_{P,X^d}\to S^{\theta}_{P,X^d}$. We have a $Q_{X^d}$-equivariant line bundle $_{\theta}\cL_{X^d}$ on $\Gr^{\theta}_{Q,X^d}$, whose fibre at 
$$
(L, L\,\iso\,\cO^n\mid_{X-x_1\cup\ldots\cup x_d})
$$ 
is 
$
\det\RG(X, \cO^n_X)\otimes \det\RG(X,L)^{-1}
$.  
As $\ZZ/2\ZZ$-graded, it is placed in degree $\flat(\theta):=\<\theta, \check{\omega}_n\> \!\!\mod 2$. The canonical $P_{X^d}$-equivariant $\ZZ/2\ZZ$-graded isomorphism 
$$
(\gs^{\theta}_{P,X^d})^*\cL\,\iso\, {_{\theta}\cL^{\otimes 2}_{X^d}\mid_{S^{\theta}_{P,X^d}}}
$$ 
yields $i^{\theta}_{P,X^d}$ via 3.1.2. Set
$$
F'^{\theta}_{X^d}(K)=(\gr^{\theta}_{P,X^d})^!(i^{\theta}_{P,X^d})^*(\tilde \gs^{\theta}_{P,X^d})^*K\;\;\;\; \mbox{and}\;\;\;\; 
F^{\theta}_{X^d}(K)=F'^{\theta}_{X^d}(K)\otimes\Qlb[1](\frac{1}{2})^{\otimes \<\theta, 2\check{\rho}-2\check{\rho}_Q\>}
$$
Note that
$$
F^{' \theta}_{X^d}(K)\,\iso\, (\gt^{\theta}_{P,X^d})_!(i^{\theta}_{P,X^d})^*(\tilde \gs^{\theta}_{P,X^d})^*K
$$
where $\gt^{\theta}_{P,X^d}: S^{\theta}_{P,X^d}\to \Gr^{\theta}_{Q,X^d}$ is the corresponding contraction map. 

 Remind the definition of the tensor category $\Sph(\Gr_{Q,x})^{\natural}$.
Equip $\Sph(\Gr_{Q,x})$ with the convolution product, associativity and commutativity constraints given by the fusion procedure, then $\Sph(\Gr_{Q,x})$ is a tensor category (\cite{BD}, 5.3.16). It has a canonical $\ZZ/2\ZZ$-grading compatible with the tensor structure, namely $\cA_{Q,\lambda}$ is even (resp., odd) if $\dim\Gr_Q^{\lambda}$ is even (resp., odd). The latter condition depends only on the connected component of $\Gr_{Q,x}$ containing $\Gr_{Q,x}^{\lambda}$. 
 
  Following (\cite{BD}, 5.3.21), we define $\Sph(\Gr_{Q,x})^{\natural}$ as the full subcategory of even objects in $\Sph(\Gr_{Q,x})\otimes\Vect^{\epsilon}$. We have an equivalence of monoidal categories $\Sph(\Gr_{Q,x})^{\natural}\to\Sph(\Gr_{Q,x})$ (i.e., it is compatible with tensor product and associativity constraints), and the commutativity constraints $A\otimes B\,\iso\, B\otimes A$ in these two categories differ by $(-1)^{\deg A\deg B}$. 

 Let $h^{\epsilon}: \Sph(\Gr_{Q,x})\to\Vect^{\epsilon}$  denote the global cohomology functor. Since $h^{\epsilon}$ is a tensor functor 
compatible with $\ZZ/2\ZZ$-gradings, it gives rise to a tensor functor 
$$
h: \Sph(\Gr_{Q,x})^{\natural}\to\Vect
$$ 
By \cite{MV}, $h$ is a fibre functor, and there is an isomorphism $\Aut^{\otimes} h\,\iso\, \check{Q}$, where $\check{Q}$ is the Langlands dual group to $Q$ (in \cite{BD}, 5.3.23 some properties of the action of $\check{Q}$ on $h$ are listed, which determine this isomorphism uniquely). Thus, $\Sph(\Gr_{Q,x})^{\natural}\,\iso\, \Rep(\check{Q})$ canonically as tensor categories. 
 
 Consider
\begin{equation}
\label{category_true_Sph}
\Sph'(\Gr_{Q,x}):=\mathop{\oplus}\limits_{\theta\in\Lambda_{G,P}} \Sph(\Gr^{\theta}_{Q,x})[\<\theta, 2\check{\rho}_Q-2\check{\rho}\>]\subset \D(\Gr_{Q,x})
\end{equation}
equipped with the convolution product, commutativity and associativity constraints given by the fusion procedure, so $\Sph'(\Gr_{Q,x})$ is a tensor category. 

\begin{Lm} There is a canonical equivalence of tensor categories $\Sph'(\Gr_Q)\,\iso\,\Sph(\Gr_Q)^{\natural}$. 
\end{Lm} 
\begin{Prf} 
Note that $2(\check{\rho}-\check{\rho}_Q)=(n+1)\check{\omega}_n\in \check{\Lambda}_{G,P}$. Consider the case of $n$ odd. In this case $\check{\rho}_Q\in\check{\Lambda}$, so all $Q(\cO)$-orbits on $\Gr_Q$ are even-dimensional and $\Sph(\Gr_G)\,\iso\,\Sph(\Gr_G)^{\natural}$. In this case the shifts in (\ref{category_true_Sph}) are even, and we are done.  
 
 Consider the case of $n$ even. The component $\Gr_{Q,x}^{\theta}$ is even iff $\<\theta, \check{\omega}_n\>$ is even. So, in (\ref{category_true_Sph}) the even (resp, odd) objects of $\Sph(\Gr_{Q,x})$ are shifted by even (resp., odd) cohomological degree. Our assertion follows. 
\end{Prf} 

\medskip

Equip $\Sph'(\Gr_{Q,x})$ with a new $\ZZ/2\ZZ$-grading such that $K\in\Sph'(\Gr^{\theta}_{Q,x})$ is placed in degree $\flat(\theta)$. This $\ZZ/2\ZZ$-grading is compatible with the tensor structure. 
Denote by $\Sph'(\Gr_{Q,x})^{\flat}$ the category of even objects in $\Sph'(\Gr_{Q,x})\otimes\Vect^{\epsilon}$, it is equipped with the induced $\ZZ/2\ZZ$-grading. 

 The proof of part ii) of the following proposition is postponed to Sect.~8.7.

 % Because of signs in the commutativity constraint, the compatibility of $F$ with the tensor structures is not clear, but we have the following. 

\begin{Pp} 
\label{Pp_tensor_functor_F}
 i) The functor $F': \Sph(\wt\Gr_{G,x})\to \Sph'(\Gr_{Q,x})^{\flat}$ given by $F'=\mathop{\oplus}\limits_{\theta\in\Lambda_{G,P}} F'^{\theta}$ is a tensor functor. \\
ii) There is a unique $\ZZ/2\ZZ$-grading on $\Sph(\wt\Gr_{G,x})$ such that $F'$ is compatible with $\ZZ/2\ZZ$-gradings.
\end{Pp}
\begin{Prf} i) Pick $F_1,F_2\in \Sph(\wt\Gr_G)$. Set $K_i=\tau^0 F_i$, 
$$
K=F^{\theta}_{X^2}(K_1\ast_X K_2)\;\;\;\;\mbox{and}\;\;\;\; K'=F'^{\theta}_{X^2}(K_1\ast_X K_2),
$$ 
where $\tau^0$ is given by (\ref{functor_tau_0}). By abuse of notation,
write also $\tau^0: \Sph(\Gr_Q)\to \Sph(\Gr_{Q,X})[-1]$ for the corresponding functor for $Q$.
 
\smallskip  
\Step 1  Recall that $U\subset X^2$ denotes the complement to the diagonal. Write $\wt\Gr_{G,X^2}(U)$ for the preimage of $U$ in $\wt\Gr_{G,X^2}$. We have a $\mu_2$-gerbe $q :(\wt\Gr_{G,X}\times\wt\Gr_{G,X})\mid_U \to \wt\Gr_{G,X^2}(U)$ (defined as the map $\tilde q_{G,X}$ in 8.3.1). The complex $q^*(K_1\ast_X K_2)$ identifies canonically with $(K_1\boxtimes K_2)\mid_U$. Denote by $i^{\theta}$ the composition
$$
S^{\theta}_{P,X^d}\,\toup{i^{\theta}_{P,X^d}}\, \tilde S^{\theta}_{P,X^d}\,
\toup{\tilde\gs^{\theta}_{P,X^d}}\, \wt\Gr_{G,X^d}
$$ 
For $\theta_1+\theta_2=\theta$ the following diagram is 2-commutative
$$
\begin{array}{ccc}
(\wt\Gr_{G,X}\times\wt\Gr_{G,X})\mid_U  & \toup{q} & \wt\Gr_{G,X^2}(U)\\
\uparrow\lefteqn{\scriptstyle i^{\theta_1}\times i^{\theta_2}} &&
\uparrow\lefteqn{\scriptstyle i^{\theta}}\\
(S^{\theta_1}_{P,X}\times S^{\theta_2}_{P,X})\mid_U & \hook{} & S^{\theta}_{P,X^2}(U),
\end{array}
$$
where the low horizontal arrow is the natural open immersion. However, the 2-morphism rending this diagram 2-commutative is well-defined only up to a sign, we normalize it as follows. 

 Write $_{\theta}\und{\cL}_{X^d}$ for the line bundle $_{\theta}\cL_{X^d}$ viewed as \select{ungraded}.
It suffices to pick an isomorphism 
$$
\epsilon^{\theta_1,\theta_2}: {_{\theta_1}\und{\cL}_X\boxtimes {_{\theta_2}\und{\cL}_X}}\,\iso\, (j^{\theta_1,\theta_2})^*{_{\theta}\und{\cL}_{X^2}},
$$
where $j^{\theta_1,\theta_2}:  (\Gr^{\theta_1}_{Q,X}\times\Gr^{\theta_2}_{Q,X})\mid_U\hook{} \Gr_{Q,X^2}^{\theta}(U)$ is the natural open immersion. 
The order of points in $X^2$ yields such $\epsilon^{\theta_1,\theta_2}$, and the usual Leibnitz rule is satisfied. 

 Namely, remind that $\sigma$ denotes the involution of $X^2$ permuting the points. For the diagram
$$
\begin{array}{ccc}
(\Gr^{\theta_1}_{Q,X}\times \Gr^{\theta_2}_{Q,X})\mid_U & \toup{j^{\theta_1,\theta_2}} & \Gr^{\theta}_{Q,X^2}(U)\\
\uparrow\lefteqn{\scriptstyle \sigma} &&  \uparrow\lefteqn{\scriptstyle \sigma}\\
(\Gr^{\theta_2}_{Q,X}\times \Gr^{\theta_1}_{Q,X})\mid_U & \toup{j^{\theta_2,\theta_1}} & \Gr^{\theta}_{Q,X^2}(U)
\end{array}
$$
the following diagram commutes
\begin{equation}
\label{diag_enfin_for_sings}
\begin{array}{ccccc}
\sigma^* (j^{\theta_1,\theta_2})^*{_{\theta}\und{\cL}_{X^2}} & \iso &  (j^{\theta_2,\theta_1})^*\sigma^*{_{\theta}\und{\cL}_{X^2}} & \iso & 
(j^{\theta_2,\theta_1})^*{_{\theta}\und{\cL}_{X^2}}\\
\uparrow\lefteqn{\scriptstyle \epsilon} && && \downarrow\lefteqn{\scriptstyle \sign}\\    
\sigma^*({_{\theta_1}\und{\cL}_X\boxtimes {_{\theta_2}\und{\cL}_X}}) & \iso & {_{\theta_2}\und{\cL}_X\boxtimes {_{\theta_1}\und{\cL}_X}} & \toup{\epsilon} & (j^{\theta_2,\theta_1})^*{_{\theta}\und{\cL}_{X^2}},
\end{array}
\end{equation}
where $\sign=(-1)^{\flat(\theta_1)\flat(\theta_2)}$, and the isomorphisms denoted by $\,\iso\,$ are the canonical ones.
    
\smallskip  
\Step 2  Note that $\Gr^{\theta}_{Q,X^2}(U)$ is the disjoint union of 
$(\Gr^{\theta_1}_{Q,X}\times\Gr^{\theta_2}_{Q,X})\mid_U$ for $\theta_1+\theta_2=\theta$.
Let us show that $K[2]$ is a perverse sheaf on $\Gr^{\theta}_{Q,X^2}$, the Goresky-MacPherson extension from $\Gr^{\theta}_{Q,X^2}(U)$. 
More precisely, we show that $\epsilon$ as above yield an isomorphism
\begin{equation}
\label{equation_epsilon_used}
(\tau^0 F'(F_1))\ast_X (\tau^0 F'(F_2))\,\iso\, F'_{X^2}(K_1\ast_X K_2)
\end{equation}

 Indeed, $\epsilon^{\theta_1,\theta_2}$ yields an isomorphism between the restriction of $K'$ to $(\Gr^{\theta_1}_{Q,X}\times\Gr^{\theta_2}_{Q,X})\mid_U$ and 
$$
\tau^0F'^{\theta_1}(F_1)\boxtimes \tau^0F'^{\theta_2}(F_2)
$$
So, $K[2]$ is a perverse sheaf over $\Gr^{\theta}_{Q,X^2}(U)$. 
Using (\ref{iso_on_diagonal}), we learn that the $*$-restriction of $K$ under the diagonal embedding $\Gr_{Q,X}\hook{}\Gr_{Q,X^2}$ identifies with $\tau^0 F^{\theta}(F_1\ast F_2)$, so it is placed in perverse degree 1. Now argue as in Proposition~\ref{Pp_exactness}, using the corresponding $\Gm$-action on $\wt\Gr_{G,X^2}$. By Proposition~\ref{Pp_hyperbolic_loc}, the $!$-restriction of $K$ under $\Gr_{Q,X}\hook{}\Gr_{Q,X^2}$ is placed in perverse degree 3. We have constructed the isomorphism (\ref{equation_epsilon_used}). 

 Restricting to the diagonal, it yields $\tau^0(F'(F_1)\ast F'(F_2))\,\iso\, \tau^0 F'(F_1\ast F_2)$.

\smallskip
\Step 3  Let us check the compatibility with the commutativity constraints. 
Using (\ref{diag_enfin_for_sings}) one shows that the diagram commutes
$$
\begin{array}{ccc}
\sigma^*(\tau^0F'(F_1)\ast_X \tau^0F'(F_2))  & \toup{\sigma^*\comp \,\epsilon} & \sigma^* F'_{X^2}(K_1\ast_X K_2)\\
\uparrow && \uparrow\\
\tau^0F'(F_2)\ast_X \tau^0F'(F_1) & \toup{\sign \comp \, \epsilon} & F'_{X^2}(K_2\ast_X K_1),
\end{array}
$$
where the vertical arrows are the canonical isomorphisms, and $\sign$ is that from Step 1. We are done. 
\end{Prf}
 
\bigskip\noindent
8.7 {\scshape The structure of $\Sph(\wt\Gr_G)$}

\medskip\noindent
Recall that $\Lambda_{G,P}$ is canonically identified with the lattice of characters of the center $Z(\check{Q})$ of the Langlands dual group $\check{Q}$ of $Q$. For a representation $V$ of $\SO_{2n+1}$ and $\theta\in\Lambda_{G,P}$ write $V_{\theta}$ for the direct summand of $V$ on which $Z(\check{Q})$ acts by $\theta$. 

 For $\lambda\in\Lambda^+$ write $V^{\lambda}$ for the irreducible representation of $\SO_{2n+1}$ of highest weight $\lambda$. Write $\omega_i\in\Lambda^+$ for the fundamental coweight of $G$ corresponding to the representation $\wedge^i V^{\alpha}$ of $\SO_{2n+1}$, $i=1,\ldots,n$.   
 Let $\Loc: \Rep(\check{Q})\to \Sph(\Gr_Q)^{\natural}$ denote the Satake equivalence, normalized to send an irreducible representation of $\check{Q}$ with highest weight $\mu$ to $\cA_{Q,\mu}$. 

\begin{Pp} 
\label{Pp_multiplicities}
Let $\lambda\in\Lambda^+$ and $\theta$ be the image of $\lambda$ in $\Lambda_{G,P}$. Then $F^{\theta}(\cA_{\lambda})\,\iso\,\Loc(V^{\lambda}_{\theta})$ canonically.
In particular, $F^{\theta}(\cA_{\omega_i})\,\iso\,\cA_{Q,\omega_i}$ for $\<\theta,\check{\omega}_n\>=i$. 
\end{Pp}
\begin{Prf} We could similarly define the functor $F^{\theta}:\Sph(\Gr_G)\to\Sph(\Gr_Q^{\theta})$.
Write $\cA_{\lambda, old}$ for the corresponding object of $\Sph(\Gr_G)$. We claim that $F^{\theta}(\cA_{\lambda})\,\iso\,F^{\theta}(\cA_{\lambda, old})$ canonically for our particular $\theta$.

 Indeed, $S^{\theta}_P\cap \ov{\Gr}_G^{\lambda}\hook{} \ov{\Gr}_G^{\lambda}$ is an open immersion, and the gerbe $\tilde S^{\theta}_P\to S^{\theta}_P$ is trivial. So, the $*$-restriction of $\cA_{\lambda}$ under $S^{\theta}_P\cap \ov{\Gr}_G^{\lambda}\to \wt\Gr_G$ is the Goresky-MacPherson extension from $S^{\theta}_P\cap \Gr_G^{\lambda}$. 
The assertion follows now from (Proposition~4.3.3 and Theorem 4.3.4,\cite{BG}).
\end{Prf}

\begin{Pp} 
\label{Pp_first_properties}
i) If $1\le i\le n$ then $\cA_{\omega_i}$ appears in $\cA_{\alpha}^{\otimes i}$.\\
ii) For $\lambda,\mu\in\Lambda$ the multiplicity of $\cA_{\lambda+\mu}$ in $\cA_{\lambda}\otimes \cA_{\mu}$ is one. 
\end{Pp}
\begin{Prf}
i) Let $\theta\in\Lambda_{G,P}$ be given by $\<\theta,\check{\omega}_n\>=i$. By Proposition~\ref{Pp_tensor_functor_F}, 
$F(\cA_{\alpha}^{\otimes i})\,\iso\, (\cA_{Q,\alpha}\oplus \cA_{Q,-\alpha})^{\otimes i}$. So, $F^{\theta}(\cA_{\alpha}^{\otimes i})\,\iso\, \cA_{Q, \alpha}^{\otimes i}$. Applying an appropriate symmetrazation functor (either invariants or anti-invariants), one gets a direct summand $\cV\subset \cA_{\alpha}^{\otimes i}$ such that $F^{\theta}(\cV)\,\iso\, \cA_{Q,\omega_i}$. 

 If $\cA_{\lambda}$ appears in  $\cV$ then $F^{\theta}(\cA_{\lambda})\subset F^{\theta}(\cV)$, because $F^{\theta}$ is exact. Besides, $\lambda\le i\alpha$ in the sense that $\Gr_G^{\lambda}\subset \ov{\Gr}_G^{i\alpha}$, so $\<\lambda,\check{\omega}_n\>\le i$. If  $\<\lambda,\check{\omega}_n\><i$ then $F^{\theta}(\cA_{\lambda})=0$ by Remark~\ref{Rem_estimate}. If $\<\lambda,\check{\omega}_n\>=i$ then, by Corollary~\ref{Cor1}, $\cA_{Q,\lambda}$ appears in $F^{\theta}(\cV)\,\iso\,\cA_{Q,\omega_i}$, so $\lambda=\omega_i$. 
The assertion follows. 

\smallskip\noindent
ii) Consider the convolution map $m:\ov{\Gr}_G^{\lambda}\tilde\times \ov{\Gr}_G^{\mu}\to \ov{\Gr}_G^{\lambda+\mu}$ as in Sect.~8.2. Its restriction to the open subscheme $\Gr_G^{\lambda}\tilde\times\Gr_G^{\mu}\to \Gr_G^{\lambda+\mu}$ is an isomorphism, as follows from (\cite{MV}, Lemma~4.3 and formula 3.6). We are done.
\end{Prf}

\medskip

\begin{Prf}\select{of Proposition~\ref{Pp_tensor_functor_F} ii)}

\smallskip\noindent
Call an object $K\in\Sph(\wt\Gr_G)$ \select{even} (resp., \select{odd}) if $F^{\theta}(K)=0$ unless $\flat(\theta)=0$ (resp., $\flat(\theta)=1$). 
Proposition~\ref{Pp_semisimplicity} combined with Proposition~\ref{Pp_first_properties} shows that $\cA_{\alpha}$ is a tensor generator of $\Sph(\wt\Gr_G)$. 
Since $\cA_{\alpha}$ is odd, we get a $\ZZ/2\ZZ$-grading on $\Sph(\wt\Gr_G)$ compatible with the tensor structure. Moreover, 
$F'$ is compatible with the gradings. The uniqueness of the $\ZZ/2\ZZ$-grading is clear, because $\cA_{\alpha}$ is irreducible. 
\end{Prf}

\begin{Def} 
\label{Def_Sph_flat}
Let $\Sph(\wt\Gr_{G,x})^{\flat}$ be the category of even objects in $\Sph(\wt\Gr_{G,x})\otimes\Vect^{\epsilon}$. 
\end{Def}

 By Proposition~\ref{Pp_tensor_functor_F}, we get a tensor functor $F': \Sph(\wt\Gr_{G,x})^{\flat}\to \Sph'(\Gr_{Q,x})$. 
Denote by $F^{\natural}$ the composition 
$$
\Sph(\wt\Gr_G)^{\flat}\toup{F'}\Sph'(\Gr_Q)\,\iso\,\Sph(\Gr_Q)^{\natural}
$$ 
Let $\tilde h: \Sph(\wt\Gr_G)^{\flat}\to\Vect$ denote the tensor functor $\tilde h=h\comp F^{\natural}$. 

\begin{Cor} 
There is an affine group scheme $\check{G}$ over $\Qlb$ such that $\Sph(\wt\Gr_G)^{\flat}$ and the category $\Rep(\check{G})$ of $\Qlb$-representations of $\check{G}$ are canonically equivalent as tensor categories. 
\end{Cor}
\begin{Prf}
By Corollary~\ref{Cor1}, for each nonzero $\lambda\in\Lambda^+$ the rank of $\tilde h(\cA_{\lambda})$ is at least 2.  By (Proposition~1.20, \cite{DM}), $\Sph(\wt\Gr_G)^{\flat}$ is a rigid abelian tensor category (cf. Definition~1.7, \select{loc.cit}) and $\tilde h:\Sph(\wt\Gr_G)^{\flat}\to\Vect$ is a fibre functor. Our assertion follows now from 
(Theorem~2.11, \select{loc.cit.}).
\end{Prf}

\medskip

  Write $W^{\lambda}$ for the representation of $\check{G}$ corresponding to $\cA_{\lambda}$, $\lambda\in\Lambda^+$.  The functor $F^{\natural}:\Sph(\wt\Gr_G)^{\flat}\to \Sph(\Gr_Q)^{\natural}$ yields a morphism
$\check{Q}\to \check{G}$. By Proposition~\ref{Pp_explicit_calculation}, $W^{\alpha}=U^{\alpha}\oplus U^{\alpha *}$, where $U^{\alpha}$ is the irreducible representation of $\check{Q}$ of highest weight $\alpha$. Since $W^{\alpha}$ is a faithful representation of $\check{Q}$, it follows that $\check{Q}\to \check{G}$ is an injection.  

 Since $W^{\alpha}$ is a tensor generator of $\Sph(\wt\Gr_G)^{\flat}$, $\check{G}$ is of finite type. We also get that $\check{G}\subset \SL(W^{\alpha})$. Indeed, the only object of rank one in $\Sph(\wt\Gr_G)^{\flat}$ is $\cA_0$, so $\check{G}$ acts trivially on $\det W^{\alpha}$. 

 Let $\cS\in\Rep(\check{G})$ be such that the strictly full subcategory of $\Rep(\check{G})$, whose objects are isomorphic to subobjects of $\oplus_{i=1}^m \cS$, is stable under the tensor structure. Then $\check{Q}$ acts trivially on $F^{\natural}(\cS)$, because $\check{Q}$ is connected. 
If $\check{Q}$ acts trivially on some $F^{\natural}(\cA_{\lambda})$ then $\lambda=0$ by Proposition~\ref{Pp_multiplicities}. So, $\cS$ is a multiple of $\cA_0$. By (\cite{DM}, 2.22), this implies that $\check{G}$ is connected. Now by (\select{loc.cit.}, 2.23), $\check{G}$ is reductive. 

 The above $\ZZ/2\ZZ$-grading on $\Sph(\wt\Gr_G)^{\flat}$ gives rise to a group homomorphism $\mu_2\to \check{G}$. 
 
\begin{Lm} 
\label{Lm_multiplicities_version_6}
For $i=1,\ldots,n$ the multiplicity of $W^{\omega_i}$ in $\wedge^i W^{\alpha}$ is one. If $W^{\lambda}$ appears in $\wedge^i W^{\alpha}$ and $\lambda\ne \omega_i$ then $\<\lambda, \check{\omega}_n\><i$.
\end{Lm}
\begin{Prf}
Let $\theta\in\Lambda_{G,P}$ be given by $\<\theta,\check{\omega}_n\>=i$. The direct summand of $\wedge^i W^{\alpha}=\wedge^i(U^{\alpha}\oplus U^{\alpha *})$, on which $Z(\check{Q})$ acts by $\theta$ is $\wedge^i U^{\alpha}$. It follows that $F^{\theta}(\wedge^i\cA_{\alpha})=\cA_{Q,\omega_i}$, where we denoted by $\wedge^i\cA_{\alpha}$ the object of $\Sph(\wt\Gr_G)^{\flat}$ corresponding to  $\wedge^i W^{\alpha}$. 

 If $W^{\lambda}$ appears in  $\wedge^i W^{\alpha}$ then $F^{\theta}(\cA_{\lambda})\subset F^{\theta}(\wedge^i\cA_{\alpha})$, because $F^{\theta}$ is exact. Besides, $\lambda\le i\alpha$ in the sense that $\Gr_G^{\lambda}\subset \ov{\Gr}_G^{i\alpha}$, so $\<\lambda,\check{\omega}_n\>\le i$. If  $\<\lambda,\check{\omega}_n\><i$ then $F^{\theta}(\cA_{\lambda})=0$ by Remark~\ref{Rem_estimate}. If $\<\lambda,\check{\omega}_n\>=i$ then, by Corollary~\ref{Cor1}, $\cA_{Q,\lambda}$ appears in $F^{\theta}(\wedge^i\cA_{\alpha})=\cA_{Q,\omega_i}$, so $\lambda=\omega_i$. 
The assertion follows. 
\end{Prf}

\medskip\noindent
\begin{Prf}\select{of Theorem~\ref{Th_Satake}}

\Step 1 Let us show that $\cA_{\alpha}\ast\cA_{\alpha}\,\iso\, \cA_{2\alpha}\oplus\cA_{\omega_2}\oplus \cA_0$ for $n\ge 2$ and $\cA_{\alpha}\ast\cA_{\alpha}\,\iso\, \cA_{2\alpha}\oplus \cA_0$ for $n=1$. 
 Indeed, by Proposition~\ref{Pp_first_properties}, $\cA_{2\alpha}\oplus\cA_{\omega_2}$ appears in $\cA_{\alpha}\ast\cA_{\alpha}$.  Let $\theta\in\Lambda_{G,P}$ be given by $\<\theta,\check{\omega}_n\>=2$. By Proposition~\ref{Pp_multiplicities}, $F^{\theta}(\cA_{2\alpha})\,\iso\,\cA_{Q,2\alpha}$ and $F^{\theta}(\cA_{\omega_2})\,\iso\,\cA_{Q,\omega_2}$. We have 
$$
F^{\theta}(\cA_{\alpha}\ast\cA_{\alpha})\,\iso\,\Loc((W^{\alpha}\otimes W^{\alpha})_{\theta})\,\iso\,\Loc(U^{\alpha}\otimes U^{\alpha})\,\iso\,\cA_{Q,2\alpha}\oplus \cA_{Q,\omega_2}
$$ 
So, 
$\cA_{\alpha}\ast\cA_{\alpha}\,\iso\, \cA_{2\alpha}\oplus\cA_{\omega_2}\oplus K$ 
for some $K\in\Sph(\wt\Gr_G)$ such that $F^{\theta'}(K)=0$ unless $\<\theta',\check{\omega}_n\><2$. Since $\cA_{\alpha}$ is odd, $\cA_{\alpha}\ast\cA_{\alpha}$ is even, so $K$ is multiple of $\cA_0$. The disired assertion follows now from $\Hom(\cA_0, \cA_{\alpha}\ast\cA_{\alpha})\,\iso\,\Hom(\cA_{\alpha}, \cA_{\alpha})\,\iso\,\Qlb$. 

\smallskip\noindent
\Step 2  Let us show that $\cA_0$ appears in $\wedge^2\cA_{\alpha}$. Assume the contrary, that is, $\cA_0$ appears in $\Sym^2 \cA_{\alpha}$. Then $n\ge 2$ and $\check{G}\subset \SO(W^{\alpha})$ for the symmetric form $\Sym^2 W^{\alpha}\to U^{\alpha}\otimes U^{\alpha *}\to\Qlb$. 

 Let $\check{U}$ (resp., $\check{U}^-$) denote the unipotent radical of the Siegel parabolic $\check{P}\subset\SO(W^{\alpha})$ (resp., $\check{P}^-\subset\SO(W^{\alpha})$) preserving the isotropic subspace $U^{\alpha}\subset W^{\alpha}$ (resp., $U^{\alpha *}\subset W^{\alpha}$). The Lie algebra $\Lie\check{G}$ is a $\check{Q}$-subrepresentation of 
$$
\so(W^{\alpha})=\gl(U^{\alpha})\oplus \Lie(\check{U})\oplus\Lie(\check{U}^-)
$$

 Since $\Lie \check{U}$ and $\Lie\check{U}^-$ are irreducible $\check{Q}$-modules, $\check{G}$ coincides with one of the groups $\check{Q},\check{P},\check{P}^-, \SO(W^{\alpha})$. Since $\check{G}$ is reductive, it is either $\check{Q}$ or $\SO(W^{\alpha})$. Since $W^{\alpha}$ is not irreducible as a representation of $\check{Q}$, $\check{G}\ne\check{Q}$, hence $\check{G}=\SO(W^{\alpha})$. 

 Now Lemma~\ref{Lm_multiplicities_version_6} shows that $\wedge^n W^{\alpha}\,\iso\, W^{\omega_n}\oplus W^{\lambda}$ for some $\lambda\in\Lambda^+$ with $\<\lambda,\check{\omega}_n\><n$. Let $\tilde U$ denote the kernel of the contraction map $\wedge^{n-1} U^{\alpha}\otimes U^{\alpha *}\to \wedge^{n-2} U^{\alpha}$, this is an irreducible $\check{Q}$-module. By the representation theory for $\SO_{2n}$, we have 
\begin{itemize} 
\item $\tilde U\subset W^{\lambda}\subset \wedge^n(U^{\alpha}\oplus U^{\alpha *})$ as $\check{Q}$-modules; 
\item if a weight $\theta$ of $Z(\check{Q})$ appears in $W^{\lambda}$ then $\<\theta,\check{\omega}_n\>\le n-2$;
\item for $\<\theta, \check{\omega}_n\>=n-2$ 
the direct summand of $W^{\lambda}$ on which $Z(\check{Q})$ acts by $\theta$ is $\tilde U$.
\end{itemize}

Let $\theta$ be the image of $\lambda$ in $\Lambda_{G,P}$, we get
$F^{\theta}(\cA_{\lambda})\,\iso\, \tilde U$. By Corollary~\ref{Cor1}, $\cA_{Q,\lambda}\,\iso\, \tilde U$. However, the highest weight of $\tilde U$ does not lie in $\Lambda_+$. This contradiction yields our statement.

\smallskip\noindent
\Step 3 We know already that $\check{G}\subset \Sp(W^{\alpha})$ for the form $\wedge^2 W^{\alpha}\to U^{\alpha}\otimes U^{\alpha *}\to\Qlb$. Let $\check{P}\subset\Sp(W^{\alpha})$ (resp., $\check{P}^-\subset\Sp(W^{\alpha})$) denote the Siegel parabolic preserving the lagrangian subspace $U^{\alpha}\subset W^{\alpha}$ (resp., $U^{\alpha *}\subset W^{\alpha}$). As in Step 2, one shows that $\check{G}$ coincides with one of the groups $\check{Q},\check{P},\check{P}^-, \Sp(W^{\alpha})$. Since $\check{G}$ is reductive, it is either $\check{Q}$ or $\Sp(W^{\alpha})$. The $\check{Q}$-representation $W^{\alpha}$ is not irreducible, so $\check{G}=\Sp(W^{\alpha})$.
\end{Prf}
 
\bigskip\bigskip
\centerline{\scshape 9. Hecke operators}

\bigskip\noindent
9.1 According to A.3, inside of $\D(\wt\Bun_G)$ we have the full triangulated subcategories $\D_{\pm}(\wt\Bun_G)$. Let us define for each $K\in\Sph(\wt\Gr_G)$ a Hecke operator 
$\H(K, .): \D(\wt\Bun_G)\to \D(X\times\wt\Bun_G)$ sending $\D_{\pm}(\wt\Bun_G)$ to $\D_{\pm}(X\times\wt\Bun_G)$.
 
  Denote by $\cH_G$ the Hecke stack classifying $(\cF_G, \cF'_G, x\in X, \beta)$, where $\cF_G,\cF'_G$ are $G$-torsors on $X$, and $\beta:\cF_G\,\iso\,\cF'_G\mid_{X-x}$ is an isomorphism. We have the diagram 
$$
\Bun_G\getsup{p}\cH_G\toup{p'}\Bun_G,
$$
where $p$ (resp., $p'$) sends the above point to $\cF_G$ (resp., to $\cF'_G$). Let $\tilde\cH_G$ be the stack obtained from $\Bunt_G\times\Bunt_G$ by the base change $\cH_G\toup{p,p'}\Bun_G\times\Bun_G$. Denote by $\tilde p,\tilde p'$ the projections that fit into the diagram
$$
\begin{array}{ccccc}
\Bunt_G & \getsup{\tilde p} & \tilde\cH_G & \toup{\tilde p'} & \Bunt_G\\
\downarrow && \downarrow && \downarrow\\
\Bun_G & \getsup{p} & \cH_G & \toup{p'} & \Bun_G  
\end{array}
$$
Recall that the 'trivial' $G$-torsor $\cF_G^0$ on $X$ is given by $M_0=\cO_X^n\oplus\Omega^n$. Write $\Bun_{G,X}$ for the stack classifying triples $(\cF_G,x\in X,\nu)$, where $\cF_G\in\Bun_G$ and $\nu:\cF_G\,\iso\,\cF^0_G\mid_{D_x}$ is a trivialization over the formal disk $D_x$ at $x\in X$. Then $\Bun_{G,X}$ is a $G_X$-torsor over $X\times\Bun_G$. Set $\Bunt_{G,X}=\Bunt_G\times_{\Bun_G} \Bun_{G,X}$. 

 Denote by $\gamma$ (resp., $\gamma'$) the isomorphism $\Bun_{G,X}\times_{G_X} \Gr_{G,X}\,\iso\,\cH_G$ such that the projection to the first term corresponds to $p$ (resp., to $p'$). Recall the line bundle $\cA$ on $\Bun_G$ (cf.~3.2).  We have canonically 
$$
 \gamma'^* p^*\cA\,\iso\, \cA\tboxtimes \cL^{-1}
$$
This yields a $G_X$-torsor $\Bunt_{G,X}\times \wt\Gr_{G,X}\to \tilde\cH_G$ extending the $G_X$-torsor 
$$
\Bun_{G,X}\times\Gr_{G,X}\to  \Bun_{G,X}\times_{G_X} \Gr_{G,X}\toup{\gamma'}\cH_G
$$ 
So, for $\cS\in\Sph(\wt\Gr_{G,X})$ and $\cT\in \D(\Bunt_G)$ we can form their twisted tensor product 
$
\cT\tboxtimes \cS\in\D(\tilde\cH_G)
$. 
Set 
$$
\H(\cS,\cT)=(\supp\times\tilde p)_!(\cT\tboxtimes \cS),
$$ 
where $\supp:\tilde\cH_G\to X$ is the projection. In a similar way, for any $\cS\in\Sph(\wt\Gr_{G,X^d})$ one defines the functor $\H(\cS,.): \D(\Bunt_G)\to\D(X^d\times \Bunt_G)$. 
   
  Recall the functor $\glob: \Sph(\wt\Gr_G)\to\Sph(\wt\Gr_{G,X})$ (cf. 8.3.1).
For $K\in\Sph(\wt\Gr_G)$ set $\H(K,\cT)=\H(\glob(K),\cT)$.  

 The Hecke functors commute with Verdier duality $\DD\H(K,\cT)\,\iso\,\H(\DD K,\DD\cT)$, because $\Gr_G$ is ind-proper. Besides, they are compatible with the convolution product on $\Sph(\wt\Gr_G)$, namely, for $\cS_1,\cS_2\in\Sph(\wt\Gr_{G,X})$ we have canonically $\H(\cS_2, \H(\cS_1,\cT))\,\iso\,\H(\cS_1\ast_X \cS_2,\cT)$.    
 
\smallskip 
  
  The \select{geometric Langlands program for the metaplectic group} would be a trial to understand the action of $\Sph(\wt\Gr_G)^{\flat}$ on $\D_-(\Bunt_G)$, that is, to look for automorphic sheaves or, more generally, for a `spectral decomposition' of $\D_-(\Bunt_G)$ under this action. 
 
\smallskip 
  
   Recall that the metaplectic representation is automorphic. In the geometric setting this is reflected in the following Hecke property of $\Aut$. Set 
$$
\St=\Qlb[2n-1](\frac{2n-1}{2})\oplus \Qlb[2n-3](\frac{2n-3}{3})\oplus\ldots\oplus\Qlb[1-2n](\frac{1-2n}{2}),
$$
so $\St$ has cohomologies in odd degrees only and $\DD(\St)\,\iso\,\St$ as a complex over $\Spec k$.                         
\begin{Th} 
\label{Th_Hecke} 
Over $X\times\Bunt_G$ we have 
\begin{eqnarray*}
\H(\cA_{\alpha}, \Aut_g)\,\iso\, \St[1](\frac{1}{2})\boxtimes\Aut_s\\
\H(\cA_{\alpha},\Aut_s)\,\iso\, \St[1](\frac{1}{2})\boxtimes\Aut_g
\end{eqnarray*}
\end{Th}  

\noindent
9.2 \select{Proof of Theorem~\ref{Th_Hecke}}.

\medskip\noindent
Let $\cH_G^{\alpha}\subset\cH_G$ be the locally closed substack given by the condition that $\cF_G$ is in the position $\alpha$ with respect to $\cF'_G$ (or, equivalently, $\cF'_G$ is in the position $\alpha$ with respect to $\cF_G$). Set $\tilde\cH_G^{\alpha}=\cH_G^{\alpha}\times_{\cH_G} \tilde\cH_G$.

\begin{Lm} 
\label{Lm_kappa}
There exist isomorphisms 
$$
\kappa, \kappa': \tilde\cH^{\alpha}_G\; \iso\; (\Bunt_G\times_{\Bun_G}\cH^{\alpha}_G)\times B(\mu_2),
$$ 
where we used $p:\cH^{\alpha}_G\to\Bun_G$ (resp., $p': \cH^{\alpha}_G\to\Bun_G$) in the fibred product, and the projection to the first term corresponds to $\tilde p:\tilde\cH^{\alpha}_G\to \Bunt_G$ (resp., to $\tilde p':\tilde\cH^{\alpha}_G\to \Bunt_G$).
\end{Lm}
\begin{Prf} 
A point of $\tilde\cH^{\alpha}_G$ is given by $(\cF_G,\cF'_G,x\in X,\beta)\in\cH^{\alpha}_G$, two 1-dimensional vector spaces $\cB,\cB'$ with $\cB^2\,\iso\,\det\RG(X,M)$,
$\cB'^2\,\iso\,\det\RG(X,M')$. Here $M,M'$ are vector bundles on $X$ obtained from $\cF_G,\cF'_G$ via the standard representation of $G$.

 The symplectic form on $M$ induces a perfect pairing  $(M+M')/M\otimes(M+M')/M'\to \Omega(x)/\Omega\,\iso\, k$ between these 1-dimensional spaces. Further,
 $$
 \frac{\det\RG(X,M)}{\det\RG(X,M')}\,\iso\,\frac{(M+M')/M'}{(M+M')/M}\,\iso\,((M+M')/M')^{\otimes 2}
 $$ 
Instead of providing $\cB,\cB'$ we may provide $\cB,\cB_0$, where $\dim\cB_0=1$, with an isomorphism $\cB_0^2\,\iso\, k$, letting $\cB'=\cB\otimes((M+M')/M')^*\otimes\cB_0$. This defines $\kappa$. The datum of $\cB',\cB_0$ defines $\kappa'$.
\end{Prf}

\medskip  
  
   As above, let $W$ denote the nontrivial local system of rank one on $B(\mu_2)$ corresponding to the covering $\Spec k\to B(\mu_2)$. For the diagram
$$
X\times\Bunt_G  \getsup{\supp\times\tilde p}\; \tilde\cH^{\alpha}_G\;\toup{\tilde p'}\, \Bunt_G
$$
the Hecke operator writes $\H(\cA_{\alpha}, K)\,\iso\,(\supp\times\tilde p)_!(\tilde p'^*K\otimes \kappa^*W)[2n+1](\frac{2n+1}{2})$. 
  
\medskip\noindent
9.2.1 \select{Stratifications}  

\smallskip\noindent
Let $(x, M)$ be a $k$-point of $X\times {_i\Bun_G}$. Denote by $Y$ the fibre of $\supp\times p:\cH^{\alpha}_G\to X\times \Bun_G$ over $(x,M)$. So, $Y$ can be identified with the variety $\bar Z-A$ of Sect.~8.5. Let $Y_k$ denote the preimage of $_k\Bun_G$ under $Y\hook{}\cH^{\alpha}_G\toup{p'}\Bun_G$. We are going to describe the stratification of $Y$ by the subschemes $Y_k$.

  Recall that $M\in\Bun_{2n}$ with symplectic form $\wedge^2 M\to\Omega$ and $\dim\H^0(M)=i$ (for brevity, in this subsection we omit the argument $X$ in the cohomology groups). For a $k$-point $M'$ of $Y$ we get
$$
\begin{array}{ccccc}
& M &  \subset & M+M' & \subset M(x)\\
& \cup && \cup\\
M(-x)\subset & M\cap M' & \subset & M'
\end{array}
$$
and $\dim(M+M')/M=1$, $\dim (M\cap M')/M(-x))=2n-1$. Actually, $(M\cap M')/M(-x)$ is the orthogonal complement to $(M+M')/M$ for the perfect  pairing 
$$
M(x)/M\otimes M/M(-x)\to \Omega(x)/\Omega\,\iso\, k
$$ 
induced by the form on $M$. Let $\pi:Y\to V=\PP(M(x)/M)$ be the map sending $M'$ to the line $M+M'/M$. Let $N$ be the image of $\H^0(M)\to M/M(-x)$. Set $j=\dim N$, so $\dim\H^0(M(-x))=i-j$. Since $M\,\iso\,M^*\otimes\Omega$,
$$
\H^0(M(-x))\,\iso\,\H^1(M(x))^*\;\;\; \mbox{and}\;\;\; \H^1(M(-x))\,\iso\,\H^0(M(x))^*
$$ 
The long exact sequence 
$$
0\to\H^0(M)\to \H^0(M(x))\to M(x)/M\to\H^1(M)\to\H^1(M(x))\to 0
$$
shows that $\dim\H^0(M(x))=i+2n-j$, because $\dim\H^1(M(x))=i-j$. We have 
$$
\H^0(M\cap M')\,\iso\,\H^1(M+M')^*\;\;\;\mbox{and}\;\;\; \H^1(M\cap M')\,\iso\, \H^0(M+M')^*,
$$
because $(M+M')^*\otimes\Omega\,\iso\, M\cap M'$. Note that $\chi(M\cap M')=-1$ and $\chi(M+M')=1$. 

We distinguish three cases
\begin{itemize}
\item[0)] $j=0$. So, $\H^0(M(-x))=\H^0(M)$ is $i$-dimensional and $\dim\H^0(M(x))=2n$. 
Then $\H^0(M(-x))\,\iso\,\H^0(M\cap M')$ is of dimension $i$, and $\dim\H^0(M+M')=i+1$. 
Clearly, for $M+M'\in\PP(M(x)/M)$ fixed we get a 1-dimensional subspace in $(M+M')/(M\cap M')$ generated by $\H^0(M+M')$. So, for $M+M'\in V$ fixed
there is a unique $M'$ with $\dim\H^0(M')=i+1$ and for the other $M'$ we have $\dim\H^0(M')=i$. 

 Thus, $\pi:Y\to V$ has a section $V\to Y$, which is the closed stratum $Y_{i+1}$. Its complement is the open stratum $Y_i$. 

\item[1)] $0<j<2n$. View $V$ as the space of hyperplanes in $M/M(-x)$. We get a nontrivial subspace $V'\subset V$ of hyperplanes that contain $N$. Distinguish two cases:

CASE 1a) $N\subset (M\cap M')/M(-x)$ then $\H^0(M\cap M')=\H^0(M)$ is of dimension $i$, so $\dim\H^0(M+M')=i+1$. In the fibre of $\pi:Y\to V$ over $M+M'/M$ we get a distinguished point corresponding to the subspace of $(M+M')/(M\cap M')$ generated by $\H^0(M+M')$. This point lies in $_{i+1}\Bun_G$, and the complement lies in $_i\Bun_G$. 

\smallskip
CASE 1b) $N\nsubseteq (M\cap M')/M(-x)$. Then $N\cap (M\cap M')$ is of dimension $j-1$. So,  $\dim\H^0(M\cap M')=i-1$ and $\dim\H^0(M+M')=i$. Since $M'\ne M$, we get $M'\in{_{i-1}\Bun_G}$.

\smallskip

 So, $Y$ has three nonempty strata in case 1). The map $\pi: \pi^{-1}(V')\to V'$ has a section, which is the closed stratum $Y_{i+1}\,\iso\, V'$. The complement to this section is the middle stratum $Y_i=\pi^{-1}(V')-V'$, and the open stratum is $Y_{i-1}=\pi^{-1}(V-V')$.
 
\item[2)] $j=2n$. Then $\H^0(M)=\H^0(M(x))$ is $i$-dimensional, so $\dim\H^0(M+M')=i$ and
$\dim\H^0(M\cap M')=i-1$. The image of $\H^0(M)\to (M+M')/(M\cap M')$ is 1-dimensional and equals $M/(M\cap M')$. So, $\dim\H^0(M')=i-1$, because $M'\ne M$. 

 In this case $Y=Y_{i-1}$.
\end{itemize}
  
\noindent
Fix in addition a vector space $\cB$ together with $\cB^2\,\iso\,\det\RG(X,M)$. 

\begin{Pp}  
\label{Pp_16_fibres}
Let $K$ denote the fibre of $\H(\cA_{\alpha},\Aut_g)$ (resp., of $\H(\cA_{\alpha},\Aut_s)$) at $(x,M,\cB)\in X\times {_i\Bunt_G}$. Then 
$K=0$ unless $i$ is odd (resp., even). If $i$ is odd (resp., even) then 
we have noncanonically $K\,\iso\,\St[1+d_G-i]$. 
\end{Pp}
\begin{Prf}
g) Consider the case where $K$ is the fibre of $\H(\cA_{\alpha},\Aut_g)$. Assume $i$ even, so only the stratum $Y_i$ of $Y$ contributes to $K$. 

 If $j=0$ then $Y_i$ is a $\Gm$-torsor over $V$, and the restriction of $\Aut_g$ to a fibre of $\pi:Y_i\to V$ is a nontrivial local system of order two, so $K=0$ in this case. 
If $j=2n$ then $K=0$ because $Y=Y_{i-1}$. If $0<j<2n$ then $Y_i$ is a $\Gm$-torsor over $V'$, and the restriction of $\Aut_g$ to a fibre of $\pi: Y_i\to V'$ is a nontrivial local system of order two, so $K=0$.

 Now let $i$ be odd, so only the strata $Y_{i-1}$ and $Y_{i+1}$ contribute to $K$.
 
If $j=0$ then the restriction of $\Aut_g$ to $Y_{i+1}$ is isomorphic to $\Qlb[d_G-i-1]$ by Theorem~\ref{Th_1}, because $Y_{i+1}\,\iso\,\PP^{2n-1}$ is simply-connected. Our assertion follows then from 
$$
\St\,\iso\, \RG(\PP^{2n-1},\Qlb)[2n-1](\frac{2n-1}{2})
$$
If $j=2n$ then the restriction of $\Aut_g$ to $Y_{i-1}$ is isomorphic to $\Qlb[d_G-i+1]$, because $Y_{i-1}$ is simply-connected. So, $K\,\iso\,\St[1+d_G-i]$. If $0<j<2n$ then the restriction of $\Aut_g$ to $Y_{i+1}$ identifies with $\Qlb[d_G-i-1]$, because $Y_{i+1}\,\iso\, V'$ is simply-connected. The contribution of $Y_{i+1}$ to $K$ is 
$$
\RG(V',\Qlb)[d_G-i+2n]
$$ 
The restriction of $\Aut_g$ to
$Y_{i-1}$ is $\Qlb[d_G-i+1]$, because any rank one local system of order two on $\pi^{-1}(V-V')$ is trivial. So, the contribution of $Y_{i-1}$ to $K$ is
$\RG_c(V-V', \Qlb)[d_G-i+2n]$. The distinguished triangle
$$
\RG_c(V-V', \Qlb)[d_G-i+2n]\to K\to \RG(V',\Qlb)[d_G-i+2n]
$$
yields the desired isomorphism.

\medskip\noindent
s) In the case where $K$ is the fibre of $\H(\cA_{\alpha},\Aut_s)$, the argument is similar. 
\end{Prf} 
  
\bigskip\noindent
9.2.2   For $k,r\ge 0$ denote by $_{k,r}\cH^{\alpha}_G$ the preimage of $_k\Bun_G\times {_r\Bun_G}$ under $p\times p':\cH^{\alpha}_G\to \Bun_G\times\Bun_G$.
Similarly, define the stack $_{k,r}\tilde\cH^{\alpha}_G$ by the cartesian square  
$$
\begin{array}{ccc}
 _{k,r}\tilde\cH^{\alpha}_G & \hook{} & \tilde\cH^{\alpha}_G\\ 
  \downarrow && \downarrow\lefteqn{\scriptstyle \tilde p\times\tilde p'}\\
_k\Bunt_G\times {_r\Bunt_G} & \hook{} & \Bunt_G\times\Bunt_G
\end{array}
$$  

  The two $S_2$-coverings over $ _{k,r}\tilde\cH^{\alpha}_G$ obtained from $_k\rho: \Cov(_k\Bunt_G)\to {_k\Bunt_G}$ and from $_r\rho: \Cov(_r\Bunt_G)\to {_r\Bunt_G}$ 
are canonically isomorphic, namely Lemma~\ref{Lm_kappa} implies the following.
\begin{Lm} 
\label{Lm_can_diag}
There is a canonical commutative diagram, where both squares are cartesian
$$
\begin{array}{ccccc}
_k\Bun_G & \gets & _{k,r}\cH^{\alpha}_G\times B(\mu_2) & \to & _r\Bun_G\\
\downarrow\lefteqn{\scriptstyle _k\rho} && \downarrow && \downarrow\lefteqn{\scriptstyle _r\rho}\\
_k\Bunt_G & \getsup{\tilde p} & _{k,r}\tilde\cH^{\alpha}_G & \toup{\tilde p'} & _r\Bunt_G 
\end{array}
\eqno{\square}
$$ 
\end{Lm}

\smallskip

  Let $\cU\subset X\times {_1\Bun_G}$ be the open substack given by $\H^0(X, M(-x))=0$. 
As in Lemma~\ref{Lm_nonempty}, one shows that $\cU$ is non empty. In general, $\cU\ne X\times {_1\Bun_G}$. Let $\tilde\cU$ be the preimage of $\cU$ in $X\times{_1\Bunt_G}$.  
  
\begin{Pp} 
\label{Pp_17}
The first isomorphism of Theorem~\ref{Th_Hecke} holds over $\tilde\cU$, the second holds over $X\times{_0\Bunt_G}$.
\end{Pp} 
\begin{Prf} 
g) Let $Y(\cU)$ be the preimage of $\cU$ under $\supp\times p: \cH^{\alpha}_G\to X\times\Bun_G$. Write $Y_k(\cU)$ for the preimage of $_k\Bun_G$ under $Y(\cU)\hook{}\cH^{\alpha}_G\toup{p'}\Bun_G$. Then $Y_0(\cU)\to\cU$ (resp., $Y_2(\cU)\to\cU$) is a fibration with fibre isomorphic to $\PP^{2n-2}$ (resp., to $\AA^{2n}$). 

 Let $Y_k(\tilde \cU)$ be the preimage of $Y_k(\cU)$ in $\tilde\cH^{\alpha}_G$. For $k=0,2$ the restriction of the local system $\tilde p'^*(_k\Aut)\otimes\kappa^*W$ descends under  $Y_k(\tilde \cU)\to\tilde \cU$ to a local system, which is canonically identified, by Lemma~\ref{Lm_can_diag}, with $\Qlb\boxtimes {_1\Aut}$. 

By Proposition~\ref{Pp_16_fibres}, $\H(\cA_{\alpha}, \Aut_g)$ vanishes over $X\times{_0\Bunt_G}$, and we denote by $K$ the restriction of this complex to $\tilde \cU$. By decomposition theorem, $K$ is a direct sum of (shifted) irreducible perverse sheaves. We get an isomorphism
$$
K\,\iso\, {_1\Aut}[d_G-2n+1](\frac{d_G-2n+1}{2})\oplus {_1\Aut}[d_G+2n-1](\frac{d_G+2n-1}{2})\otimes\RG(\PP^{2n-2},\Qlb)
$$  
The first assertion follows.

\medskip\noindent
s) Set $\cV=X\times{_0\Bun_G}$. Let $K$ be the restriction of $\H(\cA_{\alpha},\Aut_s)$ to $\tilde\cV=X\times{_0\Bunt_G}$. Let $Y(\cV)$ be the preimage of $\cV$ under $\supp\times p:\cH^{\alpha}_G\to X\times\Bun_G$. Write $Y_k(\cV)$ for the preimage of $_k\Bun_G$ under $Y(\cV)\hook{}\cH^{\alpha}_G\toup{p'}\Bun_G$. Then $Y_1(\cV)\to\cV$ is a fibration with fibre isomorphic to $\PP^{2n-1}$.   
  
   Let $Y_1(\tilde\cV)$ be the preimage of $Y_1(\cV)$ in $\tilde\cH^{\alpha}_G$. By Lemma~\ref{Lm_can_diag}, the $*$-restriction of $\tilde p'^*(_1\Aut)\otimes\kappa^*W$ descends under $Y_1(\tilde\cV)\to\tilde\cV$ to a local system canonically identified with $\Qlb\boxtimes{_0\Aut}$. 
By decomposition theorem, one gets an isomorphism
$$
K\,\iso\, {_0\Aut}\otimes\RG(\PP^{2n-1}, \Qlb)[d_G+2n](\frac{d_G+2n}{2})
$$  
We are done.
\end{Prf}  

\bigskip

  By decomposition theorem, $\H(\cA_{\alpha},\Aut)$ is a direct sum of (shifted) irreducible perverse sheaves. Proposition~\ref{Pp_17} implies that  $\St[1](\frac{1}{2})\boxtimes\Aut$ appears in it as a direct summand. But according to Proposition~\ref{Pp_16_fibres}, all the fibres of $\H(\cA_{\alpha},\Aut)$ and of $\St[1](\frac{1}{2})\boxtimes\Aut$ are isomorphic. This concludes the proof of Theorem~\ref{Th_Hecke}.

\bigskip\bigskip
\centerline{\scshape Appendix A.}

\bigskip\noindent
A.1 For the convenience of the reader we collect here some generalities on group actions.

  Let $f:\cY\to\cZ$ be a morphism of stacks, $G\to\cZ$ be a group scheme over $\cZ$. 
Write $m_G$ for the product in $G$ and $1_G:\cZ\to G$ for the unit section.
Following \cite{Br}, an action of $G$ on $\cY$ over $\cZ$ is the data of a 1-morphism $m:G\times_{\cZ}\cY\to\cY$ over $\cZ$, a 2-morphism $\mu: m\comp(m_G\times\id)\Rightarrow m\comp(\id\times m)$ making the following diagram 2-commutative
$$
\begin{array}{ccc}
G\times_{\cZ} G\times_{\cZ}\cY & \toup{m_G\times\id} & G\times_{\cZ}\cY\\
\downarrow\lefteqn{\scriptstyle \id\times m} && \downarrow\lefteqn{\scriptstyle m}\\
G\times_{\cZ}\cY & \toup{m} & \cY,
\end{array}
$$
and a 2-morphism $\epsilon: m\comp(1_G\times \id_{\cY})\to \id_{\cY}$. They should satify two axioms: an associativity condition with respect to any 3 objects in G (cf. diagram (6.1.3) in loc.cit.); $\epsilon$ is compatible with $\mu$ (cf. diagrams (6.1.4) in loc.cit.). The fact that $m$ is a $\cZ$-morphism means that the diagram
$$
\begin{array}{ccc}
G\times_{\cZ}\cY &\toup{m} & \cY\\
\downarrow\lefteqn{\scriptstyle \pr_2} && \downarrow\lefteqn{\scriptstyle f}\\
\cY & \toup{f} & \cZ
\end{array}
$$
is 2-commutative.

  For a line bundle $L$ on $\cY$ we have a notion of \select{$G$-equivariant structure} on $L$ (cf. \cite{La}, Definition~2.8). A version of this notion for an $\ell$-adic complex is as follows. 
\begin{Def} A \select{$G$-equivariant structure} on $K\in \D(\cY)$ is an isomorphism
$\lambda: m^*K\,\iso\,\pr_2^*K$ such that two diagrams commute
$$
\begin{array}{ccc}
(m_G\times\id_{\cY})^*m^*K & \toup{\lambda} & (m_G\times\id_{\cY})^*\pr_2^*K\\
\downarrow\lefteqn{\scriptstyle \mu} && \downarrow\lefteqn{\scriptstyle \lambda}\\
(\id_G\times m)^*m^*K & \toup{\lambda} & (\id_G\times m)^*\pr_2^*K=\pr_{23}^*m^*K
\end{array}
$$
and
$$
\begin{array}{ccc}
(1_G\times\id_{\cY})^*m^*K \\
\downarrow\lefteqn{\scriptstyle \lambda} & \searrow\lefteqn{\scriptstyle \epsilon}\\
(1_G\times\id_{\cY})^*\pr_2^*K & = & K,
\end{array}
$$

\medskip\noindent
where $\pr_2: G\times_{\cZ}\cY\to\cY$ and $\pr_{23}: G\times_{\cZ}G\times_{\cZ}\cY\to
G\times_{\cZ}\cY$ are the projections.
\end{Def}

\noindent
A.2 Let $f:\cY\to\cZ$ be a representable morphism of algebraic stacks, $G\to\cZ$ be a group scheme over $\cZ$ acting on $\cY$ over $\cZ$. By definition, $\cY$ is a \select{$G$-torsor over $\cZ$} if, locally in flat topology of $\cZ$, $\cY$ is isomorphic to $G$ over $\cZ$ as a $G$-scheme. 

  Assume that $\cZ$ is locally of finite type. The notion of a perverse sheaf localizes in the smooth topology, so we have a notion of a perverse sheaf on $\cZ$. For the same reason, if $G\to \cZ$ is of finite type and smooth of relative dimension $d$ then the functor $K\mapsto f^*K[d]$ is an equivalence of the category of perverse sheaves $P(\cZ)$ on $\cZ$ with the category of $G$-equivariant perverse sheaves $P_{G}(\cY)$ on $\cY$. 

\bigskip
\noindent
A.3  Let $\cA$ be a line bundle on a scheme $S$. Let $\tilde S\to S$ denote the $\mu_2$-gerbe of square roots of $\cA$ (cf. 3.3.1). Since $\mu_2$ acts on $\tilde S$ by 2-automorphisms of the identity $\id:\tilde S\to\tilde S$, $\mu_2$ acts on any $K\in\D(\tilde S)$.
Write $\pi:\tilde S\to S$ for the structural morphism. 

\begin{Lm} 
\label{Lm_Lemma_descent_on_gerb}
1) The functor $\pi^*$ is an equivalence of the category of perverse sheaves on $S$ with the category of those perverse sheaves on $\tilde S$ on which $\mu_2$ acts trivially.

\smallskip\noindent
2) The functor $\pi^*:D(S)\to D(\tilde S)$ is fully faithful, its image $D_+(S)$ is a full triangulated subcategory of $D(\tilde S)$. 

\smallskip\noindent
3) For $K\in D(\tilde S)$ the following are equivalent
\begin{itemize}
\item[i)] $-1\in\mu_2$ acts as $-1$ on each cohomology sheaf of $K$ 
\item[ii)] $\pi_! K=0$  
\item[iii)] $\pi_* K=0$. 
\end{itemize}
Let $D_-(\tilde S)\subset D(\tilde S)$ be the full triangulated subcategory of objects satisfying these conditions.

\smallskip\noindent
4) For any $K_{\pm}\in D_{\pm}(\tilde S)$ we have $\Hom_{D(\tilde S)}(K_+, K_-)=0$
and $\Hom_{D(\tilde S)}(K_-, K_+)=0$.
For $K\in D(\tilde S)$ there exist $K_{\pm}\in D_{\pm}(\tilde S)$ such that $K\,\iso\, K_+\oplus K_-$. 
\end{Lm}
\begin{Prf} 1a) In the case $\cA=\cO_S$ consider the presentation $i: S\to B(S/\mu_2)$. The functor $i^*$ identifies the category of perverse sheaves on $B(S/\mu_2)$ with
the category of perverse sheaves on $S$ equipped with an action of the group $\mu_2(S)$.

\noindent
1b)  In general we have a carthesian square
$$
\begin{array}{ccc}
\tilde S & \toup{\pi}&  S\\ 
\uparrow\lefteqn{\scriptstyle h} && \uparrow\lefteqn{\scriptstyle\pi}\\
\tilde S\times B(\mu_2) & \toup{\pr} & \tilde S,
\end{array}
$$
where $h$ sends a $T$-point $(\cB,\cB_0,\; \cB^2\iso\cA\mid_T\; \cB_0^2\iso\cO_T)$ to $\cB\otimes\cB_0$ for any $S$-scheme $T$. 

 If $F$ is a perverse sheaf on $\tilde S$ on which $\mu_2$ acts trivially, then $\mu_2\times\mu_2$ acts trivially on $h^*F$. By 1a) we then get an isomorphism $h^*F\,\iso\pr^*F$ satisfying the usual cocycle condition. So, there is an isomorphism
$F\,\iso\, \pi^*H$ for some perverse sheaf $H$ on $S$.

\smallskip\noindent
2) The map $\pi$ is smooth of relative dimension zero, and $\pi_!\Qlb\iso\Qlb$. It follows formally that $\pi^*$ is fully faithful. 

\smallskip\noindent
3) The functors $\pi_!$ and $\pi_*$ are exact with respect to the usual t-structure. So, $\pi_!K=0$
iff $\pi_!(H^i(K))=0$ for each $i$. The latter is equivalent to requiring that $-1$ acts nontrivially on 
$H^i(K)$ for each $i$. Similarly for $\pi_*$. 

\smallskip\noindent
4) Given $K_-\in D_-(\tilde S)$ and $K_+\,\iso\,\pi^*L\in D_+(\tilde S)$ we have 
$$
\Hom(K_-, K_+)\,\iso\,\Hom(K_-, \pi^!L)\,\iso\,\Hom(\pi_!K_-, L)=0
$$ 
and 
$$
\Hom(K_+, K_-)\,\iso\,\Hom(\pi^*L; K_-)\,\iso\,
\Hom(L, \pi_*K_-)=0
$$

 We claim that for each $K\in D(\tilde S)$ the adjointness map $\pi_*\pi^*\pi_*K\to \pi_*K$ is an isomorphism. Since our derived categories are bounded, by devissage we may assume that $K$ is placed in cohomological dimension zero. Then $K\iso K_0\oplus K_1$, where $-1$ acts on $K_0$ (resp., on $K_1$) as 1 (resp., as -1). Clearly, $\pi^*\pi_*K_0\,\iso\, K_0$ and $\pi_* K_1=0$, so  $\pi_*\pi^*\pi_*K\,\iso\,\pi_*K$.
 
  For $K\in D(\tilde S)$ let $K_-$ be a cone of the adjointness map $\pi^*\pi_*K\to K$ then $\pi_*K_-=0$. The triangle $\pi^*\pi_*K\to K\to K_-$ splits, because $\Hom(K_-, \pi^*\pi_*K[1])=0$.
 \end{Prf}

\bigskip

 Let $G$ be an algebraic group acting on $S$, assume that $\cA$ is equipped with a $G$-equivariant structure. Then $G$ acts on $\tilde S$, and the projection
$\tilde S\to S$ is $G$-equivariant. 
 
 The stack $\tilde S$ is equipped with the universal line bundle $\cB_u$ together with $\cB^2_u\,\iso\,\cA\mid_{\tilde S}$. One checks that $\cB_u$ is $G$-equivariant.
 
 Let $G$ act on the trivial gerbe $S\times B(\mu_2)$ as the product of the action of $G$ on $S$ with the trivial action on $B(\mu_2)$. The following lemma is straightforward.

\begin{Lm} 
\label{Lm_A_equivariant}
Let $\cB$ be a $G$-equivariant line bundle on $S$ equipped with a $G$-equivariant isomorphism $\cB^2\,\iso\,\cA$. Then $\cB$ yields a $G$-equivariant
trivialization $\tilde S\;\iso\; S\times B(\mu_2)$. \QED
\end{Lm}

\smallskip\noindent
A.4 Let $S$ be a normal variety with a $\Gm$-action, $\cA$ be a $\Gm$-equivariant line bundle on $S$. Write $\tilde S\to S$ for the gerbe of square roots of $\cA$. Let $S_0\subset S$ be the variety of fixed points. For a connected component $C$ of $S_0$ set
$$
S^+(C)=\{s\in S\mid \lim_{t\to 0} ts\in C\} \;\; \mbox{and}
$$
$$
S^-=\{s\in S\mid \lim_{t\to\infty} ts\in C\}
$$
Let $S^+$ (resp., $S^-$) denote the disjoint union of $S^+(C)$ (resp., of $S^-(C)$) indexed by the connected components of $S_0$. Write $\tilde S^+$ (resp., $\tilde S^-, \tilde S_0$) for the restriction of the gerbe $\tilde S\to S$ to the corresponding scheme. Let $f^{\pm}: \tilde S_0\to \tilde S^{\pm}$ and $g^{\pm}: \tilde S^{\pm}\to \tilde S$ denote the corresponding (representable) maps.  Following \cite{B}, define \select{hyperbolic localization} functors $\D(\tilde S)\to \D(\tilde S_0)$ by
$$
K^{!*}=(f^+)^!(g^+)^*K, \;\;\;\;\;\;  K^{*!}=(f^-)^*(g^-)^!K
$$
The following  generalization of Theorem~1 from \select{loc.cit.} is straightforward.
 
\begin{Pp} 
\label{Pp_hyperbolic_loc}
There is a natural map $i_S: K^{*!}\to K^{!*}$ functorial in $K\in\D(\tilde S)$.  Assume that there is a covering of $S$ by open $\Gm$-invariant subschemes $U_i$ and $\Gm$-invariant trivializations $\xi_i: \cA\mid_{U_i}\,\iso\,\cO\mid_{U_i}$. Then for $\Gm$-equivariant $K\in\D(\tilde S)$ the map $i_S$ is an isomorphism. 
\end{Pp}
\begin{Prf} The map is constructed as in (\select{loc.cit.}, Sect.~2). Let $\tilde U_i$ denote the restriction of $\tilde S$ to $U_i$. It suffices to show the desired map is an isomorphism over $\tilde U_i$ for any perverse sheaf $K\in P(\tilde S)$.
The trivialization $\xi_i$ induces $\Gm$-equivariant section $U_i\to \tilde U_i$ of the gerbe $\tilde U_i\to U_i$. One concludes applying Theorem~1 from  \select{loc.cit.} for $K\mid_{U_i}$. 
\end{Prf}   

\medskip

 Assume in addition that there is a $\Gm$-equivariant section $S^+\to \tilde S^+$ of the gerbe
 $\tilde S^+\to S^+$. Let $h^+:S^+\to S_0$ be the map sending $s$ to $\lim_{t\to 0} ts$. Then for any $\Gm$-equivariant object $K\in\D(\tilde S)$ we have $K^{!*}\,\iso\, (h^+\times\id)_!(g^+)^*K$ canonically. Here $h^+\times\id: \tilde S^+\iso S^+\times B(\mu_2)\to S_0\times B(\mu_2)=\tilde S_0$. 

\bigskip\bigskip
\centerline{\scshape Appendix B. Weil representation and the sheaf $S_M$}

\bigskip\noindent
B.1  Let $k=\Fq$ be a finite field with $q$ odd. Let $M$ be a symplectic space over $k$ of dimension $2d$. The sheaf $S_M$ introduced in Sect.~4.4 has its origin in the Weil representation, this is what we are going to explain. 

  Consider the Heisenberg group $H(M)=M\oplus k$ with operation
$$
(m,a)(m',a')=(m+m',a+a'+\frac{1}{2}\<m,m'\>)
$$
Fix an additive character $\psi:k\to\Qlb^*$. There exists a unique up to isomorphism irreducible representation of $H(M)$ over $\Qlb$ with central character $\psi$. Let $(\rho, S_{\psi})$ be such representation. It yields an exact sequence
 \begin{equation}
 \label{sequence_meta}
 1\to\Qlb^*\to \tilde G\to G\to 1
 \end{equation}
 with $G=\Sp(M)$. Here 
 $$
 \tilde G=\{g,M[g])\mid g\in G, M[g]\in \Aut S_{\psi},\;\;\;
 \rho(gm,a)\comp M[g]=M[g]\comp\rho(m,a)\}
 $$

 Let $\cL(M)$ denote the variety of Lagrangian subspaces of $M$. For $L\in\cL(M)$ let $\chi_L: L\oplus k\to\Qlb^*$ send $(l,a)$ to $\psi(a)$. Set
$$
S_{L,\psi}=\Ind_{L\oplus k}^{H(M)}\chi_L=\{f: H(M)\to\Qlb\mid f(xh)=\chi_L(x)f(h)\;\,\mbox{for}\; x\in L\oplus k\}
$$
For each $L\in\cL(M)$ there is a pair $(v_L\in S_{\psi}, f_L\in S_{\psi}^*)$ which is $(L\oplus k,\chi_L)$-inivariant. Normalize it by $f_L(v_L)=1$, so any such pair is $(av_L, a^{-1}f_L)$ with $a\in\Qlb^*$. Specifying such pair is equivalent to specifying an isomorphism of $H(M)$-modules
$
S_{\psi}\,\iso\, S_{L,\psi}
$  
such that the image of $f_L$ becomes the evaluation at zero $f_{L,st}\in S_{L,\psi}^*$ (resp., $v_L$ becomes the function $v_{L,st}:H(M)\to\Qlb$ supported at $L\oplus k$ with $v_{L,st}(0)=1$).

 Let $P_L\subset G$ be the Seigel parabolic subgroup preserving $L$. Restricting  (\ref{sequence_meta}) we get an exact sequence
$$
1\to \Qlb^*\to\tilde P_L\to P_L\to 1
$$
The action of $\tilde P_L$ on  $\Qlb f_L$ yields a character $\tilde P_L\to\Qlb^*$ that splits this sequence (the group $\tilde P_L$ acts on $\Qlb v_L$ by the opposite character). 

  The \select{finite-dimensional theta-function} is $\theta_L: P_L\backslash \tilde G/P_L\to \Qlb$ given by $\theta_L(g)=f_L(gv_L)$, it does not depend on the choice of the pair $(v_L,f_L)$. 

\medskip\noindent
B.2 \  Let $L_1,L_2\in\cL(M)$. For $f\in S_{L_1,\psi}$ and $z\in L_2\oplus k$ the function
$f(zh)\chi_{L_2}^{-1}(z)$ depends only on the image of $z$ in $L_2$, so we may set
$$
(F_{L_1,L_2}(f))(h)=\int_{L_2} f(zh)\chi_{L_2}^{-1}(z)dz,
$$
where $dz$ is the Haar measure on $L_2$ such that the volume of a point is one. Then  
$F_{L_1,L_2}: S_{L_1,\psi}\,\iso\,S_{L_2,\psi}$ is an isomorphism of $H(M)$-modules.   

 One checks that $F_{L_2,L_1}\comp F_{L_1,L_2}\in \Aut(S_{L_1,\psi})$ is the multiplication by
$q^{d+\dim(L_1\cap L_2)}$. 

\begin{Def} For $L_1,L_2,V\in\cL(M)$ with $V\cap L_i=0$ define $\theta(L_1,L_2,V)\in\Qlb^*$ by
$$
F_{L_2,L_1}\comp F_{V,L_2}\comp F_{L_1,V}=\theta(L_1,L_2,V) 
$$  
\end{Def}

 We have $L_1=\{(bu+u)\mid u\in L_2\}$ for uniquely defined $b: L_2\to V$. The symplectic form on $M$ yields $L_2\,\iso\, V^*$, so $b$ becomes an element of $\Sym^2 V$. From definitions it follows that
\begin{equation}
\label{formula_theta_gausse}
\theta(L_1,L_2,V)=q^d \int_{V^*} \psi(\frac{1}{2}\<bv^*,v^*\>)dv^*,
\end{equation}
where $dv^*$ is the Haar measure on $V^*$ such that the volume of a point is one. 
 
  Denote by $\tilde\cY(k)$ the set of isomorphism classes of collections $L_1, L_2\in \cL(M)$, a one-dimensinal space $\cB$ together with $\cB^{\otimes 2}\,\iso\, (\det L_1)\otimes(\det L_2)$.   
So, $\tilde\cY(k)$ is a two-sheeted covering of the set $\cY(k)$ of $G$-orbits on $\cL(M)\times\cL(M)$. Remind that $\cY(k)$ contains $d+1$ element.  
 
  Given a triple $L_1,L_2,V\in\cL(M)$ with $L_i\cap V=0$, the form on $M$ yields isomorphisms $L_1\,\iso\, V^*\,\iso\, L_2$. So, $(L_1,L_2, \cB=\det V^*)$ is a point of $\tilde\cY(k)$.  Now
Proposition~\ref{Pp_linear_algebra_last} implies that $\theta(L_1,L_2,V)$ depends only on the image of $(L_1,L_2,V)$ in $\tilde\cY(k)$, so defining a function 
$$
\theta: \tilde\cY(k)\to\Qlb
$$
which is (up to a constant) the trace of Frobenius of the sheaf $S_M$. It is well-known that
for $(L_1,L_2,\cB)\in\tilde\cY(k)$ with $i=\dim(L_1\cap L_2)$ one gets 
$$
\theta(L_1,L_2,\cB)^2= \left(\frac{-1}{q}\right)^{d-i} q^{3d+i},
$$                     
where $\left(\frac{-1}{q}\right)=\left\{
\begin{array}{rl}
1, & \mbox{if}\; -1\in k^2\\  
 -1, & \mbox{otherwise}
 \end{array}
 \right.$
 
\bigskip\noindent
B.3 Remind that we fixed a square root $q^{\frac{1}{2}}$ of $q$ in $\Qlb$ (cf. 3.1). For $L_1,L_2\in \cL(M)$ set 
$$
\cF_{L_1,L_2}=q^{\frac{1}{2}(-d-\dim(L_1\cap L_2))} F_{L_1,L_2}
$$
The following is a version of the Maslov index (cf. \cite{LV}, appendix to chapter 1).

\begin{Def} For $L_1,L_2,L_3\in\cL(M)$ define $\gamma(L_1,L_2,L_3)\in\Qlb^*$ by
$$
\cF_{L_2,L_1}\comp \cF_{L_3,L_2}\comp \cF_{L_1,L_3}=\gamma(L_1,L_2,L_3)
$$
\end{Def} 
 
 Here are its immediate properties (cf. also \select{loc.cit.}).
\begin{Pp}
1) $\gamma(L_1,L_2,L_3)=\gamma(L_1, L_3, L_2)^{-1}=\gamma(L_2,L_1,L_3)^{-1}$.

\medskip\noindent
2) $\gamma(gL_1, gL_2, gL_3)=\gamma(L_1,L_2,L_3)$ for $g\in G$.

\medskip\noindent
3) If $L_1,L_2,L_3,L_4\in \cL(M)$ then 
$$
\gamma(L_1,L_2,L_3)\gamma(L_1,L_4,L_2)=\gamma(L_3,L_4,L_2)\gamma(L_1,L_4,L_3)
\eqno{\square}
$$
\end{Pp}

This implies that the function $(g_1,g_2)\mapsto \gamma(L, g_1L, g_1g_2L)$ is a 2-cocycle of $G$. This is the cocycle defining the extension (\ref{sequence_meta}). In our case of finite field $k$ this extension splits (\cite{MVW}, chapter 2, II.1).

\bigskip\noindent
{\bf Acknowledgements.}  I am very grateful to V. Drinfeld for his comments and suggestions about the first version of this paper. The idea of introducing the sheaf $S_M$ is due to him. I thank J.-L. Waldspurger for answering my questions and G. Laumon, D. Gaitsgory for stimulating discussions. I am also very grateful to the referee who has proposed so many improvements and corrections, in particular, the signs in the commutativity constraints (Proposition~\ref{Pp_tensor_functor_F}) have been clarified due to his comments.

\end{document}